\documentclass[a4paper,11pt]{amsart}
\usepackage{amssymb,amsthm,bbm,epic,eepic,graphics,amsbsy,url}
\usepackage{amssymb,amsthm}
\usepackage{xypic}
\usepackage[applemac]{inputenc}
\usepackage{longtable}
\usepackage{array}
\newcommand{\la}{\lambda}

\def\neqnear{\mbox{\ \ $\backslash\!\!\!\!\!\!\nearrow$}}

\DeclareMathOperator{\rk}{\mathrm{rk}}
\DeclareMathOperator{\tr}{\mathrm{tr}}
\DeclareMathOperator{\Irr}{\mathrm{Irr}}
\def\ad{\mathrm{ad}}

\def\Hom{\mathrm{Hom}}

\def\End{\mathrm{End}}
\def\Ker{\mathrm{Ker}\,}
\def\Res{\mathrm{Res}}
\def\Ind{\mathrm{Ind}}
\def\Aut{\mathrm{Aut}}

\def\C{\ensuremath{\mathbbm{C}}}
\def\Q{\mathbbm{Q}}
\def\Z{\mathbbm{Z}}

\def\R{\mathbbm{R}}
\def\bar{\overline}
\def\gl{\mathfrak{gl}}
\def\sl{\mathfrak{sl}}
\def\so{\mathfrak{so}}
\def\sp{\mathfrak{sp}}
\def\osp{\mathfrak{osp}}

\def\k{\mathbf{k}}

\def\om{\omega}

\def\eps{\epsilon}
\def\g{\mathfrak{g}}
\def\h{\mathfrak{h}}

\def\un{\mathbbm{1}}
\def\onto{\twoheadrightarrow}
\def\into{\hookrightarrow}
\def\ii{\mathrm{i}}
\def\Gal{\mathrm{Gal}}
\def\Irr{\mathrm{Irr}}
\def\GL{\mathrm{GL}}
\def\SL{\mathrm{SL}}
\def\SO{\mathrm{SO}}
\def\PSO{\mathrm{PSO}}
\def\OSP{\mathrm{OSP}}
\def\GOSP{\mathrm{GOSP}}
\def\Spin{\mathrm{Spin}}
\def\Im{\mathrm{Im}\,}
\def\bla{\boldsymbol{\lambda}}
\def\bmu{\boldsymbol{\mu}}
\def\aa{\mathtt{a}}
\def\bb{\mathtt{b}}
\def\bT{\mathbf{T}}
\def\LRef{\ensuremath{\Lambda\mathrm{Ref}}}
\def\QRef{\ensuremath{\mathrm{QRef}}}
\def\Ref{\ensuremath{\mathrm{Ref}}}

\def\Id{\mathrm{Id}}
\def\RR{\mathcal{R}}
\def\HH{\mathcal{H}}
\def\XX{\mathrm{X}}
\newtheorem{theo}{Theorem}[section]
\newtheorem{theointr}{Theorem}
\newtheorem{prop}[theo]{Proposition}
\newtheorem{defi}[theo]{Definition}
\newtheorem{lemma}[theo]{Lemma}
\newtheorem{remark}[theo]{Remark}
\newtheorem{cor}[theo]{Corollary}

\def\mod{\ \mathrm{mod}\ }

\makeatletter
\DeclareRobustCommand\widecheck[1]{{\mathpalette\@widecheck{#1}}}
\def\@widecheck#1#2{%
   \setbox\z@\hbox{\m@th$#1#2$}%
   \setbox\tw@\hbox{\m@th$#1%
      \widehat{%
         \vrule\@width\z@\@height\ht\z@
         \vrule\@height\z@\@width\wd\z@}$}%
   \dp\tw@-\ht\z@
   \@tempdima\ht\z@ \advance\@tempdima2\ht\tw@ \divide\@tempdima\thr@@
   \setbox\tw@\hbox{%
      \raise\@tempdima\hbox{\scalebox{1}[-1]{\lower\@tempdima\box\tw@}}}%
   {\ooalign{\box\tw@ \cr \box\z@}}}
\makeatother
\def\CW{W_0}

\title{{\bf Infinitesimal Hecke algebras II}}
\author{Ivan Marin}
\date{July 29, 2009} 
\begin{document}

\maketitle

\bigskip
\begin{center}
Institut de Math\'ematiques de Jussieu \\
Universit\'e Paris 7 \\
175 rue du Chevaleret \\
F-75013 Paris
\end{center}
\bigskip

\bigskip

\bigskip

\noindent {\bf Abstract.}
For $W$ a finite (2-)reflection group and $B$ its (generalized)
braid group, we determine the Zariski closure of the image of $B$
inside the corresponding Iwahori-Hecke algebra. The Lie algebra
of this closure is reductive and generated in the group algebra $\C W$ by the
reflections of $W$. We determine its decomposition in simple factors.
In case $W$ is a Coxeter group, we prove that the representations
involved are unitarizable when the parameters of the representations
have modulus 1 and are close to 1. We
consequently determine the topological closure in this case.
\medskip


\section{Introduction} 

Let $V$ be a finite-dimensional complex vector space. A \emph{reflection}
of $V$ is an element of $\GL(V)$ of order 2 whose fixed points
form an hyperplane of $V$, called the reflecting hyperplane of
the reflection. A central (finite) hyperplane arrangement
$\mathcal{A}$ in $V$ is called a reflection arrangement if there exists
a set $\mathcal{R}$ of reflections, whose set of reflection hyperplanes
is $\mathcal{A}$, and which generate a finite group $W$. If this is the case,
$\mathcal{R}$ and $W$ are uniquely determined by $\mathcal{A}$
(see \cite{arrrefl} prop. 2.1).

Topological structures associated to such a reflection arrangement include
the fundamental group $P$ of the complement $X = V \setminus \bigcup \mathcal{A}$
and $B = \pi_1(X/W)$. In case $\mathcal{A}$
is the complexification of a real arrangement, $W$ is a finite Coxeter
group and conversely every finite Coxeter group occur, through
their classical reflection representation. In this case, $B$
is the so-called Artin-Tits or generalized braid group associated
to $W$. If $W$ is of Coxeter type $A_{n-1}$ then $B$ is the classical braid group
on $n$ strands.

Another well-known algebraic structure associated to a central arrangement of
hyperplanes is its holonomy Lie algebra $\g$. It admits a nice
presentation discovered by Kohno, with one generator $t_H$ for
each hyperplane $H \in \mathcal{A}$. By definition, $W$-equivariant
representations of $\g$ yields by monodromy (linear) representations of $B$.
The representations induced by the natural projection $B \to W$
correspond to trivial representations of $\g$.

A seminal remark of Cherednik is that these
representations can be deformed into a 1-parameter family
by making $t_H$ act as $h s_H$, where $h$ is a formal parameter
and $s_H \in \mathcal{R}$ denotes the reflection associated to
$H \in \mathcal{A}$. It has been shown
in most cases that one gets in this way all representations of the
generic Iwahori-Hecke algebra associated to $W$, as monodromy
representations of $B$ over the field $K_0 = \C((h))$
of Laurent series.
The remarkable phenomenon here is that
the reflections fulfill in the group algebra $\C W$ the
defining Lie-theoretic relations of $\g$.

A natural object to consider is thus the \emph{Lie} subalgebra
of $\C W$ generated by $\mathcal{R}$, that we call
the \emph{infinitesimal Iwahori-Hecke algebra} $\mathcal{H}$ associated to
the reflection group $W$ --- or to $\mathcal{A}$, since $\mathcal{A}$
determines $(W,\mathcal{R})$. It turns out to be a reductive
Lie algebra, whose center is easily determined. The determination
of its simple factors enables to determine the Zariski hull
of $P$ and $B$ in the corresponding representations (see \S 6).

This task has been already accomplished for $W$ of dihedral type and for
the symmetric group (see \cite{HECKINF,LIETRANSP}).
More precisely, infinitesimal Iwahori-Hecke algebras were first
introduced in \cite{HECKINF} in case $W$ is a finite Coxeter group,
as a generalization of the Lie algebra
of transpositions introduced in \cite{THESE} and fully decomposed in \cite{LIETRANSP}.

The purpose of this work is to extend
this achievement to arbitrary reflection arrangement.

\subsection{Structure of infinitesimal Hecke algebras}

We let $\eps : W \to \{ \pm 1 \}$ denote the sign character
of $W$, defined by $\eps(s) = -1$ for $s \in \mathcal{R}$,
or equivalently as the restriction to $W$ of the determinant $\GL(V) \to
\C^{\times}$, and let $\k \subset \C$ denote the field of
definition of $W$, namely the algebraic number field
generated by the traces of elements of $W$ on $\GL(V)$.
It is well-known (see \cite{benard,bessis,MARINMICHEL})
that all the ordinary representations of $W$ are defined over $\k$.
The group $W$ is a Coxeter group iff $\k \subset \R$, and
a crystallographic group iff $\k = \Q$.

An intermediate Lie algebra between $\mathcal{H}$ and $\k W$
is the Lie subalgebra $\mathcal{L}_{\eps}(W)$ of $\k W$ spanned
by the $g - \eps(g) g^{-1}$ for $g \in W$, introduced and decomposed in \cite{galglie}.
We have $\mathcal{H} \subset
\mathcal{L}_{\eps}(W) \subset \k W$ since $s-\eps(s)s^{-1} = 2s$
for $s \in \mathcal{R}$.

A common feature of the three Lie algebras $\mathcal{H}$,
$\mathcal{L}_{\eps}(W)$ and $\k W$ is that they are reductive and
that we get irreducible
representations for them by restriction of the irreducible
representations of
$W$.
Moreover
\begin{enumerate}
\item the corresponding Lie ideals are simple
\item all simple Lie ideals are obtained this way
\end{enumerate}
This is shown in \cite{galglie} in the case of $\mathcal{L}_{\eps}(W)$ and will be a consequence
of our results in the case of $\mathcal{H}$. A consequence is that the simple Lie
ideals of these Lie algebras define an equivalence relation on the
irreducible representations of $W$, hence a partition of the
set $\Irr(W)$ of irreducible representations of $W$.

In case of $\mathcal{L}_{\eps}(W)$, this equivalence relation is
defined by $\rho \sim \rho^* \otimes \eps$, where $\rho^*$ denotes
the dual representation of $\rho$. The ideal $\mathcal{L}_0(\rho)$
of $\mathcal{L}_{\eps}(W)$ attached to the class of $\rho \in \Irr(W)$
is either of linear or orthosymplectic type, depending on whether
$\rho$ is isomorphic to $\rho^*\otimes \eps$ or not.
Denoting $V_{\rho}$ the vector space attached to $\rho \in \Irr(W)$, we
thus have
$$
\k W \simeq \bigoplus_{\rho \in \Irr(W)} \gl(V_{\rho}), \ \ \ 
\mathcal{L}_{\eps}(W) \simeq \bigoplus_{\rho \in \Irr(W)/\sim} \mathcal{L}_0(\rho), \ \ \ 
$$

In order to describe the decomposition of $\mathcal{H}$, we need to
define a set $\Ref(W)$ of what we call reflection representations (see def. \ref{defreflrep}),
$$
\begin{array}{l}
\mathrm{QRef} = \{ \eta \otimes \rho \ | \ \rho \in \mathrm{Ref}, \eta \in
\Hom(W,\{ \pm 1 \}) \} \\
\Lambda\mathrm{Ref} = \{ \eta \otimes \Lambda^k  \rho \ | \ \rho \in \mathrm{Ref}, \eta \in
\Hom(W,\{ \pm 1 \}), k \geq 0 \} \\
\end{array}
$$
and we define an equivalence relation $\approx$ on $\Irr(W)$ by
$$ \rho_1 \approx \rho_2 \Leftrightarrow
\rho_2 \in \{ \rho_1\otimes \eta, \rho_1^* \otimes \eta \otimes \eps \ | 
\ \eta \in \mathrm{X}(\rho_1) \}
$$
with
$$
\mathrm{X}(\rho) = \{ \eta \in \Hom(W,\{ \pm 1 \}) \ | \ 
\forall s \in \mathcal{R} \ \ \eta(s) = -1 \Rightarrow \rho(s) = \pm 1 \}.
$$

In case $\rho^* \otimes \eps \simeq \rho \otimes \eta$ for some $\eta
\in \mathrm{X}(\rho)$, then $\eps \otimes \eta$ embeds in either
$S^2 \rho^*$ or $\Lambda^2 \rho^*$. This defines an orthogonal
or symplectic Lie algebra $\osp(V_{\rho}) \subset \gl(V_{\rho})$
that contains $\rho(\mathcal{H}')$. We denote $\mathcal{L}(\rho) = \osp(V_{\rho})$
in this case, and $\mathcal{L}(\rho) = \gl(V_{\rho})$ otherwise.

Note that $\rho_1 \approx \rho_2$ iff $\rho_1 \sim \rho_2$,
and $\mathcal{L}(\rho) \simeq \mathcal{L}_0(\rho)$,
whenever there is a single conjugacy class of reflection,
notably when $W$ is a Coxeter group of type ADE

 We let
$\Irr'(W) = \Irr(W) \setminus \LRef(W)$. For a special situation
that occur only then $W$ has type $H_4$, we actually need a refinement
$\approx'$ of $\approx$.
\begin{theointr} The Lie algebra $\mathcal{H}$ is reductive, its
center has for dimension the number of reflection classes $\mathcal{R}/W$,
and
$$\mathcal{H} \simeq \k^{\mathcal{R}/W} \oplus \left( \bigoplus_{\rho \in \QRef/\approx'} \sl(V_{\rho}) \right)
\oplus \left(\bigoplus_{\rho \in \Irr'(W)/\approx'} \mathcal{L}(\rho) \right)
$$
\end{theointr}

The special situation occuring for $H_4$ is related simultaneaously
to triality and to the fact that $H_4$ is the only complex reflection
group which admits $\rho \in \Irr(W)$ of dimension 8 whose
symmetric square contains $\eps$. The two representations having
this property are equivalent under $\approx'$. For all the other groups,
$\approx'$ and the relation $\approx$ described above coincide.

We finally show that the \emph{real} Lie subalgebra
$\mathcal{H}^c$ of $\C W$ generated by the $\sqrt{-1} s, s \in \mathcal{R}$
is a compact form of $\mathcal{H}$.

\medskip

The sketch of proof of theorem 1 is as follows. In \S 2 we prove
that it is only needed to prove it with $\k$ replaced by $\C$, and that
it is true provided that $\mathcal{H}(\rho) = \mathcal{L}(\rho)$
for all $\rho \in \Irr(W) \setminus \LRef$. Since in general
$\mathcal{H}(\rho) \subset \mathcal{L}(\rho) \subset \sl_N$
where $N = \dim \rho$, we need several Lie-theoretic lemmas,
proved in \S 3, in order to deal with the situation $\h \subset \g \subset \sl_N$,
where $\h, \g$ are complex semisimple Lie algebras
acting irreductibly on $\C^N$, $\h$ being known and $\g$
having to be described. Then \S 4 establishes basic facts on the
representation theory of $W$ when $W$ belongs to the main general
series $G(e,e,n)$, and in particular determines $\LRef$
in this case.

The core of the proof is then \S 5, which sets up the induction
process, assuming that $W$ contains a proper reflection
subgroup $W_0$ for which the theorem is known. We prove
there this induction step when $W_0$ has a single class of
reflections, under some conditions on the branching rule
of the pair $(W,W_0)$ ; we use that to
reduce our problem to the case where $W$ has a single
class of reflections, and then apply these results to
the general series, with $W_0 = G(e,e,n)$ and $W =
G(e,e,n+1)$. The case of exceptional groups (except $H_4$) is then tackled, using a combination of the previous
results and ad-hoc (sometimes computer-aided) methods.

We then turn (\S 6) to the applications described below
and also deal with the special case of $H_4$. The
proofs of a few technical results needed in the proof of theorem 1
are postponed at the end of the paper (\S 7).

\subsection{Applications}

By the Cherednik monodromy construction, to each $\rho \in \Irr(W),
\rho : W \to \GL(V_{\rho})$
is associated a linear representation $R$ over $K_0 = \C((h))$
of the (generalized) braid group $B$ attached to $W$.
Letting $P$ denote the pure braid group $\Ker(B \onto W)$, a
first consequence of this work is the following.

\begin{theointr} Let $R : B \to GL_N(K_0)$ be the monodromy
representation associated to $\rho \in \Irr(W)$.
Then the Zariski closure of $R(P)$ is connected with Lie
algebra $\rho(\mathcal{H})$.
\end{theointr}

We actually make more precise the above statement by computing the
closure of $R(B)$. Depending on $\rho$, it coincides with the closure
of $R(P)$ or contains it as an index 2 subgroup.

Such representations of $B$ factorizes through a morphism
$B \to H(q)^{\times}$ where $H(q)$ denotes the so-called
Hecke algebra of $W$, with (transcendant) parameter $q = \exp(\ii \pi h)\in K_0$.
It is a natural quotient of the group algebra $K_0 B$.
In the cases where this Hecke algebra is known to be isomorphic
to the group algebra $K_0 W$ 
we more generally compute the Zariski closure of the image of
$P$  and $B$ inside $H(q)^{\times} \simeq (K_0 W)^{\times}$. 

We say that the groups for which this assumption is already known to hold
are \emph{tackled}. This includes all Coxeter groups by Tits theorem.

\begin{theointr} (see theorem \ref{theozarhecke})
Assume that $W$ is irreducible and tackled. The Zariski closure of the image of $P$ inside $H(q)^{\times}$
is connected and has Lie algebra $\mathcal{H} \otimes K_0$. It is contained
in the Zarisi closure of the image of $B$ as an index 2 subgroup.
\end{theointr}

We now assume that $W$ is a \emph{Coxeter group}, and consider the
representation $R : B \to \GL_N(K_0)$ of $H(q)$ attached to some
$\rho \in \Irr(W)$. Such a $\rho$ is a complexification
of a real representation $\rho_0$ of $W$. It is known that $R$ admits a matrix model over $\R[q,q^{-1}]$,
hence any $u \in \C^{\times}$ defines a specialized representation
$R_u : B \to \GL_N(\C)$. We prove the following (theorem \ref{theotophecke}).

\begin{theointr} For $u \in \C$ with $|u| = 1$
and $u$ close enough to 1, the representation
$R_u$ is unitarizable. If in addition $u$ is a transcendental number, then
the topological closure of $R_u(P)$ is a connected
compact Lie group, with Lie algebra $\rho(\mathcal{H}^c)$.
\end{theointr}

In particular, for the case $W = \mathfrak{S}_n$, we recover
the results of M. Freedman, M. Larsen, Z. Wang in \cite{FLW}
(for generic parameters). Recall that, in this case,
the unitarizability is known by Wenzl unitary matrix models.
Such models can be obtained in type $A$ by using a rational
Drinfeld associator, as in \cite{GT}. 

\subsection{Generalizations}

In case $W$ has several conjugacy classes of reflections,
the Hecke algebra can be defined with one parameter $q_c$ for
each class. In the generic case, these parameters
are distinct indeterminates. Such representations
are obtained by monodromy construction, with $q_c = \exp(\ii
\pi \la_c h)$, the $\la_c$ being complex numbers linearly independant
over $\Q$ (then the $q_c \in K_0$ are algebraically independant over $\C$).
In particular $\la_c \neq 0$. Then the relevant Lie
algebra, being spanned by the $\la_c s$, for $s \in c \subset \mathcal{R}$,
is still $\mathcal{H}$, and the results above can be easily extended to
this case.

The case of non-transcendantal parameters (or algebraically
dependant parameters in the case of several parameters), and the specially
interesting case of roots of 1, cannot a priori be dealt with by the methods used here
(although the knowledge of the generic case should be useful for
dealing with the specialized ones).

The natural generalization of this work is thus to deal with
the pseudo-reflection groups. In that case, several parameters
immediately arise, and the general preliminary study is more complicated
than the (2-)reflection case. We will deal with this more general
setting in a forthcoming paper (which will built on the
combinatorical case-by-case arguments and the results exposed here).

\medskip
Throughout this paper, the name CHEVIE refers to the corresponding
GAP package, which can be downloaded at \url{http://www.math.jussieu.fr/~jmichel}.

\medskip

\noindent {\bf Acknowledgements.} I thank C. Cornut and
J. Michel for numerous fruitful exchanges. I also thank
C. Bonnafé, N. Matringe, D. Mauger and J.-F. Planchat for
useful discussions.

\section{Infinitesimal Iwahori-Hecke algebras} 

Let $W$ be a finite complex (2-)reflection group, and
$\eps : W \to \{ \pm 1 \}$ its
sign character, afforded by the determinant.

\subsection{Definitions and basic properties}

For now, we only assume that $\k$ is a field of characteristic 0.

\begin{defi} The infinitesimal Iwahori-Hecke algebra associated
to a finite complex (2-)reflection group $W$ is
the Lie subalgebra $\mathcal{H}_W$
of $\k W$ generated by the set of reflections of $W$.
\end{defi}

Note that, if $W_0$ is a reflection subgroup of $W$, that is if $W_0$
is generated by reflections of $W$, then
$\mathcal{H}_{W_0}$ embeds in $\mathcal{H}_W$ as a Lie subalgebra.
If $W \simeq W_1 \times W_2$ as complex reflection groups, then
$\mathcal{H}_W \simeq \mathcal{H}_{W_1} \times \mathcal{H}_{W_2}$.

Let $\mathcal{R}$ denote the set of all reflections in $W$, and
$\mathcal{S} \subset
\mathcal{R}$ a subset satisfying
\begin{enumerate}
\item $\mathcal{S}$ generates $W$
\item Every reflection of $W$ is conjugated in $W$ to some element
of $\mathcal{S}$.
\end{enumerate}
A typical example for $\mathcal{S}$ is the set of simple
reflections attached to some Weyl chamber when $W$ is a Coxeter
group. We let $W$ act on $\mathcal{R}$
by conjugation, and $\mathcal{R}/W$ the corresponding set of equivalence
classes. It is the set of conjugacy classes of $W$ whose elements
are reflections. For $c \in \mathcal{R}/W$ we denote
$T_c = \sum_{s \in c} s
\in \k W$, and introduce $p : \k G \to Z(\k G)$ defined by
$p(x) = (1/\# G) \sum_{g \in G} gxg^{-1}$.

For an arbitrary finite group $G$, Lie subalgebras of $\k G$ generated by
generating sets of $G$ have noticeable properties, that we recall
from \cite{IRRED}, lemme 3.

\begin{prop} \label{propreductivecenter} Let $S$ be a generating set of the finite group $G$. Then
the Lie subalgebra $\mathcal{L}$ of $\k G$ generated by $S$ is
reductive, and every irreducible representation of $G$
is irreducible for the action of $\mathcal{L}$. Moreover, the center of
$\mathcal{L}$ is contained in the sub-vector space
of $\mathcal{L}$ spanned by the $T_c$ for $T_c \in L$ and
$c \cap S \neq \emptyset$. If $S$ is the union of some conjugacy classes,
then $Z(\mathcal{L}) = p(\mathcal{L})$ and $\mathcal{L}'$ is generated by the
elements $s - T_{c(s)}/\# c(s)$ for $s \in S$.
\end{prop}
\begin{proof} The fact that irreducible representations of $G$
gives rise to irreducible representations of $\mathcal{L}$
comes from the fact that $S$ generates both $\mathcal{L}$
(as a Lie algebra) and $G$ (as a group). 
Now any faithful representation of $\k G$ restricts to a
faithful, semisimple representation of $\mathcal{L}$, so
$\mathcal{L}$ is reductive.
Another consequence of the fact that $S$ generates both $\mathcal{L}$
and $G$
is that $Z(\mathcal{L}) \subset Z(\k G)$.

We define $p : \k G \onto Z(\k G)$ by $p(x) = (1/\# G) \sum_{g \in G} gxg^{-1}$.
It is easily checked (see e.g. \cite{galglie}) that $(\k G)' = \Ker\, p = \cap_c \Ker \delta_c$
for $c$ a conjugacy class of $G$ and $\delta_c$ the characteristic
linear form associated to it.
Let $E$ be the vector space spanned
by $S$ and $\mathcal{L}'$ the vector space generated by iterated brackets of
elements of $\mathcal{L}$. Obviously $\mathcal{L} = E + \mathcal{L}'$ and
$\mathcal{L}' \subset (\k G)' = \cap \Ker \delta_c$ so
$\delta_c(\mathcal{L}) = \delta_c(E)$ for any conjugacy class $c$.
It follows that $\delta_c(\mathcal{L}) = \{
0 \}$
for every class such that $c \cap S = \emptyset$. Since $Z(\k G)$
is spanned by the elements $T_c$ the conclusion follows.

Finally, if $S$ is the union of conjugacy classes then $\mathcal{L}$
is stable under conjugation by $W$ hence $p$ restricts to
a linear endomorphism of $\mathcal{L}$.
Since $Z(\mathcal{L}) \subset Z(\k G)$ we know that $p$ acts identically
on $Z(\mathcal{L})$. From $\mathcal{L}' \subset (\k G)' = \Ker p$
and $\mathcal{L} = Z(\mathcal{L}) \oplus \mathcal{L}'$ we thus
deduce $p(\mathcal{L}) = Z(\mathcal{L})$. 
Now the collection $s - T_{c(s)}/\# c(s)$ for $s \in S$
generates a Lie algebra that clearly contains $\mathcal{L}'$ (because
each $T_{c(s)}$ lies in $Z(\k G)$) and
is contained in $\mathcal{L}$. Since $p(s - T_{c(s)}/\# c(s)) = 0$
it is contained in $\mathcal{L}'$, which concludes the proof. 
\end{proof}


The following facts are mostly recollected from \cite{HECKINF}. We provide a proof for
the convenience of the reader. For $s \in \mathcal{R}$
we let $s' = s - T_{c(s)}/\# c(s)$, where $c(s) \in \mathcal{R}/W$
denotes the conjugacy class of $s$, and recall from the introduction
that $\mathcal{L}_{\eps}(W)$ is the Lie subalgebra of $\k W$ spanned
by the $w - \eps(w) w^{-1}$, $w \in W$.

\begin{prop} \label{propgeninf} Let  $W$ be a finite complex reflection group. Then
\begin{enumerate}
\item $\mathcal{H}_W$ is a reductive Lie algebra.
\item The center of $\mathcal{H}_W$ has dimension $\# \mathcal{R}/W$ with
basis $\{ T_c ; c \in \mathcal{R}/W \}$.
\item The derived Lie algebra $\mathcal{H}_W'$ of $\mathcal{H}_W$
is generated by the $s' = s - T_{c(s)}/\#c(s)$ for $s \in \mathcal{R}$.
\item $\mathcal{H}_W \subset \mathcal{L}_{\eps}(W)$
\item $\mathcal{H}_W$ is stable in $\k W$ for the adjoint action of $W$,
and is generated by $\mathcal{S}$.
\item $\mathcal{H}_W'$ is stable in $\k W$ for the adjoint action of $W$,
and is generated by $\{s' \ | \ s \in \mathcal{S} \}$.
\end{enumerate}
\end{prop}
\begin{proof}
The first three claims are consequences of proposition \ref{propreductivecenter},
since $W$ is generated by $\mathcal{R}$.
If $s \in \mathcal{R}$, then $\eps(s) = -1$ hence $s = \frac{s-\eps(s)}{2}
\in \mathcal{L}_{\eps}(W)$. It follows that $\mathcal{H}_W \subset
\mathcal{L}_{\eps}(W)$. 
Finally, since $\mathcal{R}$ is stable
under conjugation then so is $\mathcal{H}_W$. For $s \in \mathcal{R}$ let
$\mathrm{Ad}(s) : y \mapsto s y s^{-1} = sys$ and $\mathrm{ad}(s)
: y \mapsto
[s,y]$ associated endomorphisms of $\mathcal{H}_W$.
It is easily checked that, for all $s \in \mathcal{R}$, $\mathrm{ad}(s)^2 = 
2 (\mathrm{Id} - \mathrm{Ad}(s))$. In particular, since $\mathcal{S}$
generates $W$, the Lie
algebra generated by $\mathcal{S}$ is stable under the action of $W$.
Since the orbit under conjugation of $\mathcal{S}$ is $\mathcal{R}$
it follows that this Lie algebra contains $\mathcal{R}$ hence
equals $\mathcal{H}_W$. The proof of (7) goes along the same lines.
\end{proof}

From now on, we assume that $W$ is irreducible, and that $\k$ contains
the field of definition of $W$. Multiplicative characters will play
a special role. We recall a consequence of Stanley theorem
(see \cite{STANLEY}, theorem 3.1), denoting $\{ \pm 1 \}^{\RR/W}$ the set
of maps from the set $\RR/W$ of reflection classes to
$\{ \pm 1 \}$.

\begin{lemma} \label{lembijconjchar} There is a natural bijection $\{ \pm 1 \}^{\mathcal{R}/W}
\to \Hom(W, \{ \pm 1 \})$ which associates to each $f \in 
\{ \pm 1 \}^{\mathcal{R}/W}$ a character $\eta : W \to \{ \pm 1 \}$
such that $\eta(s) = f(c)$ whenever $c \in \mathcal{R}/W$ and $s \in c$.
\end{lemma}

Although we do not use it now, we recall that a consequence of the Shephard-Todd classification is that
$\# \mathcal{R}/W \leq 3$ when $W$ is irreducible, with
moreover $\# \mathcal{R}/W \leq 2$ if $\rk(W) \geq 3$.
For $\rho \in \Irr(W)$ we let
$$
\mathrm{X}(\rho) = \{ \eta \in \Hom(W,\{ \pm 1 \}) \ | \ 
\forall s \in \mathcal{R} \ \ \eta(s) = -1 \Rightarrow \rho(s) = \pm 1 \}.
$$

We first need a lemma

\begin{lemma} \label{hyperbolic} Let $E$ be a finite-dimensional
$\k$-vector space, for an arbitrary characteristic 0 field $\k$.
If $E$ admits a nondegenerate symmetric bilinear form $< \ , \ >$
and an involutive skew-isometry, i.e. $s \in \End(E)$ with $s^2 = 1$ and $< 
sx,sy > = -<x,y>$ for every $x, y \in E$, then $< \ , \ >$ is hyperbolic.
\end{lemma}
\begin{proof} Whenever $x \neq 0$ we have $<sx ,x > = < sx, ssx> = - <x,sx>$
hence $<sx,x> = 0$. Since $< \ , \ >$ is nondegenerate, it admits
a non-isotropic $x \in E$ with $<x,x> \neq 0$. If $s x = \la x$
for some $\la \in \k$, then $-<x,x> = <sx,sx> = \la^2 <x,x>$ and
$\la^2 = -1$, contradicting $s^2 = 1$. Thus $x,sx$ are linearly
independant. Let $u = x+sx$, $v = (x-sx)/2<x,x>$. The plane
$H$ spanned by $u,v$ is hyperbolic, $E$ is an orthogonal direct
sum of $H$ and its orthogonal, which is $s$-invariant because
$H$ is so. The conclusion then follows by induction on the dimension of $E$.
\end{proof}

\begin{prop} \label{propospinf} If there exists $\eta \in \mathrm{X}(\rho)$ such that
$\rho^* \otimes \eps \simeq \rho \otimes \eta$, this $\eta$ is
unique. The associated embedding $\eps \otimes \eta \into
\rho^* \otimes \rho^*$ defines (up to scalar) a bilinear form
on $V_{\rho}$, and we denote by $\osp(V_{\rho}) \subset \sl(V_{\rho})$
the Lie subalgebra preserving it. We have $\rho(\HH') \subset \osp(V_{\rho})$.
If the bilinear form is symmetric, then it is hyperbolic
(provided $\dim \rho > 1$).
\end{prop}

\begin{proof} Let $\eta_1,\eta_2 \in \XX(\rho)$
with $\rho \otimes \eta_1 \simeq \rho \otimes \eta_2$,
and let $s \in \RR$. If $\rho(s) \not\in \{ \pm 1 \}$
then $\eta_1(s) = 1 = \eta_2(s)$. If $\rho(s) \in \{ \pm 1 \}$
then $\rho(s) \eta_1(s) = \rho(s) \eta_2(s)$ hence $\eta_1(s) = 
\eta_2(s)$. Since $\RR$ generates $W$ we get $\eta_1 = \eta_2$.

Let now $\eta \in \XX(\rho)$ with $\rho^* \otimes \eps \simeq
\rho \otimes \eta$ and $<\ , \ >$ the corresponding
bilinear form on $V_{\rho}$. Recall that $\rho(\HH') \subset \sl(V_{\rho})$
is generated by the $\rho(s'), s \in \RR$. For all $s \in \RR$ we have
$<\rho(s)x,\rho(s)y> = \eps(s) \eta(s) <x,  y>$ for all $x,y \in V_{\rho}$.
Since $s^2 = 1$ this means 
$<\rho(s)x,y> = \eps(s) \eta(s) <x,  \rho(s)y>$. When $\rho(s) \not\in
\{ \pm 1 \}$ we get $<\rho(s)x,y> + <x,  \rho(s)y>= 0$
hence $<\rho(s')x,y> + <x,  \rho(s')y>= 0$ and $\rho(s') \in \osp(V_{\rho})$.. When $\rho(s) \in \pm 1$
we have $\rho(s') = 0$ hence again $\rho(s') \in \osp(V_{\rho})$
and $\rho(\HH') \subset \osp(V_{\rho})$.
The last assertion is a consequence of lemma \ref{hyperbolic} as
involutive skew-isometries are provided by the
$s \in \RR$ with $\eta(s) = 1$. If not such $s$ exists,
then $\rho(W)$ is abelian and $\dim V_{\rho} = 1$.
\end{proof}

Obviously $\osp(V_{\rho})$ is either a special orthogonal Lie algebra,
when $\eps \otimes \eta \into S^2 \rho^*$, or a symplectic one,
when $\eps \otimes \eta \into \Lambda^2 \rho^*$. Note that
the hyperbolicity property implies that the orthogonal Lie algebras
possibly involved here are even ones.


Any representation $\rho$ of $W$ restricts to a representation
$\rho_{\mathcal{H}}$ of $\mathcal{H}_W$, and in particular to
a representation $\rho_{\mathcal{H}'}$ of $\mathcal{H}_W'$.
We denote $\Irr(\mathcal{H}_W)$ and
$\Irr(\mathcal{H}'_W)$ the set of (classes of)
irreducible representations of $\mathcal{H}_W$ and $\mathcal{H}_W'$,
respectively.

The following proposition is a generalization of \cite{LIETRANSP}, prop. 2.

\begin{prop} \label{propisorep} Let $W$ be a finite reflection
group. Then
\begin{enumerate}
\item If $\rho \in \Irr(W)$ then $\rho_{\mathcal{H}} \in \Irr(\mathcal{H}_W)$.
\item If $\rho^1,\rho^2 \in \Irr(W)$ then $\rho^1 \simeq \rho^2 \Leftrightarrow
\rho^1_{\mathcal{H}} \simeq \rho^2_{\mathcal{H}}$.
\item If $\rho^1,\rho^2 \in \Irr(W)$ with $\dim \rho^1 > 1$ 
then $ \rho^1_{\mathcal{H}'} \simeq \rho^2_{\mathcal{H}'}$ if and only if
$\rho^1 \simeq \rho^2\otimes \eta$ 
for some $\eta \in \XX(\rho_2) = \XX(\rho_1)$.
\item $(\rho \otimes \eps)_{\HH} \simeq (\rho_{\HH})^*$.
\end{enumerate}
\end{prop}
\begin{proof}
For $\mathcal{R}/W = \{ \mathcal{R}_1, \dots,\mathcal{R}_k \}$, we
let $T_i = T_{\mathcal{R}_i}$
and let $\eta_i : W \to \{ \pm 1 \}$ be defined by $\eta_i(\mathcal{R}_i) = -1$
and $\eta_i(\mathcal{R}_j) = 1$ for $j \neq i$.

The first two claims are immediate consequences of the fact that $\mathcal{R}$
generates $W$. The last one is a consequence of $\HH \subset \mathcal{L}_{\eps}(W)$
(see prop. \ref{propgeninf}), and of the corresponding statement for
$\mathcal{L}_{\eps}(W)$.

Now assume $\rho^1_{\mathcal{H}'} \simeq \rho^2_{\mathcal{H}'}$.
$\mathcal{H}_W'$ is generated by the elements $s - T_i/\# \mathcal{R}_i$
$1 \leq i \leq \# \mathcal{R}/W$ and $s \in \mathcal{R}_i$.
Since all irreducible representations of $W$ over $\k$
are absolutely irreducible, we know that $\rho^1(T_i)$ and $\rho^2(T_i)$
are scalars. It follows that there exists
$P \in \Hom(\rho^1, \rho^2)$ such
that $P \rho^1(s) P^{-1} = \rho^2(s) + \om_i$, with $\om_i
 = \frac{\rho^2(T_i) - \rho^1(T_i)}{\# \mathcal{R}_i} \in \k$
for all $1 \leq i \leq \# \mathcal{R}/W$ and $s \in \mathcal{R}_i$.
From $s^2 = 1$ we get
$1 = 1 + 2 \om_i \rho^2(s) + \om_i^2$. Let $I = \{ i \ | \ \om_i \neq 0 \}$.
For $i \in I$ we have $\rho^2(s) = \frac{-\om_i}{2} \in \k$
hence $\rho^1(s)
= \frac{\om_i}{2} = - \rho^2(s)$ and $P \rho^1(s) P^{-1} = - \rho^2(s)$.
For $i \not\in I$ we have $P \rho^1(s) P^{-1}= \rho^2(s)$ hence
$\rho^1 \simeq \rho^2 \otimes \eta$ where $\eta = \prod_{i \in I} \eta_i$,
because $\mathcal{R}$ generates $W$.
Moreover, if $\om_i \neq 0$ and $s \in \mathcal{R}_i$, then
$s^2 = 1$ implies $\om_i = \pm 2$ and this concludes the proof.
Conversely,
if $\rho^2 \simeq \rho^1 \otimes \eta$ and $\rho^1(s) = \pm 1$ whenever
$\eta(s)=-1$, then for such $s$ we have that $s - T_{c(s)}/\#c(s)$
has zero image under both $\rho^1$ and $\rho^2$, and
these images are obviously conjugated if $\eta(s)=1$. It follows
that $\rho^1_{\mathcal{H}'} \simeq \rho^2_{\mathcal{H}'}$. 
\end{proof}
 
The following consequence is obvious.

\begin{cor} \label{corisorep} If $W$ admits a single conjugacy class of reflections
and $\rho^1,\rho^2 \in \Irr(W)$, then $\rho^1 \simeq \rho^2$ iff
$\rho_{\mathcal{H}'}^1 \simeq \rho_{\mathcal{H}'}^2$.
\end{cor}

We use $\mathrm{X}(\rho)$ to define an equivalence relation on $\Irr(W)$.

\begin{defi} For $\rho_1, \rho_2$ we define $\rho_1 \approx \rho_2$
by $\rho_2 \in \{ \rho_1\otimes \eta, \rho_1^* \otimes \eta \otimes \eps \ | 
\ \eta \in \mathrm{X}(\rho_1) \}$.
\end{defi}

Notice that $\mathrm{X}(\rho)$ is a subgroup of $\Hom(W , \{ \pm 1 \})$
and that $\rho_1 \approx \rho_2 \Rightarrow \mathrm{X}(\rho_1)
= \mathrm{X}(\rho_2)$. It follows that $\approx$ is an equivalence
relation on $\Irr(W)$, in general coarser than the relation $\sim$ of the
introduction, generated by $\rho^* \otimes \eps \sim \rho$.
Also notice that $\eps \in \mathrm{X}(\rho)$ if and only if $\rho(s) = \pm 1$
for all $s \in \mathcal{R}$. This is possible only when $\dim \rho = 1$.
In particular, if $W$ has a single class of reflections and $\rho^1,\rho^2 \in
\Irr(W)$ whith $\dim \rho^i > 1$, then $\rho^1 \approx \rho^2$ iff
$\rho^1 \sim \rho^2$.

\begin{defi} \label{defHrho} If $\rho \in \Irr(W)$, we let $\mathcal{H}(\rho)$
denote the orthogonal of $\Ker \rho_{\mathcal{H}'}$ with respect
to the Killing form of $\mathcal{H}'$, and let $\mathcal{L}
(\rho) = \sl(V_{\rho})$ if $\rho \otimes \eta \not \simeq \rho^* \otimes \eps$
for every $\eta \in \XX(\rho)$, and $\mathcal{L}(\rho) = \osp(V_{\rho})$ otherwise.
\end{defi}

It is clear that $\mathcal{H}(\rho)$ is an ideal of $\mathcal{H}'$,
isomorphic to $\rho(\mathcal{H}')$ as a Lie algebra. We prove
that all $\mathcal{L}(\rho)$ are semisimple. By their definition,
this amounts to saying that the exceptional case $\mathcal{L}(\rho)
\simeq \so_2 \simeq \k$ with $\dim \rho = 2$ never occurs (however the non-simple
case $\so_4 \simeq \sl_2 \times \sl_2$ \emph{do} occur, see
lemma \ref{lemexcmrd} below). The reason is that, if $\dim \rho = 2$,
then $\Lambda^2 \rho^* = \det \circ \rho^* \in \Hom(W,\{ \pm 1 \})$
can be written as $\eps \otimes \eta$ for some $\eta \in \Hom(W, \{ \pm 1 \}$),
and we have $\eta \in \XX(\rho)$, as $\rho(s) \neq \pm 1$
implies $\det \rho(s) = -1 = \eps(s)$. Also notice
that $\rho(\HH')$ is semisimple and nonzero, otherwise $\rho(s) \in \{ \pm 1 \}$
for all $s \in \RR$ and $\dim \rho = 1$.
By $\sl_2(\k) \simeq \sp_2(\k)$ this implies the following.

\begin{lemma} \label{lemexcdim2} Let $\rho \in \Irr(W)$ with $\dim(\rho) = 2$.
Then $\mathcal{L}(\rho) = \sl(V_{\rho}) = \rho(\HH') \simeq \HH(\rho)$.
\end{lemma}



\subsection{Reflection representations}

From now on we assume that $W$ is an irreducible reflection group
of rank at least 2.
We need a slightly nonstandard definition.

\begin{defi} \label{defreflrep}
An irreducible representation $R : W \to \GL(V_R)$ with
$\dim R \geq 2$ is called a \emph{reflection representation}
if, for all $s \in \mathcal{R}$, 
if $R(s) \neq  \Id$ then
$R(s)$ is a reflection of $V_R$.
\end{defi}

By the irreducibility assumption on $W$, such representations
exist as soon as $\rk(W) \geq 2$, and at least one of them
is faithful. All faithful reflection representations
of $W$ have for dimension the rank of $W$ by the
classification theorem of Shephard and Todd, and they are all
deduced from the defining representation by Galois action (see \cite{MARINMICHEL}).

It is a classical fact due to Steinberg
(see \cite{BOURB456}, ch. V \S 2 exercice 3)
that the alternating powers $\Lambda^k R$
of $R$ for $1 \leq k \leq \dim R$ of a reflection representation are
distinct irreducible representations of $W$, thus leading to distinct
simple ideals of $\k W$ as an associative algebra ; indeed, it is
straightforward to check that the proof of \cite{BOURB456}
applies to our more general definition of a reflection representation.

Taking $k$-th alternating powers gives rise to two representations
of $\mathcal{H}_W$ on $\Lambda^k V_R$, namely $(\Lambda^k R)_{\mathcal{H}}$
and $\Lambda^k (R_{\mathcal{H}})$. Likewise, suppose we are given
$\eta \in \Hom(W, \{ \pm 1 \})$. Then, to any $x \in \k$ one can
associate a character $\gamma_x^{\eta}$ of $\mathcal{H}_W$
defined by $s \mapsto \eta(s) (x-1)$ for $s \in \mathcal{R}$.
Letting $R^{\eta} = R \otimes \eta$,
one may then define two twisted representations
$(\eta \otimes \Lambda^k R)_{\mathcal{H}}$ and
$\Lambda^k (R^{\eta}_{\mathcal{H}})$. Identifying $\k \otimes \Lambda^k
V_R$ with $\Lambda^k V_R$, all these representations as well as $\gamma_x^{\eta}
\otimes (\eta \otimes \Lambda^k R)_{\mathcal{H}}$ act on
$\Lambda^k V_R$.

\begin{prop} \label{proprefrep}
Let $R : W \to GL(V_R)$ be a reflection representation of a finite
irreducible reflection group $W$ and $\eta \in \Hom(W,\{\pm 1 \})$. For any
$k \in [0 ,\dim V_R ]$,
\begin{enumerate}
\item $\Lambda^k (R^{\eta}_{\mathcal{H}}) = \gamma_k^\eta \otimes
(\eta \otimes \Lambda^k R)_{\mathcal{H}}$
\item $\Lambda^k (R^{\eta}_{\mathcal{H}'}) = (\eta \otimes
\Lambda^k R)_{\mathcal{H}'}$
\end{enumerate}
and, in particular, $\Lambda^k R_{\mathcal{H}'} = (\Lambda^k R)_{\mathcal{H}'}$.
\end{prop}
\begin{proof}
Let $s \in \mathcal{R}$. If $R(s) = 1$, then $(\eta \otimes
\Lambda^k R)(s) = \eta(s)$ and
$(\Lambda^k R_{\mathcal{H}}^{\eta})(s) = k \eta(s) = \eta(s)(k-1)\Id
+(\eta \otimes \Lambda^k R)(s)$. Otherwise,
there exists a basis $e_1,\dots,e_n$ of $V_R$
such that
$s.e_1 = -e_1$ and $s.e_i = e_i$ if $i \neq 1$, where we make $W$
act on $V_R$ through $R$. One has a natural basis $(e_I)$ of
$\Lambda^k V_R$ which is indexed by the subsets of size $k$ in
$\{1,\dots,n \}$ : if $I = \{i_1,\dots,i_k \}$ with $i_1 < \dots < i_k$
then $e_I = e_{i_1} \wedge \dots \wedge e_{i_k}$. Then
$$
\left\lbrace \begin{array}{lcrr}
(\eta \otimes \Lambda^k R)(s)(e_I) & = & \eta(s) e_I & \mbox{ if } 1 \not\in I \\
 & = & -\eta(s) e_I & \mbox{ if } 1 \in I \\
\end{array} \right.$$
{}
$$ 
\left\lbrace \begin{array}{lcrr}
(\Lambda^k R^{\eta}_{\mathcal{H}})(s)(e_I) & = & k \eta(s) e_I & \mbox{ if } 1 \not\in I \\
 & = & (k-2)\eta(s) e_I & \mbox{ if } 1 \in I \\
\end{array} \right. 
$$
and it follows that $(\Lambda^k R^{\eta}_{\mathcal{H}})(s)
= \eta(s) (k-1) \mathrm{Id} + (\eta \otimes \Lambda^k R)(s)$,
thus proving (1). Statement (2) follows immediately.
\end{proof}

\begin{remark} \label{remlrefapprox}
\end{remark}
It is a standard fact that, if $R$ is faithful, then
$\eps \otimes \Lambda^k R^* = \Lambda^{n-k}
R$, where $n = \dim V_R$. In general, we easily check that $\Lambda^{n-k}
R = \eta \otimes \eps \otimes \Lambda^k R^*$
with $\eta \otimes \eps(g) = \det R(g)^{-1}$ for $g \in W$.
For $s \in \mathcal{R}$, this means that $\eta(s) = -1$ iff $\det R(s) = 1$.
In particular, $\eta \in \mathrm{X}(R)$.

\medskip

Let $L(V_R) = \bigoplus_{k=0}^{\infty} \sl(\Lambda^k V_R)$, and
$\Lambda^{\bullet} : \sl(V_R) \to L(V_R)$ be the
natural representation of $\sl(V_R)$ on the exterior algebra
of $V_R$. By considering $\Lambda^k V_R$ as the underlying vector space
of $\Lambda^k R$, there is a natural embedding $\iota : L(V_R)
\into (\k W)' \subset \k W$ and, more generally, to
any $\eta \in \Hom(W, \{ \pm 1 \})$ is associated a natural
embedding $\iota_{\eta} : L(V_R) \into (\k W)'$ by identification
of $\Lambda^k V_R$ with the underlying vector space of $\eta \otimes
\Lambda^k R$. 

Proposition \ref{proprefrep} leads to considering the following diagram
$$
\xymatrix{
\mathcal{H}'_{W} \ar[r] \ar[dr] \ar[d]_{R_{\mathcal{H}'}} & \k W \\
\sl(V_R) \ar[r]_{\Lambda^{\bullet}} & L(V_R) \ar[u]_{\iota}
}
$$
where the morphism $\mathcal{H}'_W \to L(V_R)$ is defined in the obvious
way such that the upper right triangle commutes. Proposition \ref{proprefrep}
implies that the lower left triangle also commutes. We
showed in \cite{IRRED} by a case-by-case analysis
that the morphism $\mathcal{H}_W' \to \sl(V)$ is surjective
for finite irreducible Coxeter groups.
More generally, we have the following result.

\begin{prop} \label{proprefsurj}
Let $R : W \to GL(V_R)$ be a reflection representation of a finite
irreducible reflection group $W$ and $\eta \in \Hom(W,\{ \pm 1 \})$.
Then $R^{\eta}_{\mathcal{H}'} : \mathcal{H}'_W \to \sl(V_R)$ is
surjective.
\end{prop}

To prove this, we will need the following classical lemmas.
\begin{lemma} \label{lemmAroot}A root system of type $A_n$ admits no proper subsystem
of rank $n$.
\end{lemma}
\begin{proof} We realize a root system of type $A_n$ as the elements
$e_i -e_j$ in the euclidean vector space $\R^{n+1}$ with basis
$e_1,\dots,e_{n+1}$. We let $s_{ij}$ denote the reflection corresponding
to the root $e_i-e_j$. Let $E$ be a subsystem of rank $n$ and let
$\mathcal{E} = \{ i \mid e_1-e_i \in E \}$, $\mathcal{F} = \{2,\dots,n+1 \}
\setminus \mathcal{E}$. Assume $\mathcal{F}\neq \emptyset$. Then $e_i - e_j \not\in E$
for all $i \in \mathcal{E}$ and $j \in \mathcal{F}$ otherwise
$s_{1i}(e_i-e_j)=e_1-e_j \in E$ and $j \in \mathcal{E} \cap
\mathcal{F} = \emptyset$. Let $\mathcal{E}'=\mathcal{E}\cup \{ 1 \}$.
Then $E$ is contained in the vector space generated by the $e_i - e_j$
for $i,j \in \mathcal{E}'$ or $i,j \in \mathcal{F}$, whose dimension
is at most $n-1$. By contradiction, it follows that $\mathcal{F}=\emptyset$
and $E$ contains all $e_1 - e_i$. Since the permutations $(1\ i)$
generate $\mathfrak{S}_{n+1}$, the Weyl groups of the two subsystems
are the same. Since this Weyl group acts transitively on the set
of roots and $E \neq \emptyset$ it follows that $E$ is not proper.
\end{proof}

For $V$ a finite-dimensional $\C$-vector space and $W \subset \GL(V)$
a finite group acting irreducibly, we endow $V$ with a $W$-invariant
hermitian scalar product. If $H \subset V$ is a hyperplane with
orthogonal spanned by $e_H \in V$, we let $W_H = \{ w \in W \ | \ 
w.e_H = e_H \}$. We recall the following lemma from \cite{arrrefl}.

\begin{lemma} (see \cite{arrrefl}, cor. 3.2)\label{maxparabol} Let $W \subset \GL(V)$ be a finite group generated
by pseudo-reflections and acting irreducibly on $V$. Let $\mathcal{R}$
denote its set of pseudo-reflections. There exists an hyperplane
$H$ of $V$ such that $W_H$ acts irreducibly on $H$ and such that
its image in $\GL(H)$ is generated by the images of $\mathcal{R} \cap W_H$.
\end{lemma}

We are now in position to prove proposition \ref{proprefsurj}.

\begin{proof}
Since, for all $s \in \mathcal{R}$, we have $R^{\eta}(s) = \eta(s) R(s)$,
it suffices to show that $R_{\mathcal{H}'}$ is surjective.
We will thus assume $\eta=1$. Since $R(W) \subset \GL(V_R)$ is
an irreducible reflection group admitting as reflections the image
of $\mathcal{R}$, we can assume that $R$ is faithful and is the defining
representation of $W$, that is $V_R = V$ and $W \subset \GL(V)$.

Under these assumptions, we prove the proposition by induction on $n=
\rk W = \dim V$. For convenience here, we drop the
implicit assumption $\dim V \geq 2$ in the statement
and start the induction at the trivial case $n = 1$.
We proceed by
assuming $n \geq 2$. Let $H \subset V$ an hyperplane and $W_0=
W_H \subset W$ being afforded by lemma \ref{maxparabol}. We have
$\mathcal{H}_{W_0} \subset \mathcal{H}_{W}$
hence $\mathcal{H}_{W_0}' \subset \mathcal{H}_{W}'$ and
$\Im R_{\mathcal{H}'} \supset \Im R_{\mathcal{H}_{W_0}'}
= \sl(H)$ by the induction hypothesis. It follows that there exists
a Cartan subalgebra of rank $n-2$ of $\Im R_{\mathcal{H}'}$
inside $\sl(H)$.
We are going to exhibit a semisimple element $x \in \sl(V)
\setminus \sl(H)$ centralizing $\sl(H)$ that belongs to the image
of $\mathcal{H}_W'$ : this will provide a Cartan subalgebra of rank
$n-1$ for this image, which then equals $\sl(V)$ since it
is a semisimple Lie subalgebra of $\sl(V)$ and the root system of
type $A_{n-1}$, corresponding
to $\sl(V)$, contains no proper subsystem of rank $n-1$
by lemma \ref{lemmAroot}.

We now construct $x$. Let $\mathcal{R}_0$ denote
the set of reflections of $W_0$, and define the following
elements in $\mathcal{H}_W$
$$
T = \sum_{s \in \mathcal{R}} s,\ \ \ 
T_0 = \sum_{s \in \mathcal{R}_0} s,\ \ \ 
X = (\# \mathcal{R}) T_0 - (\# \mathcal{R}_0 ) T
$$
It is clear that $X$ belongs to $\mathcal{H}_W$ and
that it commutes to $\mathcal{H}_{W_0}$. Let
$x = R_{\mathcal{H}}(X)$. We know that $T$ acts
on $V$ by the scalar $\frac{n-2}{n}\# \mathcal{R}$, because $tr R(s) = n-2$ for all $s \in \mathcal{R}$
and $R$ is irreducible. Similarly, $T_0$ acts on $H$
by $\frac{n-3}{n-1} \# \mathcal{R}_0$ and on
$D$ by $\# \mathcal{R}_0$. It follows that $x$ acts on $H$ by
$(\# \mathcal{R}) (\# \mathcal{R}_0) \frac{-2}{n(n-1)}$
and by $(\# \mathcal{R}) (\# \mathcal{R}_0) \frac{2}{n}$
on $D$. Thus $x \not\in \sl(H)$ but $x \in \sl(V)$ and $x$ is semisimple.
Since $\mathcal{H}_W$ is reductive, the intersection of $\sl(V)$
with its image in $\gl(V)$ is the image of $\mathcal{H}'_W$. It
follows that $x$ belongs to the image of $\mathcal{H}'_W$, which concludes
the proof.
\end{proof}

\subsection{Reflection ideals}

We define the following subsets of $\Irr(W)$.
$$
\begin{array}{l}
\mathrm{Ref} = \{ \mbox{reflection representations} \} \\
\mathrm{QRef} = \{ \eta \otimes \rho \ | \ \rho \in \mathrm{Ref}, \eta \in
\Hom(W,\{ \pm 1 \}) \} \\
\Lambda\mathrm{Ref} = \{ \eta \otimes \Lambda^k  \rho \ | \ \rho \in \mathrm{Ref}, \eta \in
\Hom(W,\{ \pm 1 \}), k \geq 0 \} \\
\end{array}
$$
and we let $\overline{\mathrm{QRef}} = \mathrm{QRef} / \approx$. By choosing
an arbitrary system of representatives, we identify $\overline{\mathrm{QRef}}$
with a subset of $\mathrm{QRef}$. Notice that, since $\mathrm{Ref} \neq \emptyset$,
then $\Hom(W, \{ \pm 1 \}) \subset \LRef$.
More generally, representations of small dimension belong to $\LRef$.

\begin{lemma} \label{lemref23}
Let $\rho \in \Irr(W)$. If $\dim \rho \leq 3$ then $\rho \in \LRef$.
If moreover $\dim \rho \neq 1$ then $\rho \in \mathrm{QRef}$.
\end{lemma}

\begin{proof}
Since $\Hom(W, \{ \pm 1 \}) \subset \LRef$ it is sufficient to
show the second part of the statement. Assume $\dim \rho = 2$. Let
$\mathcal{R}_0 = \{ s \in \mathcal{R} \ | \ \rho(s) = \Id \}, 
\mathcal{R}_1 = \{ s \in \mathcal{R} \ | \ \rho(s) = -\Id \}, 
\mathcal{R}_2 = \{ s \in \mathcal{R} \ | \ Sp\rho(s) = \{ -1, 1 \} \}$.
Since $\rho(s)^2 = \Id$ for all $s \in \mathcal{R}$
we have $\mathcal{R} = \mathcal{R}_0 \sqcup \mathcal{R}_1 \sqcup \mathcal{R}_2$,
and each $\mathcal{R}_i$ is stable by conjugation. Using lemma \ref{lembijconjchar}
we define $f :  \mathcal{R} \to \{ \pm 1 \}$ by $f(s) = 1$ if $s
\in \mathcal{R}_0 \sqcup \mathcal{R}_2$, $f(s) = -1$ is $s \in \mathcal{R}_1$,
and denote $\eta \in \Hom(W, \{ \pm 1 \}$ the corresponding
character. Then $\eta \otimes \rho \in \mathrm{Ref}$ hence
$\rho \in \mathrm{QRef}$. Similarly, if $\dim \rho = 3$, letting
$\mathcal{R}_0,\mathcal{R}_1,\mathcal{R}_2$ as before, we define
$\mathcal{R}_2^{\pm}$ to be the set of all $s \in \mathcal{R}_2$
such that $\pm 1$ appear in the spectrum of
$\rho(s)$ with multiplicity 2 ; we define $f$ as before on $\mathcal{R}_0
\sqcup \mathcal{R}_1$, $f(s) = \pm 1$ for $s \in \mathcal{R}_2^{\pm}$,
and again $\eta \otimes \rho \in \mathrm{Ref}$ hence
$\rho \in \mathrm{QRef}$.
\end{proof}

\begin{defi} The \emph{reflection ideal} of $\mathcal{H}'$
is
$$
\mathcal{I} = \mathcal{H}' \cap \bigoplus_{\rho \in \Lambda\mathrm{Ref}}
\End(V_{\rho}).
$$
\end{defi}

It is clear that $\mathcal{I}$ is an ideal of $\mathcal{H}'$,
and we have a natural projection $\mathcal{H}' \to \mathcal{I}$ factorizing
the inclusion $\mathcal{H}' \to \k W$.
Let $\rho \in \mathrm{QRef}$. We have $\mathcal{H}(\rho) \simeq
\sl(V_{\rho})$ by proposition \ref{proprefsurj}. Now, proposition
\ref{propisorep} and proposition \ref{proprefrep} imply that
$\mathcal{I}$ is the sum of the ideals $\mathcal{H}(\rho)$ for
$\rho \in \overline{\mathrm{QRef}}$,
since semi-simple Lie algebras are direct sums of their simple ideals
and $\sl(V_{\rho})$ is a simple Lie algebra. 
Finally, if $\rho^1,\rho^2 \in \overline{\mathrm{QRef}}$, $\mathcal{H}(\rho^1)
= \mathcal{H}(\rho^2)$ implies $\rho^1_{\mathcal{H}'} \simeq
\rho^2_{\mathcal{H}'}$ or $\rho^1_{\mathcal{H}'} \simeq (\rho^2_{\mathcal{H}'})^*$,
since $\sl_n$ admits at most two irreducible representations of dimension $n$.
It follows that $\mathcal{H}(\rho^1) = \mathcal{H}(\rho^2)
\Rightarrow \rho^1 \approx \rho^2
\Rightarrow \rho^1 = \rho^2$ by definition of $\overline{\mathrm{QRef}}$,
hence $\mathcal{I}$ is the direct sum of the ideals $\mathcal{H}(\rho)
\simeq \sl(V_{\rho})$ for $\rho \in \overline{\mathrm{QRef}}$.
In particular
we get

\begin{prop} \label{propisorefid} The reflection ideal $\mathcal{I}$
is isomorphic to
$$
\bigoplus_{\rho \in \overline{\mathrm{QRef}}} \sl(V_{\rho}).
$$
\end{prop}

\noindent {\bf Example : the reflection ideal of $G_{13}$.} We
illustrate this decomposition for $W$ of type $G_{13}$. This group
has 4 multiplicative characters $\un ,\eps, \eta, \eta\otimes \eps$,
among 16 irreducible characters. Its field of definition is $\k =\Q(\mu_{24})$. A computer computation shows that
$\dim \mathcal{H}' = 51$.

This group admits 4 faithful 2-dimensional reflection representations,
numbered $\rho_5,\rho_7,\rho_8,\rho_{10}$ in CHEVIE,
and two non-faithful ones, $\rho_6$ in dimension 2 and $\rho_{12}$
in dimension 3. We have $\rho_6(W) \simeq G(1,1,3) \simeq \mathfrak{S}_3$
and $\rho_{12}(W)$ is the Coxeter group of type $B_3$.

There are two natural actions on the 4 faithful ones, the
action by tensor product of $\Hom(W,\{ \pm 1 \})$
and the action of $\Gal(\k|\Q)$. Both are transitive, and
we have $\rho \simeq \rho^* \otimes \eps$ and $\mathrm{X}(\rho) = \{ \un \}$
for all of them. 
It follows that they all belong to $\overline{\mathrm{QRef}}$,
and that no other representation of $\mathrm{QRef}$ arise from them.
Moreover, we have $\iota(L(V_{\rho})) = \sl(V_{\rho})$
for all of them.

In dimension 2, the other one is $\rho_6 \simeq \rho_6 \otimes \eta \otimes \eps$.
We have $\rho_9 \simeq \rho_6 \otimes \eps \in \mathrm{QRef}$, and both
are defined over $\Q$. In particular $\rho_9 \simeq \eps \otimes \rho_6^*$
and we can assume $\overline{\mathrm{QRef}} \cap \{ \rho_6,\rho_9 \}
= \{ \rho_6 \}$. Another reason for $\rho_9 \approx \rho_6$
is that $\rho_9 = \rho_6 \otimes \eta$ and $\mathrm{X}(\rho_6) = \{ \un, \eta \}$.

Finally, the 3-dimensional one has a real-valued character and induces
3 other representations $\rho_{11} = \rho_{12} \otimes \eta$,
$\rho_{13} = \rho_{12} \otimes \eps$, $\rho_{14} = \rho_{12} \otimes
\eta \otimes \eps$. Since $\mathrm{X}(\rho_{12}) = \{ \un \}$,
we can choose $\overline{\mathrm{QRef}} =
\{ \rho_5,\rho_6,\rho_7,\rho_8,\rho_{10},\rho_{11},\rho_{12} \}$. It follows that 
$\mathcal{I} = \sl_2^5 \times \sl_3^2$ and 
$\dim \mathcal{I} =31$

This example shows in particular that
we cannot assume $\overline{\mathrm{QRef}} \subset \mathrm{Ref}$
in general, since $\rho_{11},\rho_{14} \not\in \mathrm{Ref}$.
Moreover, note that these four last representations of $W$
have the same alternating square, although they may represent
different simple ideals of $\mathcal{H}'$.

\subsection{Decomposition in semisimple components}

We define the following subsets of $\Irr(W)$.
$$
\begin{array}{l}
\mathcal{E} = \{ \rho \in \Irr(W) \ | \ \rho \not\in \Lambda\mathrm{Ref} \mbox{ and } \forall \eta \in \XX(\rho) \ \ \rho^* \otimes \eps  \not\simeq \rho \otimes \eta
  \} \\
\mathcal{F} = \{ \rho \in \Irr(W) \ | \ \rho \not\in \Lambda\mathrm{Ref} \mbox{ and } \exists \eta \in \XX(\rho) \ \ \rho^* \otimes \eps  \simeq \rho \otimes \eta
  \} 
\end{array}
$$
We identify $\mathcal{E}/\approx$ and $\mathcal{F}/\approx$
with subsets of $\mathcal{E}$ and $\mathcal{F}$, respectively.

Proposition \ref{propisorep} and proposition \ref{proprefrep}
imply that the inclusion $\mathcal{H}' \subset (\k W)'$ factorizes
through an injective morphism
$$
\Phi : \mathcal{H}' \into \mathcal{I} \oplus
\left( \bigoplus_{\rho \in \mathcal{E}/\approx} \sl(V_{\rho}) \right) \oplus
\left( \bigoplus_{\rho \in \mathcal{F}/\approx} \osp(V_{\rho}) \right).
$$

The central result of this article is the following

\begin{theo} \label{maintheo} Unless if $W = H_4$ the morphism
$\Phi$ is an isomorphism.
\end{theo}

For dihedral goups, all
irreducible representations except the 1-dimensional ones are
2-dimensional reflection representations, hence $\mathcal{H}'
\simeq \mathcal{I}$. This theorem has been proved for Coxeter
types $A_n$ in \cite{LIETRANSP}, and $I_2(m)$ in \cite{IRRED}.
It is sufficient to prove this theorem for $\k = \C$, which
we assume from now on.

We will show that, in order to prove this theorem for a given $W \neq F_4$,
it is sufficient to show the following :
for all $\rho \in \Irr(W) \setminus \Lambda\mathrm{Ref}$ such
that $\dim \rho > 1$, we have $\rho(\mathcal{H}') = \mathcal{L}(\rho)$
(see corollary \ref{corimpltheo}).
For this we need a few preliminary results.


Using the classification of complex reflection groups
and of their representations, we will provide later
a proof of the following lemma (see lemmas \ref{lemcasF4}, \ref{lempasuncannyeer}, \ref{lemcasA4B4}).
\begin{lemma} \label{lemexcmrd} Let $\rho \in \Irr(W)$.
\begin{enumerate}
\item If $\dim \rho = 4$ and $\eta \otimes \eps \into S^2 \rho^*$ for some $\eta \in \XX(\rho)$ then $W$
is a Coxeter group of type $F_4$.
\item If $\dim \rho = 8$ and $\eta \otimes \eps \into S^2 \rho^*$ for some $\eta \in \XX(\rho)$ then $W$
is a Coxeter group of type $H_4$.
\end{enumerate}
Moreover, if $\dim \rho = 6$ and $\eta \otimes \eps \into S^2 \rho^*$ for some $\eta \in \XX(\rho)$,
then $W$ has rank 4 and $\rho \in \LRef$.
\end{lemma}




\begin{prop} \label{propidsimpl}
Assume $W \not\in \{ F_4,H_4 \}$. If $\rho \in \Irr(W)$ satisfies $\mathcal{L}(\rho) \simeq \mathcal{H}(\rho)$
then $\mathcal{H}(\rho)$ is a simple ideal. If $\rho_1,\rho_2
\in \Irr(W)\setminus \LRef$ satisfy $\mathcal{L}(\rho_i) \simeq \mathcal{H}(\rho_i)$
for $i \in \{ 1,2 \}$
, then
$\rho_1 \approx \rho_2 \Leftrightarrow \mathcal{H}(\rho_1) = 
\mathcal{H}(\rho_2)$.
\end{prop}

\begin{proof} The first part of the statement comes from the
fact that the exceptional (nonsimple) cases $\so_2$ and
$\so_4 \simeq \sl_2 \times \sl_2$ do not occur when $W$ is not of type $F_4$,
by lemmas \ref{lemexcdim2} and \ref{lemexcmrd}.
Assume now $\rho_1,\rho_2 \not\in \LRef$. If
$\rho_1 \approx \rho_2$ then $(\rho_1)_{\mathcal{H}'}
\simeq (\rho_2)_{\mathcal{H}'}$
or 
$(\rho_1)_{\mathcal{H}'}
\simeq ((\rho_2)_{\mathcal{H}'})^*$
by proposition \ref{propisorep} hence
$\mathcal{H}(\rho_1) = 
\mathcal{H}(\rho_2)$. Assume now $\mathcal{L}(\rho_i)
\simeq \mathcal{H}(\rho_i)$ and $\mathcal{H}(\rho_1) = \mathcal{H}(\rho_2)$.
Then $\dim \rho_1 = \dim \rho_2 = N$, because the exceptional case
$\so_6 \simeq \sl_3$ is excluded by the condition $\rho_i \not\in \LRef$,
by lemma \ref{lemexcmrd}, and the case $\so_3 \simeq \sl_2$
is also excluded because only $\so_{2n}$ occur, by the hyperbolicity
property.
We first consider the
case $\mathcal{H}(\rho_1) = \mathcal{H}(\rho_2) \simeq \sl_N(\k)$
with $N \geq 3$. This Lie algebra admits only two non isomorphic
irreducible representations of dimension $N$. But $\rho_1,
\rho_1^* \otimes \eps$ and $\rho_2$ are non-isomorphic
irreducible representations of $\mathcal{H}'$ that factorize through
$\mathcal{H}(\rho_1) = \mathcal{H}(\rho_2)$, unless $\rho_1 \approx
\rho_2$. We now assume $N = 2$ or $\mathcal{H}(\rho_1)=
\mathcal{H}(\rho_2)$ is either orthogonal or symplectic. Then this Lie algebra
admits only one irreducible representation of dimension $N$, hence
$\rho_1 \approx \rho_2$, except if it is of Cartan type $D_4$.
But this does not happen here by lemma \ref{lemexcmrd},
since $W$ is not of type $H_4$.
\end{proof}

\begin{cor} \label{corimpltheo} Assume $W \not\in \{ F_4,H_4 \}$. If $\forall \rho \in \Irr(W)
\setminus \LRef \ \ \mathcal{H}(\rho) \simeq \mathcal{L}(\rho)$, then
$\Phi$ is an isomorphism.
\end{cor}
\begin{proof}
Being semisimple, $\mathcal{H}'$ is a direct sum of its simple
ideals. Moreover $\Phi$ is injective and factorizes, by definition
of $\approx$, through the sum of the $\mathcal{H}(\rho)$,
for $\rho \in \Irr(W)$. We first prove that, if $\rho_1 \in \LRef$ and
$\rho_2 \not\in \LRef$, then $\mathcal{H}(\rho_1) \neq \mathcal{H}(\rho_2)$.
Firstly we have $\rho_1 \not\approx \rho_2$ since, by
remark \ref{remlrefapprox}, $\LRef$ is the union of equivalence
classes for $\approx$. By definition of $\LRef$ and proposition \ref{proprefrep}
we can assume $\rho_1 \in \mathrm{QRef}$,
hence $\mathcal{H}(\rho_1) \simeq \sl_n$
with $n = \dim \rho_1$. If $\mathcal{H}(\rho_2) \simeq \sl_n$
then $n = \dim \rho_2$ since by lemma \ref{lemexcmrd}
the situation $\sl_3 \simeq \so_6$ is excluded by $\rho_2 \not\in \LRef$.
Since $\sl_n$ admits at most two irreducible
representations of dimension $n$, 
$\mathcal{H}(\rho_1) = \mathcal{H}(\rho_2)$ implies that, as
a representation of $\mathcal{H}'$, $\rho_2$ is isomorphic
either to $\rho_1$ or its dual, hence $\rho_1 \approx \rho_2$,
a contradiction.

If $\rho_1,\rho_2 \not\in \LRef$, then proposition \ref{propidsimpl}
says $\rho_1 \approx \rho_2 \Leftrightarrow \mathcal{H}(\rho_1)
= \mathcal{H}(\rho_2)$. By proposition \ref{propisorefid}
it follows that $\mathcal{H}'$ is the direct sum of the $\mathcal{H}(\rho)$
for $\rho \in \overline{\mathrm{QRef}} \sqcup (\mathcal{E}/\approx)\sqcup
(\mathcal{F}/\approx)$. Since $\mathcal{H}(\rho) \simeq \mathcal{L}(\rho)$
for these $\rho$ by proposition \ref{propisorefid} and the hypothesis,
we get that $\mathcal{H}'$ is isomorphic to the direct
sum of these $\mathcal{L}(\rho)$. By equality of dimension
this proves that $\Phi$ is an isomorphism.
\end{proof}

\subsection{Further Lie-theoretic properties}

The following general lemma on $\mathcal{H}(\rho)$ will
be technically useful for our purpose.

\begin{lemma} \label{lempasrk1} Let $W$ be an irreducible complex
reflection group whose reflections are all conjugated.
Assume furthermore that, for all $\rho \in \Irr(W)$,
$\dim \rho \leq 3$ implies $\dim \rho = 1$. Then,
for all $\rho \in \Irr(W)$, $\mathcal{H}(\rho)$ does
not admit a simple ideal of rank 1.
\end{lemma}
\begin{proof}
Let $\rho : W \to \GL(V_{\rho})$ be an irreducible
representation of $W$. We let $\g = \rho(\mathcal{H}') \simeq \mathcal{H}(\rho)$, and assume by contradiction
that $\g = \sl_2(\C) \times \g_0$ as a Lie algebra.
We thus have $\sl_2(\C) \subset \g \subset \gl(V_{\rho})$. Under
the action of $W$, the Lie algebra $\g$ is an invariant subspace
of $\gl(V_{\rho}) = V_{\rho} \otimes V_{\rho}^*$. Moreover, since $\sl_2(\C)$
is an ideal of $\g$, it is invariant by the adjoint action of $\g$
hence by $W$ because of the relation $\mathrm{ad}(s)^2 = 2(\Id - \mathrm{Ad}(s))$
for $s \in \mathcal{R}$. It follows that $\sl_2(\C)$
is a 3-dimensional representation of $W$. By assumption,
it is a direct sum of 1-dimensional ones. On the other hand,
the multiplicity of a 1-dimensional representation
of $W$ in $V_{\rho} \otimes V_{\rho}^*$ is at most 1 by irreducibility
of $V_{\rho}$ and Schur's lemma. It follows that $\sl_2(\C)$
can be decomposed as a direct sum of 3 distinct 1-dimensional
representations of $W$. But there are at most two such
representations since $\# \Hom(W,\{ \pm 1 \}) = 2$, hence a contradiction.
\end{proof}

The argument to prove the theorem is by induction. For this
we need to relate the induction table between two
reflection groups with the induction table of
the corresponding Lie algebras. This is partly done by
the following general lemma. Recall that, when $W$ has a single
(conjugacy) class of reflections, then $\XX(\rho) = \{ \un \}$
as soon as $\dim \rho > 1$.

\begin{lemma} \label{lemmultiplie} Let $W$ be an irreducible reflection group. Assume
that it admits an irreducible reflection subgroup $W_0 \subset W$,
different from $H_4$ and $F_4$, which has a single class
of reflections, and for which the statement of theorem \ref{maintheo}
holds true. Let $\rho \in \Irr(W)$ and $r \geq 1$ such that
$\forall \varphi \in \Irr(W_0) \ (\Res_{W_0} \rho | \varphi) \leq r$.
Then, for any simple Lie ideal $\h$ of $\rho(\mathcal{H}'_{W_0})$,
there exists $J \in \Irr(\h)$ such that $\dim \Hom_{\h} (J, \Res_{\h}
\rho_{\mathcal{H}'} ) \leq r$.
\end{lemma}
\begin{proof} For convenience we denote here $\mathcal{H}'_0
= \mathcal{H}'_{W_0}$, $\mathcal{H}' = \mathcal{H}_W$ and
let $\mathcal{H}_0(\varphi)$ denote the orthogonal of $\Ker
(\varphi_{\mathcal{H}_0'})$
for the Killing form of $\mathcal{H}'$, for $\varphi \in \Irr(W_0)$.
Let $\h$ be a simple Lie ideal of $\rho(\mathcal{H}'_0)$. We let
$a_1 \rho^1 + \dots a_t \rho^t$ be a decomposition of $\Res_{W_0} \rho$
in nonisomorphic simple components, with $1 \leq a_i \leq r$. By proposition
\ref{propisorep} (1) we know that each $\rho^i_{\mathcal{H}'_0}$
is irreducible. We can assume that $\dim \rho^i > 1$ for $i \leq s$
and $\dim \rho^i= 1$ for $s<i\leq t$. The simple Lie
ideals of $\rho(\mathcal{H}'_0)$ are the $\mathcal{H}_0(\rho^i)$
for $i \leq s$ by proposition \ref{propidsimpl}, since $W_0 \not\in
\{ F_4,H_4 \}$ and theorem \ref{maintheo} holds for $W_0$. 
Assuming that there exists a simple Lie ideal $\h$ of $\rho(\mathcal{H}'_0)$,
we have $s \geq 1$ and we can assume $\h = \mathcal{H}_0(\rho^1)$.
We have $\rho^1_{\mathcal{H}'_0}\simeq \rho^i_{\mathcal{H}'_0}$
iff $\rho^i = \rho^1$ by proposition \ref{propisorep} (3), since
$\Hom(W_0, \{ \pm 1 \}) = \{ \un, \eps \}$. Since $\rho^1_{\mathcal{H}'_0}$
is irreducible and $\rho^i_{\mathcal{H}'_0}$ factorizes through
a well-defined simple Lie ideal when $i \leq s$,
it follows that $\dim \Hom_{\h} (\Res_{\h} \rho^1_{\mathcal{H}_0'}, \Res_{\h} \rho^i_{\mathcal{H}_0'})
= 0$ for $i \neq 1$. In, particular, letting $J$ be an irreducible
component of $\Res_{\h} \rho^1$,
we get $\dim \Hom_{\h} (J, \Res_{\h} \rho_{\mathcal{H}'}) = a_1 \leq r$.
\end{proof}

\subsection{Compact form}

In this section, for $\k = \C$, we investigate the compact real form
of the reductive Lie algebra $\mathcal{H}_W$.

Let $\mathcal{H}_W^{\R} \subset \R W$ denote
the infinitesimal Hecke algebra defined over $\R$,
and $\mathcal{H}_W^{c}$ the \emph{real} Lie subalgebra of $\C W$
generated by the elements $\ii s$, for $s \in \mathcal{R}$. 
Note that conjugation by $s \in \mathcal{R}$ can be
written $\Id + (1/2) \ad( \ii s)^2$, hence $\mathcal{H}_W^c$
is also generated by the elements $\ii s$ for $s \in \mathcal{S}$,
for $\mathcal{S} \subset \mathcal{R}$ as above. In particular,
if $W$ is a Coxeter group, it is also generated by the $\ii s$ for $s$
running among the simple reflections w.r.t. a chosen Weyl chamber.

We let $W^{\pm} = \{ w \in W \ | \ \eps(W) = \pm w \}$. Clearly,
$W = W^+ \sqcup W^-$ and $\R W = \R W^+ \oplus \R W^-$.

\begin{prop} \label{compactform} $\mathcal{H}_W^c$ is a compact real form of $\mathcal{H}_W$,
and
$$
\mathcal{H}_W^c = \left( \mathcal{H}_W^{\R} \cap \R W^+ \right)
\oplus \ii \left( \mathcal{H}_W^{\R} \cap \R W^- \right).
$$
\end{prop}
\begin{proof} Let $\mathcal{H}^{c'}_W = ( \mathcal{H}_W^{\R} \cap \R W^+ )
\oplus \ii ( \mathcal{H}_W^{\R} \cap \R W^-)$, and let $\mathcal{L}$ be
the real linear span of the $w - \eps(w) w^{-1}, w \in W$ inside
$\R W$. It is easily checked that $\mathcal{L}$ is a Lie subalgebra of $\R W$.
One clearly has $\mathcal{L} = (\mathcal{L} \cap \R W^+) \oplus
(\mathcal{L} \cap \R W^-)$. The Lie algebra $\mathcal{H}_W$ is
spanned over $\C$ by elements of the form $[s_1,[s_2,\dots,[s_{r-1},s_r]\dots]$
for $s_i \in \mathcal{R}$. Such elements belong either to $\R W^+$
($r$ even) or $\R W^-$ ($r$ odd). We thus get
$\mathcal{H}_W^{\R} = (\mathcal{H}_W^{\R} \cap \R W^+) 
\oplus (\mathcal{H}_W^{\R} \cap \R W^-)$ and 
$\mathcal{H}_W^{\R} \otimes_{\R} \C = \mathcal{H}_W^{c'} \otimes_{\R} \C = \mathcal{H}_W$.
Moreover, $\mathcal{H}_W^c \subset \mathcal{H}_W^{c'}$. Since the $\ii s$ for
$s \in \mathcal{R}$ generate $\mathcal{H}_W$ over $\C$, we get equality of
dimensions hence $\mathcal{H}_W^c = \mathcal{H}_W^{c'}$. It remains to
show that this real Lie algebra is compact.

Consider the natural $W$-invariant quadratic form on $\R W$
given by $(w_1,w_2) = \delta_{w_1,w_2}$ if $w_1,w_2 \in W$.
We extend it to an hermitian form still denoted $(\ , \ )$ over
$\C W$. It defines a real bilinear form over $\C W$, for which
$\C W^+$ and $\C W^-$ are orthogonal (real) subspaces. Finally, we
restrict it to $\R W^+ \oplus \ii \R W^-$. It is easily checked
to be symmetric and positive definite. In order to prove that $\mathcal{H}_W^c$
is compact, we prove that this form is invariant under the adjoint action
of $\mathcal{H}_W^c \subset \R W^+ \oplus \ii \R W^-$.

Since $\mathcal{H}_W^c \subset (\mathcal{L} \cap \R W^+) \oplus \ii (\mathcal{L}
\cap \R W^-)$, it is enough to prove its invariance under the elements
of the Lie algebra $\C W$ of the form $w -w^{-1}$ and $\ii(w+w^{-1})$,
for $w \in W$. This calculation is straightforward and left to the reader.
\end{proof}

\section{Inductive properties of semisimple Lie algebras} 


Let $\h,\g$ be semisimple complex Lie algebras such that $\h \subset \g$,
and 
choose an irreducible \emph{faithful} representation
of $\g$. Equivalently, one may assume that $\g$
is a Lie subalgebra of $\sl(U)$, for some finite-dimensional complex
vector space $U$. Note that this implies $\rk \g \leq
\dim(U) - 1$, where $\rk \g$ denotes the semisimple rank of $\g$.

We will use the following lemma to show that $\g = \sl(U)$.

\begin{lemma} (see \cite{KRAMMINF}, proposition 3.8) \label{critsl} Let $U$
be a finite-dimensional complex vector space,
and $\g$ be a sub-Lie-algebra of $\sl(U)$ acting
irreductibly on $U$. If $\rk(\g) > \frac{\dim(U)}{2}$ then
$\g = \sl(U)$.
\end{lemma}

In order to identify $\g$ with an orthogonal or symplectic Lie
algebra, we will first use the following two lemmas in order to show that $\g$
is simple, and then use the classification of simple Lie algebras
and of their representations.

\begin{lemma} \label{lemredsimp} (see \cite{LIETRANSP}, lemme 15) Let $U$
be a finite-dimensional
$\C$-vector space, and
$\h \subset \g$ be two semisimple Lie subalgebras of $\sl(U)$.
Assume the following properties
\begin{enumerate}
\item $U$ is irreducible as a $\g$-module.
\item The restriction of $U$ to every simple ideal of $\h$ admits
an irreducible component of multiplicity 1.
\item One has $\rk(\g) < 2 \rk(\h)$.
\end{enumerate}
Then $\g$ is a simple Lie algebra.
\end{lemma}

In more subtle cases, we will use the following stronger form.

\begin{lemma} \label{lemredsimpstrong} Let $U$ be a finite-dimensional
$\C$-vector space, and
$\h \subset \g$ be two semisimple Lie subalgebras of $\sl(U)$.
We let $r = \rk \h$, and $\h = \prod_{j \in J} \h_j$
the decomposition of $\h$ as a product of its simple ideals.
For all
$j \in J$, we let
$$ m_j = \mathrm{gcd} \{ \dim \Hom_{\h_j} (R,\Res_{\h_j} U) \ \mid \
R \in \Irr(\h_j) \}, r_j = \rk \h_j
$$
If $U$ is an irreducible $\g$-module and $\dim(U) < (r+1)^2$, then
$\g$ is a simple Lie algebra if either
$$
(I)  \forall j \in J \ \ m_j = 1 \ \ \ \ 
\mbox{ or }
(II) \left\lbrace \begin{array}{lc}
1) & \forall j \in J \ \ m_j \leq 2 \\
2) & \g \ \mbox{does not have a simple ideal of rank 1}
\end{array} \right.
$$
\end{lemma}
\begin{proof}
We decompose $\g = \prod_{i \in I} \g_i$ as a product of simple
ideals. We let
$q_i$ be the composite of $\h \into \g \onto \g_i$, and
we let $\g^{(i)} = \prod_{j \neq i} \g_j$. We say that $\g_i$
is \emph{special} if $\Ker q_i \neq \{ 0 \}$.

Assume that there are no special ideal. Then
$\h$ embeds into every $\g_i$. By decomposing the irreducible $\g$-module
$U$ as a tensor product of irreducible $\g_i$-modules $U_i$, we get
$\dim U_i \geq \rk \g_i + 1 \geq r+1$ hence
$$
\dim U = \prod_{i \in I} \dim U_i \geq (r+1)^{\# I }
\Rightarrow \# I \leq \frac{\log \dim U}{\log (r+1)} <
\frac{ \log (r+1)^2}{\log (r+1)} = 2.
$$
It follows that $\# I = 1$, which shows that $\g$ is
a simple Lie algebra.
We are thus reduced to show that $\g$ has no special ideal.

If $\g_i$ is such a special ideal, then 
$\Ker q_i$ contains one of the simple ideals $\h_j$ of $\h$. Since
$U$ is an irreducible $\g$-module, it can be decomposed as
$U \simeq U_i \otimes U^{(i)}$ with $U_i$ and $U^{(i)}$ irreducible
faithful representations of $\g_i$ and $\g^{(i)}$, respectively.
In particular, $\dim U_i \geq 2$. Moreover,the image of $\h_j$ in $\g_i$
is $\{ 0 \}$, thus $\h_j$ act trivially on $U_i$ and
$
\Res_{\h_j} U = (\dim U_i) \Res_{\h_j} U^{(i)}.
$
It follows that $\dim U_i$ divides every
$\dim \Hom_{\h_j} (R,\Res_{\h_j} U)$
for $R$ an irreducible $\h_j$-module, hence their gcd $m_j$.

Under assumption $(I)$, this implies $\dim U_i = 1$,
a contradiction. Hence $(I)$ implies that there are no special
ideal. We now assume $(II)$. We have $\dim U_i = 2$. But
%
%
only rank 1 simple Lie algebras admit an irreducible
representations of dimension 2, so we get a contradiction.

%
\end{proof}

Our conventions on the numbering of heighest weights are the ones of \cite{FH}.
\begin{lemma} \label{caractBCD} (see \cite{LIETRANSP}, lemme 14) The couples $(\g,\varpi)$ with $\g$ 
a rank $n$ simple complex Lie
algebra and $\varpi$ the highest weight of a $N$-dimensional irreducible representation
of $\g$ such that $2n \leq N < 4n$ are,
for all $n \geq 2$, $(B_n,\varpi_1)$ if $N = 2n+1$, $(C_n,\varpi_1)$
if $N = 2n$,
and, for all $n \geq 3$, $(D_n,\varpi_1)$ if $N = 2n$ plus,
for $n \leq 6$,
the following exceptional couples.
\begin{itemize}
\item[Rank 1 :] $(A_1,\varpi_1)$ for $N = 2$, $(A_1,2\varpi_1)$ for
$N = 3$.
\item[Rank 2 :] $(G_2,\varpi_1)$ for $N = 7$, $(B_2,\varpi_2)$ for
$N = 4$, $(C_2,\varpi_1)$ for $N = 5$.
\item[Rank 3 :] $(B_3,\varpi_3)$ for $N = 8$, $(A_3,\varpi_2)$ for
$N = 6$.
\item[Rank 4 :] $(D_4,\varpi_3)$ and $(D_4,\varpi_4)$ for $N = 8$, $(A_4,
\varpi_2)$ and $(A_4,\varpi_3)$ for
$N = 10$.
\item[Rank 5 :] $(A_5,\varpi_2)$ and $(A_5,\varpi_4)$ for $N = 15$,
$(D_5,\varpi_4)$
and $(D_5,\varpi_5)$ for
$N = 16$.
\item[Rank 6 :] $(A_6,\varpi_2)$ and $(A_6,\varpi_5)$ for
$N = 21$.
\end{itemize}
\end{lemma}

For further use, we note the following facts in the
above list of exceptions. Starting from rank 4 we have $N \geq 8$, and the representations
are not selfdual except for $(D_4,\varpi_3)$ and $(D_4,\varpi_4)$ : their dual is a different exception also
noted in the list. On the contrary, before rank 4 we have $N \leq 8$
and all representations are selfdual. In particular, if the representation
is self dual then the rank is at most 4 and $N \leq 8$, otherwise
the rank is at least 4 and $N \geq 10$.

Moreover, note that some of the exceptions noted above
are simply reminiscent from the well-known exceptional
isomorphisms $A_3 \simeq D_3$, $B_2 \simeq C_2$.

When we know in advance that the rank of $\g$ is large enough,
we will use the following lemma.

\begin{lemma} \label{lemdim5}
Let $\g$ be a simple (complex) Lie subalgebra of $\sl(U)$ which
acts irreductibly on $U$, and such that $\rk \g > \frac{\dim U}{5}$. Then
\begin{enumerate}
\item if $\rk \g  \geq 10$ and $U \not\simeq U^*$, then $\g = \sl(U)$;
\item if $\rk \g  \geq 6$ and $U \simeq U^*$, then $\g \simeq \sp(U)$
or $\g \simeq \so(U)$.
\end{enumerate}
\end{lemma}
\begin{proof} It is easily checked that, for $ n \geq 11$,
$\sl_n(\C)$ admits no other irreducible representation of dimension
less than $5(n+1)$ than $\C^n$ and its dual. Similarly, only the 
$n$-dimensional irreducible representations of the Lie algebras of Cartan
type $B_n,C_n,D_n$ have dimension less than $5n$, provided that $n \geq 6$ ;
moreover, these representations are selfdual. This proves (1),
as there are no exceptional simple Lie algebra of rank at least 10.
Moreover, the irreducible representations of $\sl_n(\C)$, $ 6 \leq n \leq 9$,
which have dimension at most $5(n+1)$, are easily checked not to be
selfdual. In order to prove (2) we thus only need to exclude the
exceptional types $E_6,E_7,E_8$. The representations of $E_7$ and $E_8$
corresponding to the fundamental weights have dimension at least $56 > 5\times 8$,
so these are excluded. The simple Lie algebra of type $E_6$ admits exactly
two representations in this range, which are dual one to the
other, and thus do not satisfy the assumptions of (2).
\end{proof}

In small rank, we will need a specific result for 6-dimensional representations.

\begin{lemma} \label{lemliedim6} Let $\g \to \sl(U)$ be a
faithful 6-dimensional irreducible representation
of a semisimple Lie algebra $\g$.
Assume that there
exists $\h \subset \g$ with $\h \simeq \sl_2$ such that
$\Res_{\h}
U$ admits multiplicities and also an irreducible
component occuring with multiplicity 1. Then $\g$ is simple. Moreover,
$\g \to \sl(U)$ is not selfdual iff $\g \simeq \sl(U)$ and otherwise
either $\g \simeq \so(U)$ or $\g \simeq \sp(U)$.
\end{lemma}

\begin{proof}
Let $\g = \g^1 \times \dots \times \g^p$ be a decomposition of
$\g$ in simple factors. Then $U$ can be decomposed as
$U_1 \otimes \dots \otimes U_p$ with 
$\g^i$ acting faithfully and irreducibly on $U_i$.
Since $\g^i$ is simple we have $\dim U_i \geq 2$
and get $\dim U \geq 2^p$ hence $p = 2$ if $\g$ is not simple.
In that case, since $\dim U = 6$ we can assume $\dim U_1 = 2$
and $\dim U_2 = 3$. It follows that $\g^1 \simeq \sl_2$. Let $\pi_i :
\g \onto \g^i$ and $\iota : \h \into \g$ denote the canonical
morphisms, and let $\varphi_i = \pi_i \circ \iota$. If $\varphi_1 = 0$,
we would have $\h \subset \g^2$ and
$$\Res_{\h} U =
\Res_{\h} \Res_{\g^2} U = \Res_{\h} (\dim U_1) U_2 = (\dim U_1)
\Res_{\h} U_2$$
contradicting the existence of a multiplicity 1 component. It follows
that $\varphi_1 \neq 0$ and similarly we show $\varphi_2 \neq 0$.
Since $\h \simeq \sl_2$ is simple it follows that $\varphi_1$ and
$\varphi_2$ are injective.

We let $[p]$ denote the $(p+1)$-dimensional irreducible
representation of $\h \simeq \sl_2$. Since $\Res_{\h} U_1$ and
$\Res_{\h} U_2$ are
faithful we have $\Res_{\h} U_1 = [1]$ and
$\Res_{\h} U_2\in \{[2], [0]+ [1]\}$, hence
$$
\Res_{\h} U_1\otimes U_2 \in \{ [1] \otimes [2], [1] \otimes
\left( [0] + [1] \right) \} = \{ [1] + [3], [1] +[0]+[2] \}
$$
contradicting the presence of multiplicities in $\Res_{\h} \rho$.

By contradiction, it follows that $\g$ is simple. Since $\g \subset
\sl_6$ we have $1 \leq \rk \g \leq 5$. If $\rk \g = 1$ we would
have $\g = \h$ contradicting either the irreducibility of the action of
$\g$ or the presence of multiplicities in $\Res_{\h} U$. If $\rk \g > 3$
then $\g = \sl(U)$ by lemma \ref{critsl} and we get the
conclusion.

We can thus assume $\rk \g \in \{ 2, 3 \}$. If $\rk \g = 3$,
$\g$ cannot be of type $B_3$ because $\so_7$ does not admit
an irreducible representation of dimension 6, hence is either of
Cartan type $A_3 \simeq D_3$ or $C_3$. The only possibilities for
$U$ to be 6-dimensional and irreducible imply $U \simeq U^*$
as $\g$-modules and  $\g \simeq \so(U)$ or $\g \simeq \sp(U)$.

There remains to rule out the case $\rk \g = 2$. The cases $G_2$
and $C_2 = B_2$ are excluded because they do not admit
6-dimensional irreducible representations. We thus have $\g \simeq \sl_3$.
Now an embedding $\sl_2 \into \sl_3$ corresponds to a 3-dimensional
faithful representation of $\sl_2$, either $[2]$ or $[0]+[1]$.
Letting $U_0$ denote one of the standard (3-dimensional) representations
of $\sl_3\simeq \g$ we have $\Res_{\h} U_0 \in \{ [2],  [0]+[1] \}$.
On the other hand, 6-dimensional irreducible representations
of $\sl_3$ are isomorphic either to $S^2 U_0$ or its dual. But
$S^2 [2] = [4]+[0]$ and
$$
S^2 \left( [0] + [1] \right) = S^2[0] + [0]\otimes [1] + S^2 [1]
= [0] + [1] + [2]
$$
are both multiplicity-free, contradicting our assumption. This
concludes the proof of the lemma.
\end{proof}

\section{Representations of $G(e,1,r),G(e,e,r)$ and $G(2e,e,r)$} 

Recall from the Shephard-Todd classification that an irreducible
(finite) \emph{pseudo-}reflection group $W$ of rank at least 2 belongs
either to a finite list of 34 exceptions, labelled
from $G_4$ to $G_{37}$, or to an infinite series $G(de,e,r)$
with 3 integral parameters with $(d,e,r) \neq (1,2,2)$. The group $G(de,e,r)$ is defined
to be the set of $r \times r$ monomial complex matrices with entries
in $\mu_{de} = \mu_{de}(\C)$ whose product of nonzero entries belongs
to $\mu_d$. It has order $r! d^r e^{r-1}$, and rank $r$ if $(d,e) \neq (1,1)$.

The reflection groups in this infinite series are the groups $G(e,e,r)$
and $G(2e,e,r)$. The pseudo-reflection groups $G(e,1,r)$ are used
as auxiliary tools in the construction of their representations, as
$G(e,e,r)$ is a normal index $e$ subgroup of $G(e,1,r)$. Similarly,
$G(2e,2e,r)$ is an index $2$ subgroup of $G(2e,e,r)$, which has
index $e$ in $G(2e,1,r)$. The number of conjugacy classes is given by
the following lemma, which is standard and easy to check.

\begin{lemma} \label{lemgeersingle} If $W$ is an irreducible reflection group
of type $G(e,e,r)$ then $\# \mathcal{R}/W \leq 2$, and
$\# \mathcal{R}/W = 1$ for $r \geq 3$. If it has type
$G(2e,e,r)$ then $\# \mathcal{R}/W \leq 3$ and
$\# \mathcal{R}/W \leq 2$ if $r \geq 3$.
\end{lemma}

\subsection{Preliminaries}

Let $r \geq 1$. As usual, we label irreducible representations
of the symmetric group $\mathfrak{S}_r$ by partitions
$\la = [\la_1,\dots,\la_s]$ with $\la_i \in \Z_{>0}$ and $\la_i \geq \la_{i+1}$
of total size $|\la| = \sum \la_i = r$, choosing for convention
that $[r]$ labels the trivial representation of $\mathfrak{S}_r$.
We let $\la'$ denote
the conjugate partition of $\la$, defined by $\la'_i = \# \{ j \ ; \ \la_j \geq i \}$.
We denote $\emptyset$ the only partition of $r$, identify when needed
a partition with its Young diagram, using the convention that
the diagram $[3,2]$ has two rows and three columns. When convenient, we also
identify partitions with the irreducible representation
of the adequate symmetric group labelled by it.

Letting $\la,\mu$ be two partitions, we use the notation $\mu \subset \la$
if $\forall i \ \mu_i \leq \la_i$ and $\mu \nearrow \la$ for
$\mu \subset \la$ and $|\la| = |\mu| + 1$. In terms of Young diagrams,
this means that $\mu$ is deduced from $\la$ by removing one box.
Young's rule states that, under the usual inclusion $\mathfrak{S}_r \subset
\mathfrak{S}_{r+1}$, the restriction of an irreducible representation
$\la$ of $\mathfrak{S}_{r+1}$ to $\mathfrak{S}_r$ is the direct sum
of all $\mu \nearrow \la$, and in particular the number of irreducible
components of this restriction is $\delta(\la) = \# \{ i \ | \ \la_i > \la_{i+1} \}$.

Let then $e \geq 1$. The group is isomorphic to the wreath product
$\Z/e\Z \wr \mathfrak{S}_r = \mathfrak{S}_r \ltimes (\Z/e\Z)^r$, whence
its irreducible representations are indexed by multipartitions
$\bla = (\la^0,\dots,\la^{e-1})$ of total size
$|\bla| = |\la^0|+\dots+|\la^{e-1}|=r$. We let $p(\bla) = \# \{ i \ | \ 0 \leq i 
\leq e-1, \la^i \neq \emptyset \}$, and define $\delta(\bla) = \sum \delta(\la^i)$.
It is the number of irreducible compoents in the (multiplicity-free)
restriction of $\bla$ to its natural subgroup $G(e,e,r-1)$, if $|\bla| = r$.
Indeed, these components corresponds to multipartitions $\bmu$ of
$r-1$ such that $\mu^i = \la^i$ for all but one $i$, for which $\mu^i \nearrow
\la^i$. We use the notation $\bmu \nearrow \bla$ in that case.

The cyclic group $\Z/e\Z$ acts on the set of irreducibles of $G(e,1,r)$
by cyclically permuting the parts of a multipartion $\bla$, and so
does its subgroup $2 \Z/e\Z$ when $e$ is an even integer. We denote
$\Aut(\bla)$ the fixer of $\bla$ under the former action, and $\Aut^0(\bla)$
the fixer under the latter. Let $A(\bla)$ and $B(\bla)$ denote
the orders of $\Aut(\bla)$ and $\Aut^0(\bla)$, respectively.
By Clifford theory, $A(\bla)$ is the number of irreducible components
of the restriction to $G(e,e,r)$ of the irreducible representation
$\bla$ of $G(e,1,r)$. Similarly, if $\bla$ is an irreducible representation
of $G(2e,1,r)$, then $B(\bla)$ is the number of irreducible
components of its restriction to $G(2e,e,r)$. Starting from $\bla \in
\Irr(G(2e,1,r))$, we clearly have $A(\bla)/B(\bla) \in \{ 1, 2 \}$.
Moreover, if $\rho$ is an irreducible component of the restriction
of $\bla$ to $G(2e,e,r)$, then $B(\bla) = A(\bla)$ if and only if
$\rho$ has irreducible restriction to $G(2e,2e,r)$, and otherwise
has two distinct components.

If $\rho$ is an irreducible representation of $G(e,e,r)$, it is
an irreducible component of the restriction of some
$\bla = (\la^0,\dots,\la^{e-1})$ of $G(e,1,r)$. We leave to the
reader to check that $\rho^* \otimes \eps$ is then
an irreducible component of $((\la^{e-1})',\dots,(\la^0)')$.

\subsection{Basic facts}

We fix $e \geq 1$. In order to make estimates of dimensions, we will
make special use of \emph{binary representations} associated to a
subset $I$ of $\{ 0,1,\dots,e-1 \}$. The corresponding
representation $b(I) = \bla$ of $G(e,1,\# I)$ is defined
by $\la^i = [1]$ is $i \in I$ and $\la^i = \emptyset$ otherwise.
Letting $r = \# I$, the restriction of $b(I)$ to $G(e,1,r-1)$ is
the sum of all $b(I')$ for $I' \subset I$ with $\# I' = r-1$. By
induction we get $\dim b(I) = r! = (\# I)!$. We define a partial
order $\bmu \subset \bla$ on multipartitions by $\mu^i \subset \la^i$
for all $0 \leq i \leq i-1$.

From the branching rule, the following three facts are clear. Removing
the empty parts of $\bla$, we get a representation of $G(p(\bla),1,|\bla|)$
of the same dimension ; also, the dimension increases with respect
to $\subset$ ; finally, since any $\bla$ contains a (unique) binary
multipartition $b(I)$ with $\# I = p(\bla)$, we get
$\dim \bla \geq p(\bla) !$ by restriction to $G(e,1,p(\bla))$,
with equality only if $\bla$ itself is a binary multipartition.
Similarly, if $0 \leq i \leq e-1$ is such that $\la^i \neq
\emptyset$, then the restriction to $G(e,1,|\la^i|)$
of $\bla$ contains at least $(p(\bla)-1)!$ copies of $\bmu$
such that $\mu^i = \la^i$ and $\mu^j = \emptyset$ if $j \neq i$.
In particular $\dim (\bla) \geq (p(\bla)-1)! \dim (\la^i)$.

In order to deal with the automorphism group of $\bla$, we will
need the following lemma. 

\begin{lemma} \label{brisuregeer} Let $\bla$ be a irreducible
representation of $G(e,1,r)$.
\begin{enumerate}
\item $A(\bla)$ divides $e$, $r$ and $p(\bla)$.
\item If $A(\bla) \neq \{ 1 \}$ and $\bmu \nearrow \bla$ then
$A(\bmu) = \{ 1 \}$
\end{enumerate}
\end{lemma}
\begin{proof}
$A(\bla)$ divides $e$ because $\Aut(\bla)$ is a subgroup of $\Z/e\Z$.
Let $A(\bla) = e/b$. Then $\la^{i+b} = \la^i$ hence
$r = A(\bla) \sum_{i=0}^{b-1} |\la^i|$ and $p(\bla) = A(\bla)
\# \{ i \ | \ 0 \leq i |eq b-1, \la^i \neq \emptyset \}$, whence
the conclusion of (1).
Part (2) is proved in
\cite{BRANCH}, prop. 3.1.
\end{proof}

As a consequence, we get the following rough estimates on the dimensions
of irreducible representations of $G(e,e,r)$.

\begin{lemma} \label{lemineqsmall} Let $\rho$ be an irreducible component of the restriction
to $G(e,e,r)$ of a representation $\bla$ of $G(e,1,r)$.
Then $\dim \rho \geq (p(\bla)-1) !$ and
$\dim \rho \geq p(\bla) !$ as soon as $\exists i \ \la^i \not\in \{ \emptyset, [1] \}$.
\end{lemma}
\begin{proof}
If $A(\bla) = 1$ we have $\dim \rho = \dim \bla \geq p(\bla)!$ and we
are done. We thus can assume $A(\bla) \neq 1$.

By lemma \ref{brisuregeer} for
all $\bmu \nearrow \bla$ we have $A(\bmu) = 1$. It follows
that the restriction to $G(e,e,r-1)$ of $\rho$ contains the
restriction of all these $\bmu \nearrow \bla$, for which
$\dim(\bmu) \geq p(\bmu)!$. Moreover,
$p(\bmu) = p(\bla)$ unless 
$\bmu$ is deduced from $\bla$ by removing a box at some place
$i$ with $\la^i = [1]$, in which case $p(\bmu) = p(\bla) -1$. The conclusion follows.
\end{proof}

Finally, the following is
lemma is standard and easily proved.

\begin{lemma} \label{lemdim2rang2}If $W$ has type $G(e,e,2)$ or $G(2e,e,2)$,
then $\dim \rho \leq 2$ for all $\rho \in \Irr(W)$.
In particular $\Irr(W) = \Lambda \Ref(W)$.
\end{lemma}

\subsection{Quasi-reflection representations of $G(e,e,r), \ r \geq 3$}

Recall from lemma \ref{lemgeersingle} that $G(e,e,r)$ admits only one conjugacy class of reflections
for $n \geq 3$.
The goal of this section is to determine the representations
in $\LRef$ for the groups $G(e,e,r)$.

We first recall that the groups $G(e,e,r)$ admit a distinguish
set of generating reflections $s'_1,s_1,s_2,\dots,s_r$, where
$s_k$ is the permutation matrix corresponding to the transposition
$(k, k+1)$ and $s'_1$ is the conjugate of $s_1$ by the diagonal
matrix $\mathrm{diag}(\exp(\frac{2 \ii \pi}{e}),1,\dots,1)$.
If $e_1$ divides $e_2$, there exists
a natural morphism from $G(e_2,e_2,r)$ to $G(e_1,e_1,r)$
that maps each generator to the generator with the same
name. In matrix terms it is deduced from the ring homomorphism
of $\Z[\zeta]$ that maps $\zeta$ to $\zeta^{\frac{e_2}{e_1}}$, where
$\zeta = \exp(2 \ii \pi/e_2)$. If $\bla$ is an
irreducible representation of $G(e_1,1,r)$ which is restricted
to $G(e_1,e_1,r)$, then the representation of $G(e_2,e_2,r)$
deduced from this morphism is the restriction of $\bmu
\in \Irr G(e_2,1,r)$ described by $\mu^{\frac{e_2}{e_1} k} = \la^k$
for all $k$, and
$\mu^k = \emptyset$ if $k$ is not divisible by $e_2/e_1$. Conversely,
the restriction of $\bmu \in \Irr G(e_2,1,r)$ factorizes
through $G(e_1,e_1,r)$ iff $j-i$ is divisible by $e_2/e_1$
whenever $\mu^i, \mu^j \neq \emptyset$. All this is
easily checked from the explicit formulas of
\cite{ARIKI} or \cite{MARINMICHEL} (see formulas (3.2) there).
Obviously, through these natural morphisms, elements of $\LRef$
for $G(e_1,e_1,r)$ induce elements of $\LRef$ for
$G(e_2,e_2,r)$.

We prove the following.

\begin{prop} The quasi-reflection representations $\rho$ of
$G(e,e,r)$ for $r \geq 3$ are given in table \ref{tabqrefrep},
where $\bla$ is an irreducible representation of $G(e,1,r)$
whose restriction to $G(e,e,r)$ contains $\rho$.
\end{prop}

\begin{table}
$$
\begin{array}{|c|c|c|c|c|l|l|l|}
\hline
e & r & \dim \rho &  A(\bla) & p(\bla) & \la^{i_1}\mbox{ or }(\la^{i_1})' & \la^{i_2}\mbox{ or }(\la^{i_2})' & \la^{i_3}\mbox{ or }(\la^{i_3})' \\
\hline
e \geq 1 & r \geq 3 & r-1 & 1 & 1 & [r-1,1] &  &  \\
\hline
e \geq 2 & r \geq 3 & r & 1 & 2 & [r-1] & [1] &  \\
\hline
e \geq 2 & 4 & 2 & 1 & 1 & [2,2] &  &  \\
\hline
3|e & 3 & 2 & 3 & 3 & [1] & [1] & [1] \\
\hline
2|e  & 4 & 3 & 2 & 2 & [2] & [2] &  \\
\hline
\end{array}
$$
\caption{Quasi-reflection representations of $G(e,e,r)$ for $r \geq 3$.}
\label{tabqrefrep}
\end{table}

\begin{proof}
If $\rho$ is a representation of $G(e,e,r)$, we let $\aa(\rho)$ and $\bb(\rho)$
denote the multiplicity of $1$ and $-1$, respectively, in the spectrum of a
reflection. By abuse of notation, if $\bla$ is a representation of
$G(e,1,r)$, we denote $\aa(\bla)$ and $\bb(\bla)$ the multiplicities
for the restriction to $G(e,e,r)$.

First note that $\aa(\bla)$ and $\bb(\bla)$ are increasing functions of $\bla$
for the order $\subset$. Also note that they only depend on the collection
of non-empty partitions (with multiplicities) that form $\bla$.
Finally, note that, if $\rho_1,\rho_2$ are two components of the
same $<\bla>$, then it readily follows from the formulas 3.3 of \cite{MARINMICHEL}
that $\aa(\rho_1) = \aa(\rho_2)$ and $\bb(\rho_1) = \bb(\rho_2)$.

We assume that $\rho$ is a quasi-reflection representation.
If $\rho = < \bla >$ is irreducible then this case has been dealt with in
\cite{MARINMICHEL} -- specifically, the arguments in \S 5 there justify the
first three lines of table \ref{tabqrefrep}. Otherwise, $\rho$ embeds in some representation
$\bla$ of $G(e,1,r)$ with $a = A(\bla) \neq 1$. We let $\alpha = e/a$,
hence $\bla = (\la^0, \dots, \la^{\alpha-1} , \la^0, \dots)$.

Assume that there exists distincts $i,j \in [0,\alpha-1]$ such
that $\la^i \neq \emptyset$ and $\la^j \neq \emptyset$.
We refer to \cite{ARIKIKOIKE} or \cite{MARINMICHEL} for the
definitions of multitableaux and standard multitableaux of shape $\bla$,
and recall that the representation $\bla$ of $G(e,1,r)$
has a basis indexed by all standard multitableaux of shape $\bla$. We also
recall that $s_1$ acts by $1$ (resp. $-1$) on such a standard
multitableau $\bT$ if $1$ and $2$ are placed on the same line
(resp. the same column) of some part of $\bT$, and that otherwise
the multitableau $\bT'$ deduced from $\bT$ by exchanging 1 and 2
is also standard, the plane spanned by $\bT$ and $\bT'$ is stable
under $s_1$, which acts with eigenvalues $1, -1$ on this plane.

We define
$2\alpha$ standard multitableaux $\bT^{u,\pm}$ of shape $\bla$,
for $u \in [0,a-1]$, in the following way. For $\bT^{u,+}$, we
place $1$ in position $i + \alpha u$ and $2$ in position $j + \alpha u$.
For $\bT^{u,-}$, we
place $1$ in position $j + \alpha u$ and $2$ in position $i + \alpha u$.
Then we fill in the remaining boxes, with numbers from $3$ to $r$,
in a uniform way. The subspace of $\bla$ generated by these is stable
under $s_1 \in G(e,1,r)$, and $1,-1$ are eigenvalues with multiplicity
$a$. Since $a \geq 2$ we can also choose multitableaux $\bT^+$ and $\bT^-$
such that $1$ lies in position $i$ and $2$ in position $i+\alpha$
for $\bT^+$ and in the reverse order for $\bT^-$ -- the filling for
the remaining numbers being the same. The action of $s_1$ on the plane
generated by $\bT^+$ and $\bT^-$ has two eigenvalues, $1$ and $-1$.
Since this plane is in direct sum with the precedingly defined subspace,
it follows that $\aa(\bla) \geq a+1$ and $\bb(\bla) \geq a+1$.
In particular $\aa(\rho) = \aa(\bla)/a \geq 1 + \frac{1}{a}$
and similarly $\bb(\rho) \geq 1 + \frac{1}{a}$, contradicting
the assumption that $\rho$ is a quasireflection representation.

We thus can assume that only one $i \in [0,\alpha-1]$ satifies
$\la^i = \la \neq \emptyset$. Assume $\la \supset [2,1]$ and 
consider $2 \alpha$ multitableaux $\bT^{u,\pm}$
for $u \in [0,a-1]$ defined by placing $1$ and $2$ in position
$i + \alpha u$, the number $2$ being placed above or on the right of $1$
depending on $\pm$, the filling for the remaining numbers
being uniform. As before, we get $1$ and $-1$ as eigenvalues for the action
of $s_1$ on this subspace with
multiplicity $a$. Considering the plane generated by $\bT^+$ and $\bT^-$
defined in the same way, we get the same contradiction.

It follows that, up to tensorization by the sign character, we can assume
$\la = [m]$ for some $m \geq 1$. In that case there are two types
of multitableaux $\bT$ of shape $\bla$ :
\begin{enumerate}
\item $1$ and $2$ are placed in the same position $i + \alpha u$ for some
$u \in [0,a-1]$
\item $1$ and $2$ are placed in positions $i + \alpha u$ and
$i + \alpha u'$ with $u,u' \in [0,a-1]$ and $u \neq u'$.
\end{enumerate}
Multitableaux of type (1) are fixed par $s_1$. The multiplicities of
$1$ and $-1$ for the action of $s_1$ on the subspace generated by all
multitableaux of type (2) are the same. If $m \geq 2$, define
$\bla'$ by replacing $[m]$ with $[m-2]$ in position $i$ of $\bla$.
If $m = 1$ we note $\dim \bla' = 0$. Define $\bla''$ by replacing
$[m]$ with $[m-1]$ in position $i$ and $i + \alpha$ in $\bla$.
The number of multitableaux of type (1) is $a \dim \bla'$, and the
number of multitableaux of type (2) is $a(a-1) \dim \bla''$. It
follows that $\aa(\rho) = \frac{a-1}{2}
\dim \bla'' +  \dim \bla'$ and $\bb(\rho) = \frac{a-1}{2}
\dim \bla''$. In particular $\aa(\rho) \geq \bb(\rho)$ and
$\rho$ is a quasireflection representation iff $\bb(\rho)=1$.

If $m = 1$, then $\bla''$ is a binary representation of size $a-2$,
hence $(a-1) \dim \bla = (a-1)!$ and $\bb(\rho) =1$ iff $(a-1)! = 2$,
that is $a = 3$ hence $r = 3$. It follows that $\rho$ factorizes
through one of the (2-dimensional) irreducible components
of the restriction to $G(3,3,3)$ of $([1],[1],[1])$.
If $m=2$, then $p(\bla'') = p(\bla) = a$ hence $(a-1)! \dim \bla'' \geq
a! (a-1) = 2$ iff $a = 2$ and $\bla''$ is a binary representation,
whence $m=2,r=4$ and $\rho$ factorizes through one of the (3-dimensional)
components of the restriction of $([2],[2])$ to $G(2,2,4)$.

\end{proof}

We are now able to determine $\LRef$.
\begin{prop} \label{listeLRef}
If $\rho$ is an irreducible representation
of $G(e,e,r)$ for $r \geq 3$ with $\dim \rho > 1$ which
belongs to $\LRef \setminus \QRef$, then
$\rho$ is the restriction of an irreducible representation
$\bla$ of $G(e,1,r)$ and, either $p(\bla) = 1$
and there exists $\la^i = [r-p, 1^p]$ for some $p \in [2,r-3]$,
or $p(\bla) = 2$ and there exist $i \neq j$ with $\la^i = [r-p]$, $\la^j = [1^p]$
for some $p \in [2,r-2]$.
\end{prop}
\begin{proof}
Let $R$ be a reflection representation of $G(e,e,r)$ and $\rho = \eta \otimes
\Lambda^k R \in \LRef$. Assuming $\rho \not\in \mathrm{QRef}$ and $\dim \rho > 1$,
this implies $\dim R \geq 4$. By table \ref{tabqrefrep} it follows
that $R$ is the restriction of some representation $\bla$ of $G(e,1,r)$
with $A(\bla) = 1$, $p(\bla) \leq 2$ and
$\dim \bla \in \{ r-1,r \}$, which implies $r \geq 4$.
We can assume $\la^0 \neq \emptyset$. There are two cases to consider. If
$p(\bla) = 1$, then $R$ factors through the representation of $G(1,1,r)
=\mathfrak{S}_r$ corresponding to the partition $\la^0 = [r-1,1]$, and
it is well-known that $\Lambda^k [r-1,1] = [r-k,1^k]$ as representation
of $\mathfrak{S}_r$, which proves half of the proposition. If
$p(\bla) = 2$, we can assume $\la^0 = [r-1]$ and $\la^s = [1]$
for some $1 \leq s \leq e-1$.

We let $\bmu(k,r,s)$ denote the multipartition $\bmu$ with $|\bmu|= r$,
$p(\bmu) \leq 2$ such that $\mu^0 = [r-k]$, $\mu^s = [1^k]$. In
particular, $\bla = \bmu(1,r,s)$. We prove that, as a representation
of $G(e,1,r)$, $\Lambda^k \bmu(1,r,s) = \bmu(k,r,s)$ for $0 \leq k \leq r$.
This will prove the second part of the proposition. Since
$\bmu(0,r,s)$ is the trivial representation of $G(e,1,r)$, we can assume
$k \geq 1$. The proof is then by induction on $r$, the case $r = 2$
being easily checked from the matrix models of \cite{ARIKIKOIKE} or
\cite{MARINMICHEL}. We thus assume $r \geq 3$. If $k < r$, then
by the branching rule the restriction to $G(e,1,r-1)$ of $\bmu(1,r,s)$
is $\un + \bmu(1,r-1,s)$, hence the restriction of $\Lambda^k \bmu(1,r,s)$
is $\Lambda^k \bmu(1,r-1,s) + \Lambda^{k-1} \bmu(1,r-1,s)$,
which is $\bmu(k,r-1,s)+\bmu(k-1,r-1,s)$ by the induction hypothesis.
On the other hand, by Steinberg theorem we know that $\Lambda^k
\bmu(1,r,s)$ is irreducible hence corresponds to some multipartition
of size $r$ that contains both $\bmu(k,r-1,s)$ and $\bmu(k-1,r-1,s)$. The
only possibility being $\bmu(k,r,s)$ this concludes the proof,
the case $k=r$ being similar and left to the reader.
\end{proof}

\section{Induction process} 

In this section we assume that $W$ is an irreducible reflection
group and $\CW$ is a proper reflection subgroup of $W$,
meaning that $\CW$ is a proper subgroup of $W$ which is generated by
reflections of $W$,
%
for which the theorem holds. We let $\mathcal{R}_0 = \CW \cap
\mathcal{R}$ denote the set of reflections of $\CW$.

For convenience, for $\rho \in \Irr(W)$ we let $\overline{\rho}$
or $\Res \rho$ denote the restriction of $\rho$ to $\CW$. We denote $(\ | \ )$
the standard scalar product on the representation ring, and the notation
$\varphi \nearrow \rho$ means $(\overline{\rho}| \varphi) \geq 1$ for
$\rho \in \Irr(W)$ and $\varphi \in \Irr(W_0)$. Note that this notation
is consistent with the one introduced in the previous section for (multi-)partitions,
when $W = G(d,1,r)$ and $W_0 = G(d,1,r-1)$. We let $\LRef = \LRef(W)$ and
denote $\Ind \LRef = \{ \varphi \in \Irr(W) \ | \ \exists \varphi \in \LRef(W_0) \ 
\varphi \nearrow \rho \}$. If $\overline{\rho} = a_1 \rho_1 + \dots + a_r \rho_r$,
with $a_i \in \Z_{>0}$ and $\rho_i$ non-isomorphic irreducible representations
of $W_0$, we let $soc(\overline{\rho})$ denote $\rho_1 + \dots + \rho_r$. We
use the notation $\mathcal{H}'_0, \mathcal{H}_0(\varphi)$
for $\mathcal{H}'_{\CW}, \mathcal{H}_{\CW}(\varphi)$ and extend the
definition \ref{defHrho} of $\mathcal{H}_0(\varphi)$ to non-irreducible
representations so that, if $\overline{\rho}$ is decomposed as above, then
$\mathcal{H}_0(\overline{\rho}) = \mathcal{H}_0(soc(\bar{\rho})) =
\mathcal{H}_0(\rho_1) + \dots + \mathcal{H}_0(\rho_r) \subset \mathcal{H}'_0$.
We still have $\mathcal{H}_0(\bar{\rho}) \simeq \bar{\rho}(\mathcal{H}'_0)$,
and $\mathcal{H}_0(\bar{\rho})$ naturally embeds into $\mathcal{H}(\rho)$.

In \S 1 we prove several basic induction lemmas under the assumption
that $W_0$ admits a single class of reflections.
In \S 2 we use them to prove that, once we know how to prove the theorem
for any group $W$ with a single class of reflections, then we can
prove the general case. This enables us to assume
that, if $W$ belongs to the infinite series, then $W$ has type $G(e,e,r)$ with $r \geq 3$. In \S 3 we
show a few preliminary results for specific induction patterns
and in \S 4 we deal with the case $r = 3$, so that we can assume
$r \geq 4$ in the sequel. Then we put $W_0 = G(e,e,r-1)$ and assume
that the theorem holds for $W_0$. We prove under this assumption that
$\mathcal{H}(\rho) \simeq \mathcal{L}(\rho)$, first for
$\rho \in \Ind \LRef \setminus \LRef$ in \S 5, then for $\rho \not\in \Ind\LRef$
in \S 6. This implies that the theorem holds for $W$
by corollary \ref{corimpltheo}, and proves the theorem
for the groups $G(e,e,r)$ by induction on the rank.
Finally, \S 7 deals with the 15 exceptional groups.

\subsection{Induction results in the single case}

In the following lemmas we assume that $\CW$ admits only
one conjugacy class of reflections. A consequence of this assumption
is that non-isomorphic irreducible representations of $\CW$
of dimension at least 2 correspond to non-isomorphic
irreducible representations of
$\mathcal{H}_{\CW}$ by corollary \ref{corisorep}. Another
consequence is that $\CW$ is generated by the conjugacy class
in $W_0$ of any reflection $s \in W_0$.

\begin{lemma} \label{lemlrefind}Assume that $\CW$ admits a single
conjugacy class of reflections. Then $\LRef \subset \Ind \LRef$.
\end{lemma}
\begin{proof}
Let $\rho \in \LRef$ with $\dim \rho > 1$, the case
$\dim \rho = 1$ being obvious. We have $\rho = \eta \otimes \Lambda^k \rho_0$
for some $\rho_0 \in \Ref$, $\eta : W \to \{ \pm 1 \}$ and $1 \leq k <
\dim \rho_0$. The restriction
$\bar{\eta}$ of $\eta$ to $\CW$ is a morphism from $\CW$ to $\{ \pm 1 \}$.
We choose $s \in \mathcal{R}_0$ and let
$\bar{\rho_0}$ denote the restriction of $\rho_0$ to $\CW$.

Let us write $\bar{\rho_0} = \rho_1 \oplus \dots \oplus \rho_r$ with
$\rho_1,\dots, \rho_r$ irreducible. We can assume $\dim \rho_1 \geq \dim
\rho_i$ for $i \geq 1$. Since $\rho_0(s)$ is either a reflection or a scalar, there
exists at most one $i$ such that $\rho_i(s) \not\in \{ -1, 1 \}$.
Moreover, if $\rho_i(s) =  \pm 1$ then $\dim \rho_i = 1$
since $\CW$ is generated by the conjugacy class of $s$. It follows
that $\dim \rho_2 = \dim \rho_3 = \dots = \dim \rho_r=1$. We need to
find a subrepresentation of $\bar{\rho}$ which belongs to $\LRef(W_0)$.

We let $\C e_i$, $U$, and $U_0$
denote the underlying vector spaces of $\rho_i$ for $i \geq 2$,
$\rho_1$ and $\rho_0$, respectively.
If $k \leq r-1$ then the vector $e_2 \wedge \dots \wedge e_{k+1} \in \Lambda^k
U_0$ is non-zero, and $\CW$ acts on it by $\bar{\eta} \otimes (\rho_2 \otimes
\dots \otimes \rho_{k+1})$, which is 1-dimensional hence
belongs to $\LRef(W_0)$.
If $k \geq r$ then the subspace $\Lambda^{k+1-r} U \wedge e_2 \wedge \dots
\wedge e_r$ of $\Lambda^k U_0$ is non-zero ; indeed, we have $k+1-r \geq 1$
and $k \leq \dim U_0 = \dim U + r -1 \Rightarrow k+1-r \leq \dim U$.
Moreover, it is $\CW$-stable and $\CW$ acts on it by
$(\Lambda^{k+1-r} \rho_1) \otimes (\bar{\eta} \otimes \rho_2 \otimes \dots
\otimes \rho_r )$, which belongs to $\LRef(\CW)$. It follows
that $\rho \in \Ind \LRef$ in both cases.
\end{proof}

\begin{lemma} \label{lemnormsimp} Assume that $W_0$ admits
a single class of reflections. Let $\rho \in
\Irr(W) \setminus \Ind \LRef$. If $\forall \rho' \in \Irr(\CW) 
\ (\bar{\rho}|\rho') \leq 1$ then $\mathcal{H}(\rho)$ is a simple Lie
algebra. The same conclusion holds if  $\forall \rho' \in \Irr(\CW) 
\ (\bar{\rho}|\rho') \leq 2$ provided that
$W$ admits a single conjugacy class of reflections
and $\forall \rho \in \Irr(W) \ \dim \rho \leq 3 \Rightarrow
\dim \rho = 1$.
\end{lemma}
\begin{proof} Let $soc(\bar{\rho})$ denote the
(direct) sum of the irreducible components of $\bar{\rho}$ and
$r = \rk \mathcal{H}_0(\bar{\rho})$. By
assumption we have $\dim \rho \leq 2 \dim soc(\bar{\rho})$. Let
$\h$ denote one of the simple Lie ideals of $\mathcal{H}_0(\bar{\rho})$.
By hypothesis on $\CW$ we have $\h = \mathcal{H}_0(\varphi)$ for some
$\varphi \nearrow \rho$. If $\varphi \simeq \varphi^* \otimes \eps$,
we have $\dim \varphi = 2 \rk \h$. Otherwise, either $\varphi^*
\otimes \eps \neqnear \rho$ and $\dim \varphi = \rk \h + 1
\leq 2 \rk \h$ since $\rk \h \geq 1$, or $\varphi^* \otimes \eps
\nearrow \rho$ and $\dim \varphi + \dim \varphi^* \otimes \eps
= 2 \rk \h + 2$. It follows that $\dim soc(\bar{\rho}) \leq
2r + d$
where $d = \# \{ \varphi \nearrow \rho \ | \ \varphi^* \otimes \eps
\not\simeq \varphi \mbox{ and } \varphi^* \otimes \eps \nearrow \rho \}$.
Recall from lemma \ref{lemref23} that the 2-dimensional irreducible representations of $W$
are always quasi-reflection representations. Thus the assumption
$\rho \not\in \Ind \LRef$ implies $\dim \varphi \geq 3$ for
all $\varphi \nearrow \rho$. In particular $\dim soc(\bar{\rho}) \geq
3d$ hence, by using $\dim soc(\bar{\rho}) \leq 2r + d$, we get $\frac{2}{3} \dim soc(\bar{\rho}) \leq 2r$.
In case $\forall \rho' \in \Irr(\CW) 
\ (\bar{\rho}|\rho') \leq 1$ it follows that $\dim \rho \leq 3r 
< (r+1)^2$, and the conclusion follows from lemma
\ref{lemredsimpstrong} using
assumption (I). In the other case, if in addition $\rk \mathcal{H}_0(\bar{\rho})
\geq 4$, we also get $\dim \rho \leq 6r 
< (r+1)^2$ and the conclusion follows from lemma \ref{lemredsimpstrong}
using assumption (II) and lemma \ref{lempasrk1}.

We can thus assume that $r \leq 3$, that $\# \mathcal{R}
/W = 1$ and that $\forall \rho \in \Irr(W) \ \dim \rho \leq 3 \Rightarrow
\dim \rho = 1$.
We notice
that, if $\h = \mathcal{H}_0(\varphi)$ is a simple Lie ideal of
$\mathcal{H}_0(\bar{\rho})$ afforded by $\varphi \nearrow \rho$,
we have $\rk \mathcal{H}_0(\varphi) \geq 2$ because otherwise
by the hypothesis on $\CW$ we would have $\dim \varphi = 2$
contradicting $\rho \not\in \Ind\LRef$. Since $r \leq 3$ it follows that
$\mathcal{H}_0(\bar{\rho}) = \h$ is simple and that $r \geq 2$. 
In particular $soc(\bar{\rho}) = \varphi$ or $soc(\bar{\rho}) = \varphi
+ \varphi^*\otimes \eps$ for some $\varphi \in \Irr(\CW)$ with
$\dim \varphi \geq 3$. It follows that $\dim \rho \leq 4 \dim \varphi$.
On the other hand, since $\varphi \not\in \LRef$ we have
$\dim \varphi \in \{ r +1, 2r \}$ hence $\dim \varphi \leq 2r$
and $\dim \rho \leq 8 r \leq 24$.

Now, if $\mathcal{H}(\rho)$ is not
simple, since it cannot have simple Lie ideals of rank 1
by lemma \ref{lempasrk1}, its rank is
at least 4. If we decompose $\mathcal{H}(\rho) \simeq \h_1 \times
\h_2 \times \dots \times \h_m$ with $m \geq 2$,
$\h_i$ simple and $\rk \h_i \geq 2$, then $\rho_{\mathcal{H}'} \simeq V_1 \otimes \dots
\otimes V_m$ with $V_i$ an irreducible faithful representation
of $\h_i$. Since $\rk \h_i \geq 2$ we have $\dim V_i \geq 3$ and
$\dim \rho \geq 3^m$. Since $m \geq 2$ and $\dim \rho \leq 24$
it follows that $m = 2$. Moreover, $\rho \simeq \rho^* \otimes
\eps \Leftrightarrow V_1 \simeq V_1^*$ and $V_2 \simeq V_2^*$.
In that case, we actually have $\dim V_i \geq 4$,
hence $\dim \rho \geq 16$. We will also use the following
elementary facts : since $\rk \h_i \geq 2$, if $\dim V_i = 3$
then $\h_i$ has Cartan type $A_2$ ; if $\dim V_i = 4$, then either
$\h_i$ has Cartan type $B_2=C_2$ if
$V_i \simeq V_i^*$ or it has Cartan type $A_3$.

We denote by $\psi = \psi_1 \times \psi_2 : \mathcal{H}_0(\bar{\rho})
\to \mathcal{H}(\rho) \simeq \h_1 \times \h_2$ the natural
inclusion. Notice that, since $\mathcal{H}_0(\bar{\rho})$
is simple, $\dim \mathcal{H}_0(\bar{\rho}) > \dim \h_i$
or $\rk \mathcal{H}_0(\bar{\rho}) > \rk \h_i$ imply $\psi_i = 0$.
On the other hand, $\bar{\rho}$ admits irreducible
components with multiplicity at most 2 by assumption.
Since $\CW$ has a single class of reflections, the same
is true for $\bar{\rho}_{\mathcal{H}'_0} \simeq \Res_{\mathcal{H}'_0} V_1
\otimes V_2$. Let $i \in \{ 1, 2 \}$. Since $\dim V_i \geq 3$
it follows that $\mathcal{H}_0(\bar{\rho})$ is not included in any
simple ideal of $\mathcal{H}(\rho)$, otherwise $\bar{\rho}$
would admit irreducible components of multiplicity at least 3, hence
$\psi_i$ is non-zero, hence injective. It follows that
$\dim \mathcal{H}_0(\bar{\rho}) \leq \dim \h_i$ and 
$\rk \mathcal{H}_0(\bar{\rho}) \leq \rk \h_i$. In particular,
if $r = 3$ we have $\rk \h_i \geq 3$ hence $\dim V_i \geq 4$
and $\dim \rho \geq 16$.

Since the cases $\bar{\rho} = \varphi$
and $\bar{\rho} = \varphi+\varphi^*\otimes \eps$ with
$\varphi \not\simeq \varphi^*\otimes \eps$ have already been tackled,
only the following possibilities still have to be investigated :
$\bar{\rho} = 2 \varphi$ or  $\bar{\rho} \in \{
2 \varphi + 2\varphi^* \otimes \eps,2 \varphi + \varphi^* \otimes \eps\}$
with
$\varphi \not\simeq \varphi^* \otimes \eps$.

We will use the well-known fact that,
if $U_1,U_2$ are $n$-dimensional irreducible representations
of $\sl_n$, namely the standard representation or its dual,
then $U_1 \otimes U_2$ has 2 distinct irreducible components.
It follows that,
if $\mathcal{H}_0(\bar{\rho})$ has type $A_r$, then the case
$\dim V_1 =\dim V_2 = r+1$ is excluded. Moreover,
in that case the
possibility $\dim V_1 = r+1$, $\dim V_2 = r+2$ is also excluded
because a $(n+1)$-dimensional faithful representation
of $\sl_n$ has to be a direct sum of the trivial
representation and of a $n$-dimensional irreducible
representation, thus leading to 3 irreducible components. 

We can now proceed to the separate study of the special cases.

\begin{itemize}
\item $\bar{\rho} = 2 \varphi$. If $\varphi \not\simeq \varphi^*\otimes \eps$
we have $\dim \rho = 2r+2 \leq 8 < 9 \leq (r+1)^2$ hence
this case is handled by lemma \ref{lemredsimpstrong} (II)
for $\mathcal{H}_0(\bar{\rho})$. We thus can assume
$\varphi \simeq \varphi^*\otimes \eps$, hence $\dim \rho = 4r \leq 12$.
If $r = 2$ we still have $4r = 8 < 9 = (r+1)^2$ and we conclude
as before. There only remains the case $r = 3$ and $\dim \rho = 12 < 16$, a contradiction.

\item $\bar{\rho} = 2 \varphi + \varphi^* \otimes \eps$ with
$\varphi \not\simeq \varphi^* \otimes \eps$. Since
$\varphi \not\simeq \varphi^* \otimes \eps$ we have $\dim
\rho = 3(r+1)$. If $r =3$, $\dim \rho < (r+1)^2$ and we conclude
by lemma \ref{lemredsimpstrong} (II). We thus assume $r = 2$. Then
$\dim \rho = 9$ and $\dim V_1 = \dim V_2 = 3=r+1$. It follows that
$\mathcal{H}(\rho)$ has Cartan type $A_2 \times A_2$. Its Lie subalgebra
$\mathcal{H}_0(\bar{\rho})$ has Cartan type $A_2$, thus leading to
a contradiction by the remark above.

\item $\bar{\rho} = 2 \varphi + 2\varphi^* \otimes \eps$. If $r = 2$,
we have $\dim \rho = 12$, hence we can assume $\dim V_1 = 3$, hence $\h_1 \simeq
\sl_3$, and $\dim V_2 = 4$ ; moreover, $\mathcal{H}_0(\bar{\rho})$
has Cartan type $A_2$. Since
$\psi_1,\psi_2$ are injective, $\psi_1$
is an isomorphism. Let $E, E^*$ denote the two smallest faithful
representations of $\sl_3$. Since they
have dimension 3, as a
representation of $\mathcal{H}_0(\bar{\rho})\simeq \sl_3$
the representation $V_2$ can be decomposed as the trivial representation
plus either $E$ or $E^*$.
Similarly, as a representation of
$\mathcal{H}_0(\bar{\rho})$, $V_1 \simeq E$ or $V_1 \simeq E^*$.
It
follows that $\bar{\rho}$ should have 3 irreducible components,
a contradiction.

We thus can assume $r = 3$,
hence $\dim \rho = 16$, $\dim V_1 = \dim V_2 = 4$. 
Moreover $\mathcal{H}_0(\bar{\rho})$ has Cartan type $A_3$, hence
$\rk \h_i \geq 3$ and $\h_i \simeq \sl_4$. It follows that
the restriction of $V_1 \otimes V_2$ to $\mathcal{H}_0(\bar{\rho})$
is one of the inner tensor products $F \otimes F$, $F \otimes F^*$,
$F^* \otimes F^*$ with $F$ the standard representation of $\sl_4$.
The fact that these tensor products have only two irreducible
components yields a contradiction.
\end{itemize}

\end{proof}

\begin{lemma} \label{lemindmult1} Assume that $\CW$ admits a single
class of reflections. Let $\rho \in \Irr(W) \setminus \Ind\LRef$. If
$\forall \rho' \in \Irr(\CW)$ we have $(\bar{\rho}|\rho') \leq 1$,
then
$\mathcal{H}(\rho) \simeq \mathcal{L}(\rho)$.
\end{lemma}
\begin{proof} We decompose $\bar{\rho}$ in irreducible components
$$
\bar{\rho} = \sum_{i=1}^r \rho_i + \sum_{i=1}^s \varphi_i + \varphi_i^* \otimes
\eps + \sum_{i=1}^t \psi_i
$$
with $\rho_i^* \otimes \eps \neqnear \rho$, $\varphi_i \not\simeq
\varphi^* \otimes \eps$ and $\psi_i \simeq \psi_i^* \otimes \eps$.
We know that $\mathcal{H}(\rho)$ is simple by lemma \ref{lemnormsimp}.
Moreover,
$$
\rk \mathcal{H}_0(\bar{\rho}) = 
\sum_{i=1}^r \left(\dim \rho_i - 1 \right) +
\sum_{i=1}^s \left( \dim \varphi_i - 1 \right)
+ \sum_{i=1}^t \frac{\dim \psi_i}{2}
= \frac{\dim \rho}{2}
+ \left( \sum_{i=1}^r \frac{\dim \rho_i}{2} \right) -r-s.
$$
It follows that $\rk \mathcal{H}_0(\bar{\rho}) > \frac{\dim \rho}{4}$
as soon as
$$\frac{\dim \rho}{4} > r + s - \frac{1}{2} \sum_{i=1}^r \dim \rho_i
$$
Since $\rho \not\in \Ind\LRef$ we have $\dim \rho_i \geq 3$ thus
$\frac{\dim \rho_i}{2} \geq \frac{3}{2} > 1$ hence
$r - \frac{1}{2} \sum_{i=1}^r \dim \rho_i < 0$ if $r \geq 1$.
Similarly $\dim \varphi_i \geq 3$ hence $\dim \rho \geq 2 \sum_{i=1}^s
\dim \varphi_i \geq 6s$ and $s \leq \frac{\dim \rho}{6} < \frac{\dim \rho}{4}$.
It follows that $\rk \mathcal{H}_0(\bar{\rho}) > \frac{\dim \rho}{4}$.
Then $\mathcal{H}(\rho) = \mathcal{L}(\rho)$ by lemmas \ref{critsl}
and \ref{caractBCD}, at least if $\dim \rho \geq 22$ or $\rk \mathcal{H}(
\rho) \geq 7$.
To complete the proof, we can thus assume that
$\dim \rho \leq 21$ and $\rk \mathcal{H}(
\rho) \leq 6$. We first assume that $\rho_{\mathcal{H}'}$ is \emph{not}
selfdual. In particular $\rho \not\simeq \rho^* \otimes \eps$.

Our first task is to exclude the case $r = 0$, that is $\bar{\rho}
\simeq \bar{\rho}^* \otimes \eps$. Assume by contradiction that $r = 0$.
A first consequence is that, since each $\psi_i$ has even dimension,
then $\dim \rho$ is also even. Since $\rho_{\mathcal{H}'}$ is not selfdual, the
remaining exceptions are
\begin{enumerate}
\item $\rk \mathcal{H}(\rho) = 5$, $\mathcal{H}(\rho) \simeq \so_{10}$,
$\dim \rho = 16$.
\item $\rk \mathcal{H}(\rho) = 4$, $\mathcal{H}(\rho) \simeq \sl_{5}$,
$\dim \rho = 10$.
\end{enumerate}
Moreover, we have $2 (s+t) \leq \rk \mathcal{H}_0(\bar{\rho}) \leq 5$
hence $s+t \leq 2$. Note that $s+t$ is the number of simple Lie
ideals of $\mathcal{H}_0(\bar{\rho})$. We first consider the case
$s+t = 2$, that is $\mathcal{H}_0(\bar{\rho}) \simeq \h_1 \times \h_2$.
We have  $\rk \h_i \in \{ 2, 3 \}$ and $\rk \h_1 + \rk \h_2 \in \{ 4, 5 \}$.
Exhausting all possibilities we find that this is not possible
if $\dim \rho = 16$, and that the only possibilities for $\dim \rho = 10$
are of Cartan type $C_3 \times C_2$, $D_3 \times C_2$ and $A_2 \times C_2$.
Since $\dim \rho = 10 \Leftrightarrow \rk \mathcal{H}(\rho) = 4$
and $\rk \mathcal{H}_0(\bar{\rho}) \leq \rk \mathcal{H}(\rho) = 4$,
we have $\mathcal{H}_0(\bar{\rho}) \simeq \sl_3 \times \sp_4$. But
$\sl_3 \times \sp_4$ does not embed into $\sl_5$ (for instance because
the smallest faithful representation of $\sl_3 \times \sp_4$ has
dimension $3+4=7$), a contradiction.

Now we know that $r \geq 1$. By the calculation above we have
$\rk \mathcal{H}(\rho) > \frac{\dim \rho}{2}$ as soon as
$\sum_i \dim \rho_i > 2(r+s)$. Since $r \geq 1$ and $\dim \rho_i \geq 3$,
if $s=0$ then $\mathcal{H}(\rho) = \sl(\rho) = \mathcal{L}(\rho)$ by lemma
\ref{critsl}. We thus have $r \geq 1, s \geq 1$. Since every representations
of $\CW$ involved here have dimension at least 3 we have $\dim \rho \geq 9$.
On the other hand, since $\rho \not\in \Ind \LRef$ we know that
$\mathcal{H}_0(\bar{\rho})$ does not contain Lie ideals of rank 1 hence
$\rk \mathcal{H}(\rho) \geq \rk \mathcal{H}_0(\bar{\rho}) \geq 2(r+s+t)$.
On the other hand $\rk \mathcal{H}(\rho) \leq 6$ hence $r+s+t \leq 3$.
Since $r, s \geq 1$ we also have $\rk \mathcal{H}(\rho) \geq 4$. Moreover 
the case $\rk \mathcal{H}(\rho) = 4$ can be ruled out,
for it implies $r = s = 1$ and $t = 0$,
thus $r+s = 2 < \frac{\dim \rho}{2}$ except if $\dim \rho_1 = 4$.
But $\dim \rho_1 = 4$ implies $\rk \mathcal{H}(\rho) \geq \rk \mathcal{H}(\rho_1)
+ \rk \mathcal{H}(\rho_2) \geq 3+2=5$, a contradiction.

We are left with two cases : either $\rk \mathcal{H}(\rho) = 5$
or $\rk \mathcal{H}(\rho) = 6$. We first deal with the former one.
If $r \geq2$ or $s \geq 2$ or $t \geq 1$ we would
have $\rk \mathcal{H}(\rho) \geq 6$, so we know that $r=1,s=1,t=0$,
with $\rk \mathcal{H}(\rho_1) = \mathcal{H}(\rho_2) = 2$, hence
$\dim \rho = 3+3+3 = 9$. But there are no 9-dimensional representations
in the list of lemma \ref{caractBCD}.

It follows that $\rk \mathcal{H}(\rho) = 6$ hence $\mathcal{H}(\rho) \simeq
\sl_7$ and $\dim \rho = 21$. If $t \neq 0$, from $6 \geq 2(r+s+t)$ we
get $r=s=t=1$ and $\rk \mathcal{H}(\rho_1) = \mathcal{H}(\varphi_1)
= \mathcal{H}(\psi_1) = 2$ whence $\dim \rho = 3+2 \times 3+4 = 13 \neq 21$.
Hence $t = 0$.If $r+s = 3$ we also have $\rk \mathcal{H}(\rho_i)
= \mathcal{H}(\varphi_j) = 2$ hence $\dim \rho_i = \dim \varphi_j = 3$
and $\dim \rho \leq 3(r+2s) \leq 9+3s \leq 15 < 21$ since $r \geq 1$
and $s \leq 2$. It follows that $r=s=1$ and $t=0$. Letting
$\alpha = \rk \mathcal{H}(\rho_1)$ and $\beta = \rk
\mathcal{H}(\varphi_1)$ we have $4 \leq \alpha + \beta \leq 6$
and $\dim \rho = \alpha+1 + 2 (\beta+1) \leq 2(\alpha+\beta)+3 \leq 15 < 21$,
again a contradiction.

We now assume that $\rho_{\mathcal{H}'}$ is selfdual. This implies
that $\bar{\rho}_{\mathcal{H}'}$ is selfdual, hence $\bar{\rho}\simeq \bar{\rho}^* \otimes \eps$,
since $\CW$ has a single class of reflections, and $r=0$,
$\mathcal{H}(\rho)$ is included either in $\sp_N$ or
$\so_N$ with $N = \dim \rho$.
The only selfdual representations in the list of lemma
\ref{caractBCD} are for $\dim \rho \leq 8$ and $\rk
\mathcal{H}(\rho) \leq 4$. Moreover, if $\dim \rho = 8$ or
$\rk \mathcal{H}(\rho) = 4$, this means that $\mathcal{H}(\rho)$
has type $D_4$ and in that case again $\mathcal{H}(\rho) \simeq \mathcal{L}
(\rho)$ since it is included in $\so_8$ or $\sp_8$ and
there are no inclusion from $\so_8$ to $\sp_8$. It follows
that $\dim \rho \leq 6$ and $\rk \mathcal{H}(\rho) \leq 3$. Since
$\mathcal{H}_0(\bar{\rho})$ has no ideal of rank 1 it follows that
$s+t = 1$.

If $s = 1$, that is $\bar{\rho} = \varphi_1 + \varphi_1^* \otimes \eps$, from
$\dim \rho \leq 6$ we get $\dim \varphi_1 \leq 3$ hence $\varphi_1 \in
\LRef(W_0)$ by lemma \ref{lemref23}, and $\rho \in \Ind\LRef$, a contradiction.
We thus have $s=0, t = 1$, $\bar{\rho} = \psi_1$.

If $\rk \mathcal{H}_0(\psi_1)
= 2$ we have $\dim \psi_1 = \dim \rho = 4$, and the only exception is
of type $B_2$, which is a fake exception since $B_2 \simeq C_2$.
If $\rk \mathcal{H}_0(\psi_1) = 3$ we have $\dim \psi_1 = \dim \rho = 6$.
But the only 6-dimensional irreducible representation of $\sl_4$
is orthogonal, hence $\mathcal{H}(\rho) \simeq \so_6 \simeq \sl_4$, hence
$\mathcal{L}(\rho) \simeq \mathcal{H}(\rho)$ in that case too.
\end{proof}

If $\bar{\rho} = a_1 \rho_1 + \dots + a_r \rho_r$ with $\rho_i \neq \rho_j$,
we denote by
$\delta(\rho) = a_1 + \dots + a_r$ the number of irreducible components
(counting multiplicities) of $\overline{\rho}$.

\begin{lemma} \label{lemindmult2one} Asssume that $W$ and $\CW$ admits a single class of reflections
and that $\forall \rho \in \Irr(W) \ \dim \rho \leq 3 \Rightarrow \dim \rho = 1$.
Let $\rho \in \Irr(W) \setminus \Ind \LRef$, and $\zeta \nearrow \rho$.
Assume $\delta(\rho) \geq 3$, $(\zeta | \bar{\rho}) \leq 2$,
and $\rho' \nearrow \rho \Rightarrow \dim \rho' \geq 5$. If $\forall \rho'
\in \Irr(\CW) \ \rho' \not\simeq \zeta \Rightarrow (\rho' | \bar{\rho}) \leq 1$
then $\mathcal{H}(\rho) = \mathcal{L}(\rho)$.
\end{lemma}
\begin{proof}
If $(\zeta | \bar{\rho}) = 1$ then $\mathcal{H}(\rho) = \mathcal{L}(\rho)$
by lemma \ref{lemindmult1}, so we can assume $(\zeta| \bar{\rho}) = 2$.
By lemma \ref{lemnormsimp} we know that $\mathcal{H}(\rho)$ is a simple
Lie algebra. We decompose
$$
\bar{\rho} = \zeta + \sum_{i=1}^r \rho_i + \sum_{i=1}^s \left( 
\varphi_i + \varphi_i^* \otimes \eps \right) +
\sum_{i=1}^t \psi_i
$$
with $\zeta$ isomorphic to one of the $\rho_i, \varphi_i, \psi_i$, with
$\rho_i \not\simeq \rho_i^* \otimes \eps$, $\varphi_i
\not\simeq \varphi_i^* \otimes \eps$, $\rho_i^* \otimes \eps \neqnear
\rho$ and $\psi_i \simeq \psi_i^* \otimes \eps$. By assumption,
we have $r + s+t \geq 2$.

Like in the proof of lemma \ref{lemindmult1}, we have
$\rk \mathcal{H}_0(\bar{\rho}) > \frac{\dim \rho}{4}$ iff
$\frac{\dim \rho}{4} > A$ where
$$ A = r + s + \frac{\dim \zeta}{2} - \frac{1}{2} \sum_{i=1}^r \dim \rho_i.
$$
From this proof, we also recall that $r - \frac{1}{2} \sum_{i=1}^r
\dim \rho_i < 0$ if $r \geq 1$, since $\dim \rho_i \geq 5 \geq 3$.
Likewise, $\dim \rho \geq 10s$ since $\dim \varphi_i \geq 5$ hence
$s \leq \frac{\dim \rho}{10}$. Finally, $\dim \rho \geq 15$
since $\delta(\rho) \geq 3$.

There are three possibilities
\begin{enumerate}
\item $\zeta \simeq \rho_{i_0}$ for some $i_0 \in [1,r]$. Since
$\dim \rho_i \geq 5$ we have
$$
r + \frac{\dim \zeta}{2} - \frac{1}{2} \sum_{i=1}^r \dim \rho_i =
r - \frac{1}{2} \sum_{i \neq i_0} \dim \rho_i \leqslant
r - \frac{5}{2}(r-1) = \frac{-3}{2} r + \frac{5}{2} = \frac{5-3r}{2}
< 0
$$
as soon as $r\geq 2$. In this case, $A < s \leq \frac{\dim \rho}{10} <
\frac{\dim \rho}{4}$. If $r= 1$, then
$A = s+1 \leq 1 + \frac{\dim \rho}{10} < \frac{\dim \rho}{4}$ since
$\dim \rho \geq 15$.
\item $\zeta \simeq \psi_{i_0}$ for some $i_0 \in [1,t]$.
In that case, $\dim \rho \geq 2 \dim \zeta + 10s$ hence
$\frac{\dim \rho}{4} \geq \frac{\dim \zeta}{2} + s + \frac{3s}{2}$
and $A < \frac{\dim \rho}{4}$ as soon as $s > 0$ or $r > 0$. If
$s = r= 0$ then $t \geq 2$, $\dim \rho \geq 2 \dim \zeta + 5 > 2
\dim \zeta$ and $A = \frac{\dim \zeta}{2} < \frac{\dim \rho}{4}$.
\item $\zeta \simeq \varphi_{i_0}$ for some $i_0 \in [1,s]$.
Then $\zeta^* \otimes \eps \nearrow \rho$ and
$\dim \rho \geq 3 \dim \zeta + 10(s-1)$. It follows that
$$
\frac{\dim \rho}{6} \geqslant \frac{\dim \zeta}{2} + \frac{5}{3}(s-1)
= \frac{\dim \zeta}{2} + s + \frac{2}{3} s - \frac{5}{3} =
s + \frac{\dim \zeta}{2} + \frac{2s-5}{3}
$$
hence $A \leq \frac{\dim \rho}{6} + \frac{5-2s}{3} < \frac{\dim \rho}{4}$
since $\dim \rho \geq 15$ and $s \geq 1$.
\end{enumerate}
To conclude we apply lemmas \ref{critsl} and \ref{caractBCD}. The
exceptions in lemma \ref{caractBCD} that we have to rule out
are for $\dim \rho \in \{ 15,16,21 \}$, with $\rho \simeq \rho^* \otimes \eps$,
and $\rk \mathcal{H}(\rho) \in \{ 5,6 \}$. Since $r+s+t \geq 2$,
we have $\rk \mathcal{H}(\rho) \geq 4r + 4s + 3t > 3(r+s+t) \geq 6$
as soon as $(r,s) \neq (0,0)$. If $(r,s) = (0,0)$ then $\rk \mathcal{H}_0(\bar{\rho})
\geq 3t \geq 6$, with equality holding iff $t=2$, $\rk \mathcal{H}(\psi_i)
= 3$, $\dim \psi_i = 6$ and $\dim \rho = 18$, contradicting $\dim \rho \in
\{ 15,16,21 \}$.
\end{proof}

\subsection{Reduction to one conjugacy class}

The goal of this section is to show that, using the lemmas
above, we can assume that $W$ admits a single conjugacy class
of reflections (except when $W$ is a Coxeter group of type $F_4$ or $H_4$).

We use the Shephard-Todd classification to prove the following.

\begin{lemma} If $\Phi$ is an isomorphism for every irreducible $W \neq H_4$ with
$\# \mathcal{R}/W = 1$, then $\Phi$ is an isomorphism for every irreducible $W \not\in \{ F_4, H_4 \}$.
\end{lemma}

\begin{proof}
Let $W \not\in \{ F_4, H_4 \}$ with $\# \mathcal{R}/W > 1$. The only case in the
exceptional series is for $W$ of type $G_{13}$, for which we check
by a direct computation of the dimensions that $\Phi$ is surjective.
We thus assume that $W$ belongs to the infinite series $G(2e,e,r)$
or $G(e,e,r)$. When $r = 2$, lemma \ref{lemdim2rang2} states
that $\Irr(W) = \LRef(W)$, and we know that $\Phi$ is an isomorphism
in that case. Assuming $r \geq 3$, it follows from lemma \ref{lemgeersingle}
that $W$ has type $G(2e,e,r)$ for some $e \geq 1$. Let $\CW$ be its
natural reflection subgroup of type $G(2e,2e,r)$, for which the
theorem is assumed to hold. By corollary \ref{corimpltheo} we have to prove
that, for any $\rho \in \Irr(W) \setminus \LRef$, then $\mathcal{H}(\rho)
\simeq \mathcal{L}(\rho)$. Now $\CW$ has index 2 in $W$, hence
by Clifford theory $\bar{\rho}$ has at most two irreducible components,
and this decomposition is multiplicity-free. If $\rho \not\in \Ind\LRef$
the conclusion follows from lemma \ref{lemindmult1}. We thus assume
$\rho \in \Ind\LRef$ and $\rho \not\in \LRef$.

Let $\bla$ be an irreducible representation of $G(2e,1,r)$ whose
restrictions contains $\rho$. In particular, its restriction to
$G(2e,2e,r)$ belongs to $\LRef(W_0)$, and is thus listed by table \ref{tabqrefrep} and proposition \ref{listeLRef}.
We check that, for $r \geq 5$, this implies $\rho \in \LRef$, and that
the only cases to consider are when $\bla$ has one of the two special
shapes at the bottom of table \ref{tabqrefrep} (the last two lines).

The first shape to deal with is for $r = 3$, in which case
we have $A(\bla) = 3$. Since $A(\bla)/B(\bla) \in \{ 1, 2 \}$
it follows that $B(\bla) = 3$ hence $\rho$
has dimension 2 and belongs to $\QRef \subset \LRef$.
The second one is for $r = 4$, in which case we have
$A(\bla) = 2$ and $B(\bla) \in \{1,2 \}$. First
assume $B(\bla) = 1$, meaning that the restriction
$\rho$ of $\bla$ to $W$ is irreducible. Then $\bar{\rho}$ is
the sum of two quasireflection representations $\rho_1,\rho_2$
of dimension 3 of $W_0$. Since $W_0$ has a single reflection
class, these correspond to distinct simple
ideals provided that $\rho_2 \not\simeq \rho_1^* \otimes \eps$. In
order to check this it is enough to check $\bar{\rho} \not\simeq
\bar{\rho}^* \otimes \eps$, which is true since $\bar{\rho}^* \otimes \eps$
is the restriction of an irreducible representation $\bmu$ of
$G(2e,1,r)$ with $\{ \mu^i \} = \{ (\la_i)' \}$, and either
$\{ \la^i  \} = \{ \emptyset, [2] \}$ or $\{ \la^i \} = \{ \emptyset, [1,1] \}$.
Thus $\rk \mathcal{H}(\rho) \geq 4 > 6/2$ and
$\mathcal{H}(\rho)\simeq \sl(V_{\rho}) \simeq \mathcal{L}(\rho)$
by lemma \ref{critsl}.
Now assume $B(\bla) = 2$. Then $\bar{\rho}$ is an irreducible
quasireflection representation of $\CW$, hence
$\mathcal{H}(\rho) \simeq \sl(\rho) \simeq \mathcal{H}_0(\bar{\rho})$,
which concludes the proof.
\end{proof}

\subsection{Induction process for special cases}

\subsubsection{Three irreducible components}

We assume that $W_1 \subset
W_0 \subset W$ is a chain of inclusions between irreducible
reflection groups which respect the reflections, and such that
$W_1,W_0$ and $W$ admit a single conjugacy class of
reflections. We assume that the theorem holds for $W_0$.

Let $\rho \in \Irr(W)$, and
assume that the restriction diagram of $\rho$ with respect to $W_1 \subset
W_0 \subset W$ has the following form
$$
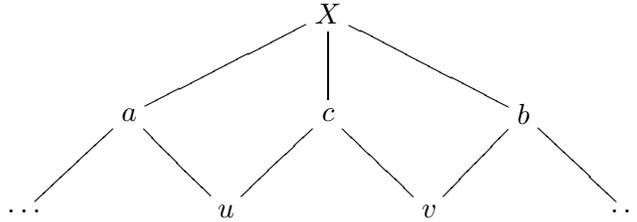
\xymatrix{
 & & & X & & \\
& a \ar@{-}[urr] \ar@{-}[dl] \ar@{-}[dr] & & c \ar@{-}[u] \ar@{-}[dl]
\ar@{-}[dr] & & b \ar@{-}[dl] \ar@{-}[dr] \ar@{-}[ull] & \\
\dots & & u & & v & & \dots \\
}
$$
meaning that the restriction of $\rho$ to $W_0$ has (exactly) three
(irreducible) distinct components $a,c,b$, that the restriction of $c$ to $W_1$
has two distinct components $u,v$, that $u$ is a component of
the restriction of $a$ to $W_1$, that $v$ is a component of
the restriction of $b$ to $W_1$, and that  the restrictions of $a$ and $b$
to $W_1$ are not irreducible.

\begin{lemma} \label{lemres3} In the configuration above,
\begin{enumerate}
\item if $\mathcal{H}_0(a) \simeq \sl(V_a)$ and $\mathcal{H}_0(b) \simeq \sl(V_b)$,
with $\mathcal{H}_0(a)$, $\mathcal{H}_0(b)$ and $\mathcal{H}_0(c)$
corresponding
to distinct ideals of $\mathcal{H}'_{0}$, then $\mathcal{H}(\rho) \simeq \sl(
V_{\rho})$
as soon as $\rk \mathcal{H}_0(c) \geq 2$.
\item if $\mathcal{H}_0(a) \simeq \sl(V_a)$, $(\rho_{\mathcal{H}'})^* 
\simeq \rho_{\mathcal{H}'}$, $b_{\mathcal{H}'_0} = (a_{\mathcal{H}'_0})^*$,
and $\mathcal{H}_0(a)$, $\mathcal{H}_0(c)$ are
distinct ideals of $\mathcal{H}'_{0}$, then $\mathcal{H}(\rho) \simeq
 \osp(V_{\rho})$
as soon as $\rk \mathcal{H}_0(a)+\rk \mathcal{H}_0(c)\geq 5$ and $\rk \mathcal{H}(c) \geq 1$.
\end{enumerate}

\end{lemma}
\begin{proof}
By lemma \ref{lemredsimpstrong}
we know that $\mathcal{H}(\rho)$ is a simple Lie algebra.
We first prove (1). Since $\mathcal{H}(\rho)$
contains a copy of $\mathcal{H}_0(a) \times \mathcal{H}_0(b)
\times \mathcal{H}_0(c)$ one has
$$
\begin{array}{lclcl}
\rk \mathcal{H}(\rho) & \geqslant & \rk(\mathcal{H}_0(a))+\rk(\mathcal{H}_0(b)) +
\rk(\mathcal{H}_0(c)) & \geqslant & \dim(a) + \dim(b)
+ \rk(\mathcal{H}_0(c)) - 2 \\
& \geqslant & \dim(a) + \dim(b) &= & \dim(\rho) - \dim(c) \\
\end{array}
$$ 
Since $\dim(\rho) = \dim(a) + \dim(b) + \dim(c) \geq \dim(c) + \dim(u)
+ \dim(v) +2= 2 \dim(c)$ it follows that $\rk \mathcal{H}(\rho) >
\frac{\dim(\rho)}{2}$
hence $\mathcal{H}(\rho) = \sl(V_{\rho})$ by lemma \ref{critsl}.

We then prove (2). We have
$$
\begin{array}{lclcl}
\rk \mathcal{H}(\rho) & \geqslant & \rk(\mathcal{H}_0(a))+\rk(\mathcal{H}_0(c)) 
 & \geqslant & \dim(a) 
+ \rk(\mathcal{H}_0(c)) -1 \\
& \geqslant & \dim(a)  & > & \frac{\dim(\rho)}{4} \\
\end{array}
$$ 
because
$$
\begin{array}{lclcl}
\dim (\rho) & = & \dim(a) + \dim(b) + \dim(c) & = & 2 \dim(a) + \dim(u) + \dim(v) \\
& \leqslant & 2 \dim(a) + \dim(a) + \dim(b) -2 & < & 4 \dim(a).
\end{array}
$$
This implies $\dim \rho < 4(r+1)$ where $r = \rk \rho(\mathcal{H}'_0)$,
hence $\dim \rho < (r+1)^2$ since $r \geq \rk \mathcal{H}_0(a)
+ \rk \mathcal{H}_0(c) \geq 5$. Then lemma \ref{lemredsimpstrong} (I)
claims that $\mathcal{H}(\rho)$ is simple and, since $\rk \mathcal{H}(\rho)
\geq r \geq 5$, lemma \ref{caractBCD} with $(\rho_{\mathcal{H}'})^*
\simeq \rho_{\mathcal{H}'}$ together imply $\mathcal{H}(\rho) \simeq
\osp(V_{\rho})$.

\end{proof}

For future reference, we state as a lemma the following
simple fact.

\begin{lemma} \label{lemres2}If $\dim(\rho) > 4$ and the restriction of $\rho$
to $W_0$ has exactly two irreducible components $a,b$ corresponding
to distinct ideals $\mathcal{H}_0(a) \simeq \sl(V_a)$,
$\mathcal{H}_0(b) \simeq \sl(V_b)$, then $\mathcal{H}(\rho) \simeq \sl(V_{\rho})$.
\end{lemma}
\begin{proof} We have $\rk \mathcal{H}(\rho) \geq \dim(a) + \dim(b) -2
= \dim(\rho) - 2$ and $\dim(\rho) - 2 > \frac{\dim(\rho)}{2}$ as soon
as $\dim(\rho) > 4$. The conclusion follows from lemma \ref{critsl}.
\end{proof}

\subsection{Groups $G(e,e,r)$ of small rank}

The case $r = 2$ is known by lemma \ref{lemdim2rang2}.
Also note that we have $\rho \simeq \rho^* \otimes \eps$
for all 2-dimensional irreducible representations of $G(e,e,2)$;
hence distinct such representations define distinct ideals
of $\mathcal{H}$.

We now assume $r = 3$. Let $\rho$ be an irreducible representation
of $G(e,e,3)$, and $\bla$ a representation of $G(e,1,3)$
such that $\rho$ embeds in the restriction of $\bla$. If $p(\bla)\leq 2$,
then $\bla$ is a quasi-reflection representation and so is
$\rho$, whence $\mathcal{H}(\rho) = \mathcal{L}(\rho)$.
We thus can assume $p(\bla) = 3$, hence $\dim \bla = 6$. Since $A(\bla)$ divides
$p(\bla)$ it follows that $A(\bla) \in \{1,3 \}$. If $A(\bla) = 3$,
then $\dim \rho = 2$ and $\rho$ is a quasi-reflection representation.
Hence we can assume $A(\bla) = 1$ and $\dim \rho = 6$.
We can assume $\la^0 = \la^i = \la^j = [1]$ with $\# \{ 0, i,j \} = 3$.
The restriction of $\bla$ to $G(e,1,2)$ is the sum of three irreducible components
$\bla_1$, $\bla_2$ and $\bla_3$ of dimension 2, with $p(\bla_i) = 2$.
In particular $A(\bla_i) \in \{ 1, 2 \}$. Let $\rho_i = < \bla^i>$.
Notice
that $\mathcal{H}(\rho)$ cannot have
Cartan type $C_2$, as $\sp_4 \simeq \so_5$ does not
admit an irreducible 6-dimensional representation. It cannot be
of type $G_2$ or $B_3$ for the same reason.
We separately consider the following possibilities.
\begin{itemize}
\item If $A(\bla_1)=A(\bla_2)=A(\bla_3)=1$ and
the restrictions to $G(e,e,3)$ of the $\bla_i$ are irreducible
representations. Then the restriction of $\rho$ to
$G(e,e,2)$ is multiplicity free, hence by lemma \ref{lemredsimpstrong}
(I) $\mathcal{H}(\rho)$ is simple. If $\rk \mathcal{H}(\rho) \geq 4$ then $\mathcal{H}(\rho)
= \sl(V_{\rho})$ by lemma \ref{critsl} and we are done. On the other hand, $\mathcal{H}(\rho) \supset \sl_2^3$
hence we can assume $\rk \mathcal{H}(\rho) = 3$. Since $(\sl_2)^3$
does not embed
in $\sl_4$ or in a Lie algebra of Cartan type $B_3$,
we get that $\mathcal{H}(\rho)$ must have Cartan type $C_3$. Since $\dim \rho = 6$
it follows that $\rho_{\mathcal{H}'}$ is selfdual hence
$\mathcal{H}(\rho)  \simeq \mathcal{L}(\rho)$.
\item If $A(\bla_1)=A(\bla_2)=1$, $A(\bla_3) = 2$ and $\bla_1 ,\bla_2$
have isomorphic restriction to $G(e,e,3)$. This situation
can occur only if
$e = 4g$ for some integer $g$, and we can assume
$ i = g, j = 2g$. But then $\rho$ factorises
through the morphism $G(4g,4g,3) \onto G(4,4,3)$, and it
is sufficient to check that $\mathcal{H}(\rho) \simeq \mathcal{L}(\rho)$
for the restriction to $G(4,4,3)$ of the representation $([1],[1],[1],\emptyset)$,
which we do by computer.
\end{itemize} 
In all other cases, the restriction of $\bmu$ admits multiplicities,
and also contains some 2-dimensional irreducible component,
whose corresponding ideal is isomorphic to $\sl_2$, with multiplicity 1.
Lemma \ref{lemliedim6} tackles these cases, hence the theorem holds
for $r = 3$ too.

\subsection{Induction process for $W = G(e,e,r)$, $\rho \in \Ind\LRef$}
We first consider the case where the restriction of $\rho$
to $G(e,e,r-1)$ contains one of the exceptional
quasireflection representations of table \ref{tabqrefrep}. The
first case is when $r-1 = 4$, and $e$ is even. In this case,
up to $\rho \rightsquigarrow \rho^* \otimes \eps$,
wa can assume $\la^0 = \la^g = [2]$, for $e = 2g$. Notice
that the corresponding 3-dimensional representations are not
selfdual, nor dual of each other, as representations of
$\mathcal{H}'_0$. The second case is when $r-1 = 3$,
and $e$ is divisible by 3. Then the corresponding
representations of $\mathcal{H}'_0$, being 2-dimensional,
are selfdual and correspond to distinct ideals
of $\mathcal{H}'_0$, since $\# \mathcal{R}_0 / W_0 = 1$.

\subsubsection{The special cases of $G(2g,2g,4)$ and $G(3g,3g,3)$}

\medskip

\paragraph{$G(2g,2g,4)$}

Let $\rho_0,\rho_1$ denote the two irreducible (3-dimensional) components
of the restriction to $G(2g,2g,4)$ of the representation $\bla$
of $G(2g,1,4)$ defined by $\la^0 = [2]$, $\la^g = [2]$.
Note that $\rho_1 \not\simeq \rho_0^* \otimes \eps$
hence $\mathcal{H}(\rho_0) \neq \mathcal{H}(\rho_1)$.
Assume that $\rho$ is an irreducible representation of $G(2g,2g,5)$
such that $\rho_0$ embeds in $\rho$. The exists an irreducible
representation $\bmu$ of $G(2g,1,5)$ with $\bla \nearrow \bmu$
such that $\rho$ embeds
in its restriction to $G(2g,2g,5)$.
We have $p(\bmu) \in \{ 2, 3 \}$. By lemma \ref{brisuregeer},
$A(\bmu)$ has to divide $5$ and $p(\bmu)$, hence
$A(\bmu)=1$ and $\dim \rho = \dim \bmu$.

If $p(\bmu) = 2$, first assume $\mu^0 = [3]$ and $\mu^g = [2]$,
and the restriction of $\rho$ to $G(2g,2g,4)$
decomposes as $\rho_0 \oplus \rho_1 \oplus \rho_2$
with $\rho_2$ a 4-dimensional reflection representation of $G(2g,2g,4)$.
If follows that $\rk \mathcal{H}(\rho) \geq 2+2+3 = 7$.
Since $\dim \rho = 10$ lemma \ref{critsl} implies
$\mathcal{H}(\rho) = \sl(V_{\rho})$. The other case is for $\mu^0 = [2,1]$,
$\mu^g = [2]$. Then the restriction of $\rho$ to $G(2g,2g,4)$ decomposes
as $\rho_0 \oplus \rho_1 \oplus \rho_2 \oplus \rho_3$, with
$\rho_2^* \otimes \eps \simeq \rho_2$, $\rho_3^* \otimes \eps
\simeq \rho_3$ with $\rho_2,\rho_3 \not\in \LRef$, $\dim \rho_2 = 6$,
$\dim \rho_3 = 8$. It follows that $\rk \mathcal{H}(\rho) \geq 11$,
hence $\mathcal{H}(\rho) \simeq \sl(V_{\rho})$ by
lemma \ref{critsl}, since $\dim \rho = 20$ in this case.

If $p(\bmu) = 3$ then we can assume $\mu^0 = \mu^g = [2]$
and $\mu^i = [1]$ for some $i \not\equiv 0 \mod g$. It follows
that the restriction of $\rho$ to $G(2g,2g,4)$
decomposes as $\rho_0 \oplus \rho_1 \oplus \rho_2\oplus \rho_3$
with $\rho_2,\rho_3$ distinct 12-dimensional representations of $G(2g,2g,4)$
such that $\rho_3 \not\simeq \rho_2^* \otimes \eps$.
By the induction hypothesis it follows that $\rk \mathcal{H}(\rho) \geq
2+2+11+11 = 26$. Since $\dim \rho = 30$ lemma \ref{critsl} implies
$\mathcal{H}(\rho) \simeq \sl(\rho)$.

\medskip

\paragraph{$G(3g,3g,3)$}

Let $\rho_0,\rho_1,\rho_2$ denote the three irreducible components of
the restriction to $G(3g,3g,3)$ of the representation $\bla$
of $G(3g,1,3)$ defined by $\la_0 = \la_g = \la_{2g} = [1]$.
Assume that $\rho$ is an irreducible representation of $G(3g,3g,4)$
such that $\rho_0$ embeds in $\rho$. Let $\bmu$ denote
a representation of $G(3g,1,4)$ such that $\rho$ embeds
in its restriction. By lemma \ref{brisuregeer} we have $A(\bmu)=1$
We have $\bla \nearrow \bmu$, and $p(\bmu) \in \{ 3, 4 \}$.

If $p(\bmu) = 3$, then we may assume, up to tensorization by
the sign character, that $\mu^0 = [2]$ and $\mu^g = \mu^{2g} = [1]$.
Then the restriction of $\rho$ equals $\rho_0\oplus\rho_1\oplus\rho_2$
plus two irreducible components, restriction of the representations
$\bmu_1, \bmu_2$ of $G(3g,1,4)$ defined by
$\mu_1^0 = \mu_2^0 = [2]$, $\mu_1^g = \mu_2^{2g} = [1]$, $\mu_1^{2g} = 
\mu_2^g = \emptyset$. It is easily checked that these five components
correspond to distinct ideals, and that the sum of their rank is
$7 > \frac{12}{2} = \frac{\dim \rho}{2}$ hence $\mathcal{H}(\rho)
\simeq \sl(\rho)$ by lemma \ref{critsl}.

If $p(\bmu) = 4$, then we may assume that $\mu_0 = \mu_g = \mu_{2g} = \mu_i
= [1]$ for some $0 < i < g$. The restriction of $\rho$ equals $\rho_0\oplus\rho_1\oplus\rho_2$
plus three irreducible 6-dimensional components $\tilde{\rho}_1,\tilde{\rho}_2,\tilde{\rho}_3$, restriction of the representations
$\bmu_1, \bmu_2,\bmu_3$ of $G(3g,1,4)$ defined by
$\mu_1^i=\mu_2^i=\mu_3^i= [1]$, $\mu_1^0 = \mu_1^g = \mu_2^0 = \mu_2^{2g} = 
\mu_3^g = \mu_3^{2g} = [1]$.
As representations of $\mathcal{H}'_0$, $\tilde{\rho}_1$ may be
selfdual, $\tilde{\rho}_2$ may be the dual of $\tilde{\rho}_3$, but
in any case the ideal corresponding to $\tilde{\rho}_2$ is isomorphic
to $\sl_6$ and different from the ones afforded by $\tilde{rho}_1$
or $\rho_0,\rho_1,\rho_2$. It follows that $\rk \mathcal{H}_0(\bar{\rho})
\geq 1+1+1+3+5 = 11$. Since $\dim \rho = 24$, lemma \ref{lemredsimpstrong} (I)
implies that $\mathcal{H}(\rho)$ is simple. Then parts (1) and (2) of lemma
\ref{lemdim5} (1) imply that $\mathcal{H}(\rho) \simeq \mathcal{L}(\rho)$.

\subsubsection{General case}

By \S 4 we know that the representations that we are
going to deal with correspond to very special
multipartitions $\bla$, so many information on $\bla$
will clear from the context. For brevity, we use the
following shortcuts. If $\la$ is some partition, we
denote by $(\la)$ a multipartition $\bla$ with $p(\bla) = 1$ such
that $\la^i = \la$ for some $i$. If $\la, \mu$ are two
partitions, we also denote $(\la;\mu)$ a multipartition
$\bla$ with $p(\bla) = 2$ and $\la^i = \la$, $\la^j = \mu$ for
some $i \neq j$. In rare cases, we similarly use
the notation $(\la;\mu;\nu)$ for a 3-parts multipartition.
For $\bla = (\la; \mu)$ and $A(\bla) = 2$, we denote
$(\la;\la)^{\pm}$ the irreducible components of the restriction
to $G(e,e,r)$. When the restriction is irreducible and the context
clear, we identify $\bla$ with its restriction to $G(e,e,r)$.
Notice that, using these shortcuts, we may write $(\la;\mu)^* \otimes \eps
= (\mu';\la')$. In particular, we can assume $\la \geq \la'$
(for the lexicographic order) when needed.

If $r \geq 4$, for the representations $c$ of $G(e,e,r-1)$
in $\LRef \setminus \mathrm{QRef}$ listed in proposition \ref{listeLRef},
we can apply lemma \ref{lemres3} with respect to the chain
of inclusions $G(e,e,r-2) \subset G(e,e,r-1) \subset G(e,e,r)$
to all representations $\rho$ of $G(e,e,r)$
in $\Ind\LRef$ such that $c \nearrow \rho$. The corresponding
representations
$a,b,u,v$ are listed in table \ref{tabindLRef}. 
To save space in the table, we let $r' = r-1$.
If $c = [r'-p,1^p]$ we have $u = [r'-p,1^{p-1}]$ and $v = [r'-p-1,1^p]$ ;
if $c = [r'-p],[1^p]$ we have $u = [r'-p-1],[1^p]$ and $v = [r'-p],[1^{p-1}]$. Note that, when $c = [r'-p,1^p]$
and $r' = 2p+1$, or when $c = ([r'-p];[1^p])$, then the part (2)
of this lemma has to be applied.

\begin{table}
$$
\begin{array}{|l|l|l|l|l|}
\hline
p & c & \rho & a & b  \\
\hline
\hline
2\leq p \leq r'-3 & ([r'-p,1^p]) & ([r'-p,1^p];[1]) & ([r'-p,1^{p-1}];[1])
 & ([r'-p-1,1^p];[1]) \\
\cline{3-5}
 &  & ([r'-p,2,1^{p-1}]) & ([r'-p,2,1^{p-2}]) & ([r'-p-1,2,1^{p-1}]) \\
\hline
2\leq p \leq r'-2 & ([r'-p];[1^p]) & ([r'-p,1];[1^p]) & ([r'-p-1,1];[1^p])
& ([r'-p,1];[1^{p-1}])  \\
\cline{3-5}
 &  & ([r'-p];[2,1^{p-1}]) & ([r'-p-1];[2,1^{p-1}]) & ([r'-p];[2,1^{p-2}]) \\
\cline{3-5}
 &  & ([r'-p];[1^p];[1]) & ([r'-p-1];[1^p];[1]) & ([r'-p];[1^{p-1}];[1])  \\
\hline
\end{array}
$$
\caption{Restrictions for $\Ind\LRef$ of $G(e,e,r)$, $r \geq 4$, $r'=r-1$.}
\label{tabindLRef}
\end{table}

We now assume that $c \in \QRef$, and that $\bla$ is a representation of
$G(e,1,r)$ whose restriction contains $c$. We can assume $r \geq 4$,
since the case of $G(e,e,3)$ has already been settled.
If $c \in \{ ([r'-1];[1]), ([r'-1,1]) \}$, then any $\rho$
with $c \nearrow \rho$ which is not a quasireflection representation
is handled by lemma \ref{lemres3} or lemma \ref{lemres2}, along
the patterns described in table \ref{tabindRef}, (1) and (2), provided
that $\rk \mathcal{H}_0(c) \geq 2$. This is the case for $r \geq 5$ by
the induction hypothesis. If $r = 5$ and $c$ is not of the above
types, either we have $A( \bla)  = 2$ and we are in the special case
for $G(2g,2g,4)$ already described, or $c = ([2,2])$ and we can use
lemma \ref{lemres2} along the pattern (3) of table \ref{tabindRef}.

If $r=4$, then $\rho$ is of one of types
listed in table \ref{tabindrefspec} (up to taking conjugate parts). The two
possible results for $\rk \mathcal{H}_0(\bar{\rho})$ given in the
table for $\rho = ([2,1];[1])$ depend on whether $\rho \simeq \rho^* \otimes
\eps$ (first column) or not. The lemmas of \S 3 show that $\mathcal{H}(\rho) \simeq
\mathcal{L}(\rho)$ in these cases.

\begin{table}
$$
\begin{array}{|l||l|l|l|l|l|}
\hline
 & r & c & X & a & b  \\
\hline
\hline
(1) & r' \geq 4 & ([r'-1];[1]) & ([r'-1,1];[1]) & ([r'-2,1];[1]) & ([r'-1,1]) \\
\cline{4-6}
&  & & ([r'-1];[2]) & ([r'-2];[2]) &  \\
\cline{4-6}
& & & ([r'-1];[1];[1]) & ([r'-2];[1];[1]) & ([r'-1]; \emptyset;[1]) \\
\hline
(2) & r' \geq 3 & ([r'-1,1]) & ([r'-1,2]) & ([r'-2,2]) &  \\
\hline
(3)&  & ([2,2]) & ([3,2]) & ([3,1]) &  \\
\hline
(4)&  & ([2]; [2]) ^{\pm} & ([3]; [2]) & ([3];[1]) & [2][2]^{\mp} \\
\hline
\end{array}
$$
\caption{Restrictions for $\Ind\QRef$ of $G(e,e,r)$, $r \geq 4$, $r' = r-1$.}
\label{tabindRef}
\end{table}

\begin{table}
$$
\begin{array}{|l||c|c|c|c|}
\hline
\rho & \multicolumn{2}{|c|}{([2,1];[1])} & ([2];[2]) & ([2];[1];[1]) \\
\hline
\dim \rho &  \multicolumn{2}{|c|}{9} & 6 & 12 \\
\hline
\rk \mathcal{H}_0(\bar{\rho}) & 4 & 6 & 4 & 9 \\
\hline
\end{array}
$$
\caption{Special cases.}
\label{tabindrefspec}

\end{table}

\subsection{Induction for the groups $G(e,e,r)$}

Here we prove the theorem for the groups $G(e,e,r)$. The
cases $r \in \{ 2, 3 \}$ have been done separately, so we can
start the induction at $G(e,e,r)$ for $r \geq 4$. Let $\rho \in
\Irr G(e,e,r)$,
a component of $\bla \in \Irr G(e,1,r)$. If $\rho \in \Ind \LRef$
we already proved that $\mathcal{H}(\rho) \simeq \mathcal{L}(\rho)$, so
we can assume $\rho \not\in \Ind \LRef$. If the restriction of $\rho$
to $G(e,e,r-1)$ is multiplicity-free we have $\mathcal{H}(\rho) \simeq
\mathcal{L}(\rho)$
by lemma \ref{lemindmult1}. 

It is easily checked that, if $r=4$, then the only possibility for
the restriction $\bar{\rho}$ of $\rho$ to have multiplicities
is when $e = 5g$ for some integer $g$ and $\bla$ can be taken
to be $\la^0 = \emptyset$, $\la^{g} = \la^{2g} = \la^{3g} = \la^{4g} = [1]$
with $p(\bla) = 4$. In that case, $\rho$ factorizes through a representation
of $G(5,5,4)$, and we check that $\mathcal{H}(\rho) \simeq \mathcal{L}(\rho)$
by computer for this representation. We can thus assume $r > 4$. This
implies that $G(e,e,r)$ does not admit irreducible representations
of dimensions 2 or 3 (see proposition \ref{propsmallrep}). By \cite{BRANCH}
theorem 1 we know that the multiplicities in $\bar{\rho}$ are at most 2. We need
a lemma.

\begin{lemma}
If $\rho$ is an irreducible representation of $G(e,e,r)$ with $r \geq 2$ whose
restriction $\bar{\rho}$ to $G(e,e,r-1)$ admits an irreducible component
with multiplicity 2, then $\bar{\rho}$ admits, counting multiplicities,
at least 3 irreducible components.
\end{lemma}
\begin{proof}
By \cite{BRANCH} theorem 3.3, the assumption implies that $\rho$ is
the restriction of some irreducible representation $\bla$ of $G(e,1,r)$.
Let $\rho_1$ be an irreducible component of $\bar{\rho}$ occuring
with multiplicity 2. By \cite{BRANCH} proposition 3.4, $\rho_1$
is the restriction of some $\bmu \nearrow \bla$, and there is some
$\bmu^+ \nearrow \bla$ with $\bmu \neq \bmu^+$ whose restriction
to $G(e,e,r-1)$ is isomorphic to $\rho_1$. We argue by contradiction. If
there were less than 3 irreducible components, then the restriction
of $\bla$ to $G(e,1,r-1)$ would be isomorphic to $\bmu \oplus \bmu^+$.
In particular, $p(\bla) \leq 2$. Moreover, since $\bmu^+$ should be a nontrivial
cyclic permutation of $\bmu$, one has $p(\bla) > 1$ hence
$p(\bla) = 2$. Letting $\bla = (a;b)$ with $a,b$ nontrivial partitions
at suitable places $0 \leq i<j\leq e-1$, we can assume $\bmu = (a^0 ; b)$,
$\bmu^+ = (a; b^0)$ with $a^0 \nearrow a$, $b^0 \nearrow b$. Since $\bmu$
and $\bmu^+$ are cyclically permuted one from the other, we have $2(j-i) = e$.
But this implies for $\bla = (a;a)$ that $A(\bla) = 2$, contradicting
the irreducibility assumption on $\rho$.
\end{proof} 

By lemma \ref{lemindmult2one} we have $\mathcal{H}(\rho) \simeq
\mathcal{L}(\rho)$ as soon as $\bar{\rho}$ has only one
component of multiplicity 2. We thus assume that $\bar{\rho}$ admits
at least 2 multiplicity-two components. Then every irreducible
component of $\bar{\rho}$ has dimension at least 12, by the following
lemma.

\begin{lemma} If $r > 4$ and $\bar{\rho}$ admits at least 2 multiplicity-two
components, then $\rho' \nearrow \rho \Rightarrow \dim \rho' \geq 12$.
\end{lemma}

\begin{proof}
$\rho'$ is a component of some $\bmu \nearrow \bla$, and
$\dim \rho' = (\dim \bmu)/A(\bmu)$. Since $\dim \bmu \geq p(\bmu) !$,
$A(\bmu)$ divides $p(\bmu)$ and $p(\bmu) \geq p(\bla) -1$, first
note that $\dim \rho' \geq (p(\bla) -2)!$.
By \cite{BRANCH} (see end of \S 3 and proposition 7.3 there) the
assumption on the multiplicity-two components of $\bar{\rho}$ implies
that $\bla$ can be chosen such that there exists $m$ dividing $e$ with
$u = e/m \geq 5$ and
\begin{enumerate}
\item $\la^i$ only depends on the class of $i$ modulo $m$ if $m$
does not divide $i$
\item $\la^0 = a$, $\la^i = b$ for $i$ a nonzero multiple of $m$,
with $a \nearrow b$.
\end{enumerate}
Condition (2) implies $p(\bla) \geq 4$. Moreover, if there exists $o$
with $m \not| i$ and $\la^i \neq \emptyset$, then condition (1)
implies $p(\bla) \geq 9$ hence $\dim \rho' \geq 7! \geq 12$. Assuming
that this is not the case, then $\rho$ factors through $G(u,u,r)$, so that
we can assume $m = 1$ hence $e \geq 5$. If $e \geq 7$ then $p(\bla) \geq 6$
and $\dim \rho' \geq 4! \geq 12$, so we can assume $e \in \{ 5,6 \}$.
If $a \neq \emptyset$ then $|b| \geq 2$ and by the branching rule
$\dim \bmu \geq \dim ([1];[1];[2];[2];[2]) = 5040$. Since $A(\bmu) \leq e \leq
6$ this implies $\dim \rho' \geq 5040/6 \geq 12$. We thus can
assume $a = \emptyset$, $b = [1]$, $r = e-1$ and $\dim \bmu = (e-2)!$.
Since $r > 4$ we have $e = 6$, $\dim \bmu = 24$. Now $A(\bmu)$
divides both $p(\bmu) = 4$ and $e = 6$ by lemma \ref{brisuregeer}
hence $A(\bmu) \leq 2$ and $\dim \rho' \geq 24/2 = 12$.
\end{proof}

Now write
$$
\bar{\rho} = \sum_{i=1}^r \alpha_i \rho_i + 
\sum_{i=1}^s \left( \beta_i \varphi_i + \gamma_i \varphi_i^* \otimes
\eps \right) + \sum_{i=1}^t \eta_i \psi_i
$$
with $\alpha_i, \beta_i, \gamma_i \in \{ 1, 2 \}$ and the
usual conventions on $\rho_i,\varphi_i, \psi_i$. Since $\dim \rho' \geq 3$
for $\rho' \nearrow \rho$ we have $\dim \rho' - 1 \geq \frac{\dim \rho'}{2}$
and
$ \rk \mathcal{H}_0(\bar{\rho}) \geq \dim soc(\bar{\rho})/2 -s$,
where
$$
soc(\bar{\rho}) =  \sum_{i=1}^r  \rho_i + 
\sum_{i=1}^s \left(  \varphi_i +  \varphi^* \otimes
\eps \right) + \sum_{i=1}^t  \psi_i
$$
hence
$$\rk \mathcal{H}_0(\bar{\rho}) \geqslant \frac{\dim \rho}{4} -s
> \frac{\dim \rho}{5}
$$
as soon as $s < \dim \rho/20$. In particular, since
$\dim \varphi_i \geq 11$ for all $i \in [1,s]$ we have
$\rk \mathcal{H}_0(\bar{\rho}) > \frac{\dim \rho}{5}$ and $\mathcal{H}(\rho)
\simeq \mathcal{L}(\rho)$ by lemma \ref{lemdim5}, since $\rk \mathcal{H}(\rho)
\geq \rk \mathcal{H}(\varphi_i) = \dim \varphi_i - 1 \geq 10$. This concludes
the proof of the theorem for the infinite series.

\subsection{Exceptional groups}

In order to conclude the proof of the theorem for $W \neq H_4$.
we need to show that
$\mathcal{H}(\rho) \simeq \mathcal{L}(\rho)$
for all irreducible representations not in $\LRef$ for these
groups. This holds true also for $W$ of Coxeter type $H_4$.

\begin{prop} \label{proplieoutlref} If $W$ is an irreducible exceptional reflection
group and $\rho \in \Irr(W) \setminus \LRef$, then
$\mathcal{H}(\rho) \simeq \mathcal{L}(\rho)$.
\end{prop}

This section is devoted to the proof of this proposition.
We first deal with groups of small rank, before using the
above induction lemmas to prove the proposition for
groups of higher rank.

\subsubsection{Groups of small rank} For these
groups, we used a computer to show that $\rho(\mathcal{H}')$ and
$\mathcal{L}(\rho)$ have the same dimension. For the groups
labelled $G_{12}$, $G_{13}$, $G_{22}$, $G_{23} = H_3$, the representations
have dimension small enough to do the computations directly over the
field of definition of these groups, and matrix models
of all representations are known. For some representations of the groups
$G_{24}, G_{27}$,
$G_{28}=F_4$, and $G_{30} = H_4$, we used reduction of the
coefficients. Notice that it is enough to show that
$\dim \rho(\mathcal{H}') \geq \dim \mathcal{L}(\rho)$.

When
$\rho$ admits a model over $\Q$, we consider a Lie subalgebra
$A_{\rho}$ of $\sl(V_{\rho})$ generated by nonzero multiples of the
$\rho(s'), s \in \mathcal{\rho}$, so that $A_{\rho}$ is defined
over $\Z$ and $\rk_{\Z} A_{\rho} = \dim_{\Q} \rho(\mathcal{H}')
= \dim_{\Q} \mathcal{H}(\rho)$. It is then sufficient to find
a prime $p$ such that $\dim_{\mathbbm{F}_p} A_{\rho} \otimes_{\Z} \mathbbm{F}_p \geq
\dim_{\Q} \mathcal{L}(\rho)$ to prove this inequality. In most
cases, it was enough to take $p = 11$. This solves the cases of
the groups $G_{24}$ and $G_{28} = F_4$. For $G_{24}$ we had to find
matrix models of 3 rational representations which were not available
before. We used for this decomposition of tensor products of
already known representations and methods \cite{MARINMICHEL}). These new models
have since been included in CHEVIE.

Matrix models for the representations of $H_4$ were found
by Alvis, and are included in CHEVIE. In the case of $G_{27}$
we computed matrix models for the representations not in $\LRef$
(now also included in CHEVIE). As predicted by a theorem of Benard
(see \cite{benard}), the nonrational representations among them
can be realized over a quadratic field, either $\Q(\sqrt{5})$ or
$\Q(\zeta_3)$. As in the rational case, the dimension of $\rho(\mathcal{H}')$
is equal to the rank of some Lie algebra $A_{\rho}$ defined
(by chasing denominators) over
the ring of integers of the corresponding quadratic field. It is then
enough to compute the dimension of the reduction of $A_{\rho}$
modulo a suitable place. Recall that, when $\Z[u]$ is a principal
quadratic ring and $\pi = \alpha + \beta u \in \Z[u]$ has norm a
prime integer $p$, then a morphism $\Z[u] \onto \mathbbm{F}_p$ can be defined
by sending $u$ to $- \alpha \beta^{-1}$. Using these
explicit reductions, we prove by computer that $\mathcal{H}(\rho) \simeq
\mathcal{L}(\rho)$ for these additional two groups.

\subsubsection{Groups of higher rank}

The remaining exceptional groups are labelled $G_{29},G_{31},G_{33},G_{34},
G_{35} = E_6, G_{36}=E_7$ and $G_{37} = E_8$ in the Shephard-Todd
classification. We use the
character and induction tables implemented in CHEVIE to check
the following properties. All of them have a
single class of reflections, and the dimensions greater than 1 of their
irreducible representations are at least 4. By lemma \ref{lempasrk1},
it follows that, for all $\rho \in \Irr(W)$, $\mathcal{H}(\rho)$
admits no simple ideal of rank 1. We choose for $W_0\subset W$
the subgroups described in table \ref{tabpatternexc}, where the generators
correspond to the tables of \cite{BMR}. 

\begin{table}
$$
\begin{array}{|c|c|c|}
\hline
\mbox{Group} & \mbox{Type of subgroup} & \mbox{Reflection subgroup} \\
\hline
G_{29} & G(4,4,3) & <t,u,v > \\
G_{31} & G_{29} & <s,t,v,w > \\
G_{33} & G(3,3,4) & <s,t,u,w> \\
G_{34} & G_{33} & <s,t,u,v,w> \\
E_6 & D_5 & <s_1,\dots,s_5> \\
E_7 & E_6 & <s_1,\dots,s_6> \\
E_8 & E_7 & <s_1,\dots,s_7> \\
\hline
\end{array}
$$\caption{Induction patterns for exceptional groups.}
\label{tabpatternexc}
\end{table}

Notice
that they all such $W_0$ are irreducible, that they have only one class
of reflections, and all differ from $F_4$ and $H_4$. 
Assuming that the theorem holds true for $W_0$, we check
that the condition $\dim V < (r+1)^2$ of lemma \ref{lemredsimpstrong}
holds for the pairs $\h = \rho(\mathcal{H}'_{W_0})$
and $\g = \rho(\mathcal{H}')$ for all $\rho \in \Irr(\rho)
\setminus \LRef$. Using the induction table, this is easily
done by computer. Now the induction
table for these pairs does not contain multiplicities more than 2.
Then lemma \ref{lemmultiplie} shows that condition (II) of lemma \ref{lemredsimpstrong} is
fulfilled, which proves that $\mathcal{H}(\rho)$ is a simple
Lie algebra for all $\rho \in \Irr(W) \setminus \LRef$.

We check that $\dim \rho <4\rk \rho(\mathcal{H}_0')$
for all $\rho \in \Irr(W) \setminus \LRef$. This proves
that $\mathcal{H}(\rho) \simeq \mathcal{L}(\rho)$
for all these by lemma \ref{caractBCD}, provided we are not in the exceptions
listed there. For the groups here of rang at least 6,
we check that $\rk \rho(\mathcal{H}_0') \geq 14 > 6$,
hence $\mathcal{H}(\rho) \simeq \mathcal{L}(\rho)$,
for all $\rho \in \Irr(W) \setminus \LRef$. This proves the theorem
for $G_{34},E_6,E_7,E_8$ by corollary \ref{corimpltheo}.

For $G_{29}, G_{31}$ and $G_{33}$ and $W_0$ as in the table,
the following situations appear, for $\rho \not\in \LRef$.
Here we marked ``*'' the cases where $\rho^* \otimes \eps \simeq \rho$.
$$
\begin{array}{|l|l|c|c|c|c|c|c|c|c|c|c|}
\hline
G_{29} & \mathrm{dim} & 5 & 6 & 10 & 10 & 15 & 15 & 16 & 20 \\
\cline{2-10}
 &\rk \rho(\mathcal{H}'_0) & 3 & 3 & 5 & 6 & 8 & 9 & 9 & 10 \\
\cline{2-10}
& \mathrm{sym.} & & & & & & & & \\
\hline
G_{31} & \mathrm{dim} & 5 & 9 & 10 & 15 & 16 & 20 & 20 & 20 \\
\cline{2-10}
 &\rk \rho(\mathcal{H}'_0) & 4 & 7 & 9 & 14 & 7 & 9 & 18 & 19\\
\cline{2-10}
& \mathrm{sym.}& & & & & * & * & &  \\
\hline
G_{33} & \mathrm{dim} & 6 & 15 & 15 & 20 & & & & \\
\cline{2-10}
 &\rk \rho(\mathcal{H}'_0) & 5 & 11 & 13 & 17 & & & & \\
\cline{2-10}
& \mathrm{sym.} & & & & & & & & \\
\hline
\end{array}
$$
We leave to the reader to check that none of these belong to
the exceptions of lemma \ref{caractBCD}, hereby proving $\mathcal{H}(\rho) \simeq \mathcal{L}(\rho)$
(and the theorem) for these groups.

\section{Closures inside Iwahori-Hecke algebras} 

Recall that $W \subset \GL_n(\C)$ be a finite
irreducible reflection group of rank $n$, and that we mean by
reflection an involutive element of $\GL_n(\C)$
distinct from the identity which fixes an hyperplane of $\C^n$. 

\subsection{Braid groups and holonomy Lie algebras}

Let $\mathcal{R}$ be the set of reflections of $W$, and $\mathcal{A}$
the associated (central) hyperplane arrangement in $\C^n$. We
let $X = \C^n \setminus \bigcup \mathcal{A}$ denote its
complement. The (generalized) braid group associated to $W$
is defined by $B = \pi_1(X/W)$. It fits into a short
exact sequence $1 \to P \to B \to W \to 1$ where $P = \pi_1(X)$.
To every $H \in \mathcal{A}$ we associate the closed 1-form
$\frac{d \alpha_H}{\alpha_H}$, where $\alpha_H$ is an arbitrary linear
form on $\C^n$ with kernel $H$.

We recall from \cite{KOHNOKOSZUL} the construction of the holonomy Lie algebra in this
setting. The transposed of the cup-product gives a morphism
$\delta : H_2(X,\Q) \to \Lambda^2(H_1(X,\Q))$. On the other hand,
$\Lambda^2(H_1(X,\Q))$ can be identified to
the homogenous part of degree 2 of the free Lie algebra on the vector
space $H_1(X,\Q)$. The \emph{holonomy Lie algebra} $\mathcal{T}$ of $X$ is defined
as the quotient of this Lie algebra by the image $\delta\left(H_2(X,\Q)\right)$
of $\delta$. It it thus a graded quadratic Lie algebra, defined
such that the class of $\om \wedge \om$ in $H^2(X,\Q)$ vanishes,
where $\om$ is the following closed form
$$
\om = \sum_{H \in \mathcal{A}} t_H \om_H \in \Omega_1(X) \otimes \mathcal{T}.
$$
By a classical result of Brieskorn \cite{BRIESKBOURB}
the algebra generated by the forms $\om_H$ for $H \in \mathcal{A}$
embeds in the cohomology algebra, hence $\om \wedge \om = 0$.

Let us consider the action $w.t_{w(H)}$ of $W$ on $\mathcal{T}$
and introduct the semidirect product $\mathfrak{B} = \C W \ltimes
\mathsf{U} \mathcal{T}$ where $\mathsf{U} \mathcal{T}$ denotes
the universal envelopping algebra of $\mathcal{T}$. This algebra
is graded, by $\deg(t_H) = 1$ and $\deg(w) = 0$ if $w \in W$.
We denote $\widehat{\mathfrak{B}}$ its completion with respect
to the grading, $A = \C[[h]]$ the ring of formal series in one
variable, $K_0 = \C((h))$ its field of fractions, $K = \bigcup_{n \geq 0}
\C((h^{\frac{1}{n}}))$ the field of (formal) Puiseux series. Recall
that $K$ is an algebraic closure of $K_0$.
Any (finite dimensional) linear representation $\rho : \mathfrak{B}
\to M_N(\C)$ defines a closed integrable 1-form
$\om_{\rho} = h \sum_{H \in \mathcal{A}} \rho(t_H) \om_H
\in \Omega^1(X) \otimes M_N(A)$. This 1-form is $W$-equivariant,
hence induces an integrable 1-form on $X/W$.

We now fix a base point $\underline{z} \in X/W$, and consider
the differential equation $d F = \om_{\rho} F$. The monodromy of this
equation gives a representation $R : B \to \GL_N(A)$.
We say that an element of $B = \pi(X/W)$ is
a braided reflection (`generator-of-the-monodromy' in \cite{BMR})
if it is the composition of a path $\gamma$ from
$\underline{z}$ to some (neighborhood of an) hyperplane $H$,
followed by a positive half turn around $H$ in $X$, and by $\gamma$
in the opposite direction. In particular the image in $W$ of such an
element is the reflection around $H$.

We let $\mathcal{T}_{\rho} = \rho(\mathcal{T})\otimes_{\C} A$,
and $\GL_N^0(K) = \exp \left( h \End(V_{\rho})\right) \subset \GL_N(A)$.
We will use the following facts, for which we refer to \cite{KRAMCRG}
and \cite{ASSOC}.

\begin{itemize}
\item $R(B) \subset \rho(W) \ltimes \GL_N^0(K)$ and $R(P) \subset \GL_N^0(K)$.
More generally, the monodromy of $\om_{\rho}$ along any path in $X$
belongs to $\exp h \mathcal{T}_{\rho}$. 
\item If $b \in B$ is a braided reflection around $H \in \mathcal{A}$,
then $R(b)$ is conjugated to $\rho(s) \exp(\ii \pi h \rho(t_H))$ by an element of
$\exp h \mathcal{T}_{\rho} \subset \GL_N^0(K)$
and $R(b^2) \equiv 1 + 2 \ii \pi h \rho(t_{H})$ modulo $h$.
\item If the restriction of $\rho$ to $\mathcal{T}$ is irreducible,
then the restriction of $R$ to $P$ is irreducible. 
\item The Lie algebra of the Zariski-closure $\overline{R(P)}$ of $R(P)$ in $\GL_N(K)$
contains $\rho(\mathcal{T}) \otimes K$.
\end{itemize}
We will also need the following facts, which are elementary consequences
of $R(P) \subset \exp h \mathcal{T}_{\rho}$, and of the fact that
$\mathcal{T}$ is generated by the $t_H, H \in \mathcal{A}$.
\begin{itemize}
\item If there exists a bilinear form $<\ , \ >$ on $V_{\rho}$
such that $\rho(t_H)$ is antisymmetric, then the induced bilinear
form on $V_{\rho} \otimes A$ is invariant under $R$
\item If the $\rho(t_H)$ have zero trace, then $R(P) \subset \SL_N(K)$,
and $\det R(g) = \det \rho(\pi(g))$ for all $g \in B$.
\end{itemize}
Finally, the correspondance $\rho \rightsquigarrow R$ is functorial (see
\cite{ASSOC}). In particular,
if $\rho^1,\rho^2$ are two representations of $\mathcal{T}$
and $R_1,R_2$ the corresponding monodromy representations of $P$,
then
\begin{enumerate}
\item $R_1 \otimes R_2$ is the monodromy represention of $\rho_1 \otimes \rho_2$.
\item $R_1 \simeq R_2$ iff $\rho_1 \simeq \rho_2$ (see \cite{ASSOC}, proposition 5).
\end{enumerate}


\subsection{Hecke algebras}
\label{pfzarheck}
The (generic) Iwahori-Hecke algebra $H_W(q)$ associated to $W$ can be defined as
the quotient of the group algebra $K B$ by the ideal generated
by all $(\sigma -q)(\sigma+q^{-1}) = 0$, for $\sigma$ a braided
reflection, with $q$ some transcendental constant in $K$. Here
we take $q = e^{\ii \pi h}$. We recall the following facts from
\cite{BMR}. It is known for Coxeter groups
and reflection groups in the infinite series, plus a few other ones,
that these algebras are semisimple with representations afforded
by the following monodromy construction, originally due to I. Cherednik :
given a representation $\rho$ of $W$, it can be extended to a representation
of $\mathfrak{B}$ by letting $t_H$ act as $\rho(s_H)$, where $s_H \in W$ is
the reflection corresponding to $H \in \mathcal{A}$ ; the monodromy
representation associated to $\rho$ factors through an irreducible
representation of $H_W(q)$ by the
remarks above. In particular we denote $S_{\eta} : B \to K^{\times}$
the 1-dimensional representation associated to $\eta \in \Hom(W,\{ \pm 1 \}$.

We say that the reflection groups $W$ for which this property has been proved
are \emph{tackled}.
Recall from \cite{BMR} that this property is conjectured to hold
for any reflection group.


We define $V_{\rho}^h = V_{\rho} \otimes_{\C} K$, $\SL(V_{\rho}^h) = \{ x \in \GL(V_{\rho}^h) \ | \ \det(x) = 1 \}$
and $\widetilde{\SL}(V_{\rho}^h) = \{ x \in \GL(V_{\rho}^h) \ | \ \det(x) = \pm 1 \}$.
Recall from \cite{BMR} that $B$ is generated by braided reflections. Let
$\sigma$ be such a braided reflection, and $s = \pi(\sigma) \in \mathcal{R}$.
Since $R(\sigma)$ is conjugated to $\rho(s) \exp(\ii \pi \rho(s))$,
if $\tr \rho(s) = 0$ for all $s \in \mathcal{R}$, then $R(P) \subset \SL(V_{\rho}^h)$,
$R(B) \subset \widetilde{\SL}(V_{\rho}^h)$ and the following diagram
commutes, with exact rows.
$$
\xymatrix{
1 \ar[r] & \SL(V_{\rho}^h) \ar[r] & \widetilde{\SL}(V_{\rho}^h) \ar[r]^{\mathrm{det}} & \{ \pm 1 \} \ar[r] & 1 \\
1 \ar[r] & P \ar[u] \ar[r] & B \ar[u] \ar[r]^{\pi} & W \ar[u]^{\det \rho} \ar[r] & 1 
}
$$

Notice that, if $\tr \rho(s) = 0$ for all $s \in \mathcal{R}$,
then, since $s^2 = 1$, we have $\det \rho(s) = (-1)^{\frac{\dim \rho}{2}}$ for all
$s \in \mathcal{R}$. In particular, $\det\rho = \un$ or
$\det \rho = \eps$ in this case. We let $B^2 = \Ker \eps \circ \pi$. We
have $P \subset B^2 \subset B$.
Recalling that the Lie algebra of $\overline{R(P)}$
contains $\rho(\mathcal{H})\otimes_{\C} K$,
the following is a consequence of proposition \ref{proprefsurj},
as $\GL(V_{\rho}^h)$ and $\SL(V_{\rho}^h)$ are connected algebraic groups.

\begin{prop} Let $\rho \in \mathrm{QRef}$. If $\exists s \in \mathcal{R} \ 
\tr \rho(s) \neq 0$ then $\overline{R(B)} = \overline{R(P)} = \GL(V_{\rho}^h)$.
Otherwise, this implies $\dim \rho = 2$ and we have $\overline{R(P)} =
\overline{R(B^2)} = \SL(V_{\rho}^h)$, $\overline{R(B)} = \widetilde{\SL}(V_{\rho}^h)$.
\end{prop} 

Similarly, the following is a consequence of proposition
\ref{proprefrep} (notice that the morphisms defined there
are $W$-equivariant, hence are isomorphisms of $\mathfrak{B}$-modules).
\begin{prop} Let $\rho \in \LRef \setminus \mathrm{QRef}$.
If $\dim \rho = 1$, then $R = S_{\rho}$ and
$\overline{R(P)}= \overline{R(B)} = K^{\times}$. If $\dim \rho > 1$,
that is $\rho = \eta \otimes \Lambda^k \rho_0$ with $\rho_0 \in \mathrm{Ref}$,
$\eta \in \Hom(W,\{ \pm 1 \})$, $\dim \rho_0 > k \geq 2$, then
we have $R \simeq S_{\eta}^{-k} \otimes \Lambda^k R_{0,\eta}$,
with $R_{0,\eta}$ the monodromy representation of $B$ associated to $\eta\otimes\rho_0 \in \Irr(W)$.
In particular, $\overline{R(P)} = \overline{R(B)} \simeq \GL(V_{\rho_0})$
\end{prop}

If $\rho \in \Irr(W)$ satisfies $\rho \simeq \rho^* \otimes \eps$, then
we denote $\OSP(V_{\rho}^h)$ the algebraic group over $K$ of direct (determinant one) isometries
of the bilinear form $(\ | \ )$ associated to
$\eps \into \rho^* \otimes \rho^*$. It
has
for Lie algebra $\osp(V_{\rho}) \otimes K$. By the remarks above,
$R(P) \subset \OSP(V_{\rho}^h)$. Recall from proposition \ref{propospinf} that, when $\eps \into S^2 \rho^*$,
then the corresponding symmetric bilinear form is hyperbolic, hence defined over
$\Q$.

We define $\widehat{\OSP}(V_{\rho}^h) = \{ x \in \GL_N(K) \ | \ \exists \alpha_1,\alpha_2 \in \{ \pm 1 \} \ \ x^{-1} = \alpha_1 x^+,
\det(x) = \alpha_2 \}$ where $x^+$ denotes the adjoint of $x$
with respect to $(\ | \ )$. We have $\OSP(V_{\rho}^h) \subset \widehat{\OSP}(V_{\rho}^h)$,
with index 2 if $\eps \into \Lambda^2 \rho^*$, and
index 2 or 4 if $\eps \into S^2 \rho^*$, depending on $\dim \rho$ modulo 4.
We denote $p : \widehat{\OSP}(V_{\rho}^h) \to \{ \pm 1 \}^2$ the
morphism defined by $(\alpha_1,\alpha_2)$.
Since every reflection in $\mathcal{R}$ is antisymmetric w.r.t.
$(\ | \ )$ and has trace 0, we have the following commutative
diagram
$$
\xymatrix{
\OSP(V_{\rho}^h) \ar[r] & \widehat{\OSP}(V_{\rho}^h) \ar[r]^p & \{ \pm 1 \}^2 \\
P \ar[u]_R \ar[r] & B \ar[u]_{R} \ar[r]_{\pi} & W \ar[u]^{(\eps, \det \rho)}
}
$$
and define $\widetilde{\OSP}(V_{\rho}^h) = p^{-1} ((\eps,\det \rho)(W))$.
By definition, $R(B) \subset \widetilde{\OSP}(V_{\rho}^h)$. Moreover,
since $\rho \simeq \rho^* \otimes \eps$ then $\tr \rho(s) = 0$
for all $s \in \mathcal{R}$ hence $\det \rho \in \{ \un, \eps \}$.
It follows that $\OSP(V_{\rho}^h)$ always has index 2 in $\widetilde{\OSP}(V_{\rho}^h)$
and that $R(B^2) \subset \OSP(V_{\rho}^h)$.
Denoting $p_1 : \widetilde{\OSP}(V_{\rho}^h) \to \{ \pm 1 \}$ the map $x \mapsto \alpha_1$,
this leads to the following commutative diagram with exact rows
$$
\xymatrix{
1 \ar[r] & \OSP(V_{\rho}^h) \ar[r] & \widetilde{\OSP}(V_{\rho}^h) \ar[r]^{p_1} & \{ \pm 1 \} \ar[r] & 1\\
1 \ar[r] & P \ar[u]_R \ar[r] & B \ar[u]_{R} \ar[r]_{\pi} & W \ar[u]^{\eps} \ar[r] & 1
}
$$

For any $\rho \not\in \LRef$, we proved in \S 5 that $\mathcal{H}(\rho) \simeq
\rho(\mathcal{H}') \simeq \mathcal{L}(\rho)$, which is isomorphic
to $\osp(V_{\rho})$ if $\rho^* \otimes \eps \simeq \rho$, and to $\sl(V_{\rho})$
otherwise. Also note that $\rho(\mathcal{H}) = \tr(\rho(\mathcal{H})) +
\rho(\mathcal{H}') = \tr(\rho(\mathcal{R})) + \rho(\mathcal{H}')$.
From this, the proof of the following proposition
is straightforward, and follows the lines of \cite{LIETRANSP}, theorem B.

\begin{prop} Let $\rho \not\in \LRef$. If $\rho^* \otimes\eps \simeq
\rho$, then $\overline{R(P)} = \OSP(V_{\rho}^h)$ and $\overline{R(B)}
= \widetilde{\OSP}(V_{\rho}^h)$. If $\rho^* \otimes \eps \not\simeq
\rho \otimes \eta$ for every $\eta \in \mathrm{X}(\rho)$, we have
the following cases :
\begin{enumerate}
\item if there exists $s \in \mathcal{R}$ with $\tr \rho(s) \neq 0$, then
$\overline{R(P)} = \overline{R(B)} = \GL(V_{\rho}^h)$
\item if $\forall s \in \mathcal{R} \ \tr \rho(s) = 0$, then $\overline{R(P)} = \SL(V_{\rho}^h)$.
If $\det \rho = \un$ then $\overline{R(B)} = \SL(V_{\rho}^h)$, otherwise
$\overline{R(B)} = \widetilde{SL}(V_{\rho}^h)$.
\end{enumerate}
\end{prop}

The last case to consider is when $\rho^* \otimes \eps \simeq
\rho \otimes \eta$ for some $\eta \in \mathrm{X}(\rho)$ with
$\eta \neq \un$. We define as before the bilinear form
$( \ | \ )$ associated to $\eps \otimes \eta \into \rho^* \otimes
\rho^*$ and $x \mapsto x^+$ the adjunction operation.
Since $\eta \neq 1$, there exists $s \in \RR$ with $\rho(s) \neq \pm 1$,
and in particular $\tr \rho(s) \neq 0$. Moreover, any
braided reflection associated to $s$ is mapped by $R$
to the scalar $\pm q^{\pm 1 }$ with $q = \exp(\ii \pi h)$,
its square to $q^{\pm 2}$, and
the subgroups generated by either image has $K^{\times}$
for algebraic closure.

Let $\widehat{\GOSP}(V_{\rho}^h)$ denote the algebraic
group $\{ x \in \GL(V_{\rho}^h) \ | \ x^+ x \in K^{\times} \}$,
$\varphi : \widehat{\GOSP}(V_{\rho}^h) \to K^{\times}$
the character $\varphi(x) = x^+ x$, and $d :
\widehat{\GOSP}(V_{\rho}^h) \to K^{\times}$ its determinant
character. Recall that $\dim V_{\rho}=2u$ is even, and let
$\psi = \varphi^u d^{-1} : \widehat{\GOSP}(V_{\rho}^h) \to K^{\times}$.
We define $\GOSP(V_{\rho}^h) = \Ker \psi$.

It is easily checked that $\GOSP(V_{\rho}^h)$ is connected,
with Lie algebra $K \oplus \osp(V_{\rho}) \otimes K$,
and contains both $K^{\times}$ and $\OSP(V_{\rho}^h)$. Every
braided reflection is conjugated to some $y = \rho(s) \exp(\ii
\pi h \rho(s)) \in \widehat{\GOSP}(V_{\rho}^h)$
by some element in $\OSP(V_{\rho}^h)$. If $\rho(s) \neq \pm 1$
we have $\tr \rho(s) = 0$ as before (since $\rho(s)$ is
then conjugated to $-\rho(s)$) hence
$\det y = \det \rho(s) = (-1)^u$, and $\varphi(y) = \varphi(
\rho(s)) = -1$, where $\psi(y) = 1$. If $\rho(s) = \pm 1$
then $d(y) = (\pm 1)^{2u} q^{\pm 2u}$ for $q = \exp(\ii \pi h)$,
and $\varphi(y) = q^{\pm 2}$, whence $\psi(y) = 1$. It follows that
$\overline{R(P)} \subset \overline{R(B)} \subset \GOSP(V_{\rho}^h)$
and by the same arguments as above we get that $\overline{R(P)}
= \overline{R(B)} = \GOSP(V_{\rho}^h)$.

\begin{prop} Let $\rho \not\in \LRef$ with $\rho^* \otimes
\eps \simeq \rho \otimes \eta$ for some $\eta \in \mathrm{X}
(\rho) \setminus \{ \un \}$. Then 
$\overline{R(P)}
= \overline{R(B)} = \GOSP(V_{\rho}^h)$.
\end{prop}

The four propositions above together imply theorem 2 of the introduction
(notice that the algebraic groups considered here are defined over $K_0$
and have a dense set of $K_0$-points). It is worth noticing
that all situations above actually occur, most of them
already for $W = \mathfrak{S}_n$ (see \cite{LIETRANSP}).



If $\rho \not\in \LRef$, we let $G(\rho) = \OSP(V_{\rho}^h)$,
$\tilde{G}(\rho) = \widetilde{\OSP}(V_{\rho}^h)$
if $\rho^* \otimes \eps \simeq \rho$, 
$\tilde{G}(\rho) = \GOSP(V_{\rho}^h)$ if $\rho^* \otimes \eps
\simeq \rho \otimes \eta$ for some $\eta \in \mathrm{X}(\rho) \setminus
\{ \un \}$,
and $\tilde{G}(\rho) = \GL(V_{\rho}^h)$
otherwise. We assume that $W$ is \emph{tackled}.
Then $H_W^{\times}$ is naturally isomorphic to
$\prod_{\rho \in \Irr(W)} \GL(V_{\rho}^h)$. Under this isomorphism
the results above prove that the morphism $B \to H_W^{\times}$
factors through the algebraic group
$$
\prod_{\rho \in \Hom(W,\{ \pm 1 \})} K^{\times} \times \prod_{\rho \in
\overline{\mathrm{QRef}}} \GL(V_{\rho}^h) \times
\prod_{\rho \in (\Irr(W) \setminus \LRef)/\approx}
\tilde{G}(\rho)
$$
We let $\mathcal{R}/W = \{ \mathcal{R}_1,\dots,\mathcal{R}_k \}$
and denote $\eta_i : W \to \{ \pm 1 \}$ such that $\eta_i(\mathcal{R}_j)
= \{ -1 \}$ if $j=i$ and $\{ 1 \}$ if $j \neq i$. We let $\chi_{\rho}^i = (\dim \rho - \tr
\rho(s_i))/2$ for $s_i \in \mathcal{R}_i$.
If $\rho \not\in \LRef$ is such that $\rho^* \otimes \eps
\simeq \rho \otimes \eta$ for some $\eta \in \mathrm{X}(\rho)
\setminus \{ \un \}$, we notice that $\eta$ is
uniquely determined by $\eta(s) = -1 \Leftrightarrow \rho(s)
\in \{ \pm 1 \}$, and we write $\eta = \prod_{i \in I_{\rho}} \eta_i$
for some nonempty $I_{\rho} \subset \{ 1, \dots, k \}$. For
$i \in I_{\rho}$ we let $\rho(\RR_i)$ denote the common value
of the $\rho(s)$ for $s \in \RR_i$.

We introduce the subgroup $\check{G}$ of tuples $(a_{\rho})$ as
above such that
$$
(a_{\un} a_{\eps})^2 = 1, \ \ \ \det a_{\rho} = a_{\un}^{\dim \rho} \prod_{i=1}^k (a_{\eta_i} a_{\un}^{-1})^{\chi_{\rho}^i},
$$
and
$$
\phi(a_{\rho}) = (a_{\eps} a_{un})^{-1} \prod_{i \in I_{\rho}} 
(a_{\eta_i} a_{\un}^{-1})^{- \rho(\RR_i)}
$$
for $\varphi : \GOSP(V_{\rho}^h) \to K^{\times}$ defined
above.
Since a braided reflection $\sigma$ with $\pi(\sigma) = s$ has
image conjugated to $\rho(s) \exp(\ii \pi h \rho(s))$, and
since $B$ is generated by braided reflections, it is easily checked
that the image
in $H_W^{\times}$ of $B$ lies inside $\check{G}$.
There is a morphism $p_{\eps} : \check{G} \to \{ \pm 1 \}$
given by $(a_{\rho}) \mapsto a_{\un} a_{\eps}$, and also
one $p_{\rho} : \check{G} \to \{ \pm 1 \}$
for each $\rho \not\in \LRef$ such that $\rho^* \otimes \eps \simeq \rho$,
which is inherited from $p_1 : \widetilde{\OSP}(V_{\rho}^h) \to \{ \pm 1 \}$.
We define $\tilde{G} \subset \check{G}$ by
$$
\tilde{G} = \{ g \in \check{G} \ | \ p_{\eps}(g) = p_{\rho}(g), \rho \not\in \LRef, \rho^* \otimes \eps \simeq \rho \}
$$
and denote $\pi_{\eps} : \tilde{G} \to \{ \pm 1 \}$ the morphism induced by $p_{\eps}$.
The image of $B$ is contained in $\tilde{G}$, while the image of $B^2$ is contained in
$G = \Ker \pi_{\eps}$.

The following follows from theorem \ref{maintheo}, and is a
generalization to (almost all) irreducible
complex reflection groups of theorem C of \cite{LIETRANSP}.

\begin{prop} \label{propadhHW} 
Assume that $W$ is an irreducible complex reflection
group which is \emph{tackled} with $W \neq H_4$. The image of $B$ is Zariski-dense in $\tilde{G}$. The connected component
of $\tilde{G}$ is $G$, which contains the images of $B^2$ and $P$ as
Zariski-dense subgroups. The short exact sequence
$1 \to G \to \tilde{G} \to \{ \pm 1 \} \to 1$
defined by $\pi_{\eps}$ is split.
\end{prop}
\begin{proof}
Since $\tilde{G}$ contains the image of
$B$ and $G$ contains
the image of $B^2 \supset P$, it is sufficient to show
that the Zariski-closure of the image of $P$ has the same Lie algebra as $G$,
as $G$ is clearly connected. Actually,
since the former Lie algebra is obviously included in the latter, one only needs
to know that these two Lie algebras are isomorphic. This is the content
of theorem \ref{maintheo}, which 
concludes the proof of the main part of the proposition. The
extension is split by $-1 \to (\rho(s))_{\rho}$, where $s$ is some reflection
is $s$. Indeed, letting $\sigma$ denote a braided reflection with
$\pi(\sigma) = s$, its image in $H_W^{\times}$ belongs to $\tilde{G} \setminus G$ and
is conjugated to the collection of $(\rho(s) e^{\ii \pi h \rho(s)})_{\rho}$ by some element of $G$.
By checking the definition of $G$ we see that
$e^{\ii \pi h \rho(s)}$ belongs to $G$, hence $(\rho(s))_{\rho}$
belongs to $\tilde{G} \setminus G$.
\end{proof}


\subsection{The case of $W = H_4$}

The decomposition of $\mathcal{H}$ is similar to the other cases,
except that one has to decompose the factor $\mathcal{H}(\rho)$
for $\rho = \rho^a + \rho^b$ where $\rho^a, \rho^b$ are the two
8-dimensional irreducible representations of $W$. We have $\rho^i \not\in \LRef$, $\eps \into S^2 \rho^i$,
hence $\mathcal{H}(\rho^i) \simeq \rho^i (\mathcal{H})' \simeq \so_8(\k)$
for $i \in \{ a,b \}$.
Let $\CW$ be a maximal parabolic subgroup of $W$ of Coxeter type $H_3$.
We have $\Res_{\CW} \rho^a \simeq \varphi + \varphi^* \otimes \eps \simeq
\Res_{\CW} \rho^b$
where $\varphi \in \Irr(\CW) \setminus \LRef(\CW)$ is
4-dimensional, defined over $\Q$ and such that
$\varphi \not\simeq \varphi^* \otimes \eps$.

\subsubsection{The Lie algebra $\rho(\mathcal{H}')$}
 In particular,
$\varphi(\mathcal{H}'_0) \simeq \sl_4(\k)$. We thus can choose
matrix models for $\rho^a,\rho^b$ in hyperbolic $\so_8(\k)$ with
\emph{equal} restrictions to $\mathcal{H}'_0$, such that
$(\varphi + \varphi^* \otimes \eps)(\mathcal{H}'_0)$ is an
isotropic $\sl_4$ inside $\so_8$. By computer we check
that $ \dim \rho(\mathcal{H}') = \dim \so_8(\k)$. Since
we proved $\mathcal{H}(\rho^i) \simeq \so_8(\k)$ in proposition \ref{proplieoutlref}, the following
follows readily from the simplicity of the Lie algebra $\so_8(\k)$.

\begin{lemma} \label{lemlieH4} Under the models above, $\rho(\mathcal{H}') = 
\{ (x, \psi(x) \ | \ x \in \so_8(\k) \}$ for some $\psi \in \Aut(\so_8(\k))$
which is not induced by $\GL_8(\k)$-conjugation and which pointwise stabilizes the isotropic $\sl_4(\k)$.
\end{lemma}

The fact that $\psi$ is not a conjugation automorphism is derived from
the fact that $\rho^a_{\mathcal{H}'} \not\simeq \rho^b_{\mathcal{H}'}$.
The existence of such an automorphism is reminiscent from the triality
phenomenon. We describe it in more (matrix) detail. Assuming that the
quadratic form is in hyperbolic shape, elements of $\so_8(\k)$
have the form $\begin{pmatrix} s+a & b \\ c & -s- ^t a \end{pmatrix}$,
with $s \in \k, a \in \sl_4(\k)$ and $b,c \in \so_4(\k)$, where $\so_4(\k)$
designates skew-symmetric matrices. We leave to the reader to check that
any automorphism $\psi$ satisfying the conditions of the lemma
has the form
$$
\psi \left( \begin{pmatrix} s + a & b \\ c & -s- ^t a \end{pmatrix}
\right) = \begin{pmatrix} - s + a & \la^{-1} c' \\ \la b'  &  s - ^t a \end{pmatrix}
$$
where $\la$ is some fixed nonzero scalar, and $M \mapsto M'$ is defined by
$$
\begin{pmatrix}
0 &  b_{1} &  b_{2} &  b_{3} \\-b_{1} &  0 &  b_{4} &  b_{5} \\-b_{2} &  -b_{4} &  0 &  b_{6} \\-b_{3} &  -b_{5} &  -b_{6} &  0
\end{pmatrix}
\mapsto
\begin{pmatrix}
0 &  b_{6} &  -b_{5} &  b_{4} \\-b_{6} &  0 &  b_{3} &  -b_{2} \\b_{5} &  -b_{3} &  0 &  b_{1} \\-b_{4} &  b_{2} &  -b_{1} &  0\end{pmatrix}
$$
In particular, $\psi^2 = \Id$, and
the matrix model can be chosen such that $\la = 1$.

\subsubsection{Zariski closure of $R(P)$}

For $x \in \SO_8(K)$, we denote $\bar{x}$ its image
in $\PSO_8(K)$. We know compute the Zariski closure $\Gamma \subset \SO_8(K) \times
\SO_8(K)$ of $R(P)$.
We denote $\Gamma_i \simeq \SO_8(K)$ the Zariski closure of $R_i(P)$.
We let $p_i$ denote the projection of $\SO_8(K)\times \SO_8(K)$ onto
its $i$-th factor. It is a morphism of algebraic groups, hence it maps
$\Gamma$ onto a closed subgroup of $\Gamma_i$ and its restriction
$\hat{p}_i$ to $\Gamma$ is also a morphism of algebraic groups. Since
$R(P)$ is Zariski-dense in $\Gamma$ it follows that $R_i(P) = \hat{p}_i\circ
R(P)$ is Zariski-dense in $p_i(\Gamma)$ hence $p_i(\Gamma) = \Gamma_i$.

We now let $\Gamma^i= \Ker \hat{p}_j \subset \Gamma$ for $\{i,j \} = \{ 1,2 \}$.
Then $p_i$ identifies $\Gamma^i$ with a normal subgroup of $\Gamma_i$, wich
is connected. Since $\Gamma$ and the $\Gamma_i, \Gamma^i$
are linear algebraic groups, we have a natural isomorphism
of algebraic groups $(\Gamma/\Gamma^i)/(\Gamma^i \Gamma^j/\Gamma^i)
\simeq (\Gamma/\Gamma^j)/(\Gamma^i \Gamma^j/\Gamma^j)$ (see e.g.
\cite{SPRINGER} 5.5.11). Since by definition $\Gamma_j = \Ker \hat{p}_i$
we have $\Gamma/\Gamma^j \simeq p_i(\Gamma) = \Gamma_i$
and similarly $\Gamma/\Gamma^i \simeq \Gamma_j$. Finally, since
$\Gamma^i \cap \Gamma^j = \{ 1 \}$ we have $\Gamma^i \Gamma^j / \Gamma^j
\simeq \Gamma^i$. It follows that this natural isomorphism
defines an isomorphism $\Psi : \Gamma_1/\Gamma^1 \to \Gamma_2/\Gamma^2$.

Since $\SO_8(K)$ does not admit any proper connected normal subgroup,
the connected component of $p_1(\Gamma^1)$ is either trivial
or equal to $\Gamma_1$. The latter case implies
$p_1(\Gamma^1) = \Gamma_1$, $\Ker \hat{p}_1 = \{ 1 \} \times \SO_8(K)$.
Since $p_1(\Gamma) = \SO_8(K)$ this implies $\Gamma = \SO_8(K) \times \SO_8(K)$.
In this case, $R_1 \otimes R_2$ would be irreducible under the action
of $P$, meaning that the tensor product $(\rho^a)_{\mathcal{T}} \otimes (\rho^b)_{\mathcal{T}}$ would be irreducible
under the action of the holonomy Lie algebra $\mathcal{T}$.
A direct computation shows that this is not the case, for instance because
$\dim (\rho^a)_{\mathcal{T}} \otimes (\rho^b)_{\mathcal{T}}(\mathsf{U}\mathcal{T}) < 64^2$.
If follows that the connected components of $p_1(\Gamma^1)$, and
similarly of $p_2(\Gamma^2)$, are trivial. It follows that each $\Gamma^i$
is a subgroup of $Z(\Gamma_i) = \{ \pm 1 \}$. By $\Gamma_1/\Gamma^1 \simeq
\Gamma_2 / \Gamma^2$ it the same subgroup for $i \in \{ 1 , 2 \}$, as
$\SO_8(K) \not\simeq \PSO_8(K) = \SO_8(K) / \{ \pm 1 \}$.

We prove that this subgroup is $\{ \pm 1 \}$, because otherwise $\Psi$
would induce an automorphism of $\SO_8(K)$ ; since automorphisms
of $\SO_8(K)$ are conjugation by some element in $\mathrm{O}_8(K)$,
this would provide an intertwinner between the representations $R_1,R_2$,
hence also between the representations $(\rho^a)_{\mathcal{H}},
(\rho^b)_{\mathcal{H}}$. By proposition \ref{propisorep} this implies that $\rho^a,\rho^b \in \Irr(W)$
are isomorphic, which is not the case.

In particular $\Psi$ is an automorphism $\PSO_8(K) \simeq \Gamma_1 / \{ \pm 1 \}
\to \Gamma_2 / \{ \pm 1 \} \simeq \PSO_8(K)$, and $\Gamma =
\{ (x,y) \in \SO_8(L)^2 \ | \ \bar{y} = \Psi(\bar{x}) \}$.
This automorphism cannot be induced by $\mathrm{PO}_8(K)$, otherwise
there could exist $x \in \GL_8(K)$ such taht $x R_2(g) x^{-1} = 
\pm R_1(g)$ for all $g \in P$, whence $x \rho^b(g) x^{-1} = \pm
\rho^a(g)$ for all $g \in \mathcal{R}$, hence for all $g \in W$.
Then $\eta : W \to \{ \pm 1 \}$ defined by $x \rho^b(g) x^{-1}
= \eta(g) \rho^a(g)$ is a character of $W$ such that $\rho^b \simeq \rho^a
\otimes \eta$. But $\Hom(W,\{ \pm 1 \}) = \{ \un, \eps \}$
and $\rho^b \not\simeq \rho^a, \rho^b \not\simeq \eps \otimes \rho^a$,
a contradiction. It is thus a triality automorphism. Moreover, the
matrix models of $\rho^a,\rho^b$ have been chosen such that $\rho^a(w)
= \rho^b(w)$ for all $w \in W_0$. It follows that $\Psi$ pointwise
stabilizes the isotropic $\SL_4(K)$, hence the induced automorphism
$\psi$ of $\Aut(\so_8(K))$ belongs to the ones determined by the lemma
above, and has order 2. Since $\PSO_8(K)$ is connected it
follows that $\Psi^2 = \Id$.

From this description and the classification of reductive algebraic groups,
it is clear that $\Gamma$ is isomorphic to $\Spin_8(K)$, as $\Gamma$ is connected with Lie algebra $\so_8(K)$
and has at least the 4 elements $(\pm 1, \pm 1 )$ in its center. We
thus proved the following.

\begin{prop} There exists $\Psi \in \Aut(\PSO_8(K))$
inducing $\psi \in \Aut(\so_8(K))$, such that
the Zariski-closure of $R(P)$ is
$\{ (x , y) \in \SO_8(K)^2 \ | \ \bar{y} = \Psi(\bar{x}) \}$.
This group is isomorphic to $\mathrm{Spin}_8(K)$.
\end{prop}

\subsubsection{Zariski closure of $R(B)$}
Let $\widetilde{\SO}_8(K) = \SO_8(K) \times \{ \pm 1 \}$, with $\pi_0 : \widetilde{\SO}_8(K)$
the projection onto the second factor. Note that $\widetilde{\OSP}(\rho^i) \simeq \widetilde{\SO}_8(K)$.
We let $\widetilde{\PSO}(K) = \PSO_8(K) \times \{ \pm 1 \}$. One clearly has
$R(B) = \{ (x,y) \in \widetilde{\SO}_8(K)^2 \ | \ \pi_0(x) = \pi_0(y) \}$.
The surjective morphism $\eps \circ \pi : B \to \{ \pm 1 \}$ can thus be extended
to a morphism $\tilde{\Gamma} \to \{ \pm 1 \}$, restriction of $\pi_0 \times \pi_0$ to $\tilde{\Gamma}$, $\tilde{\Gamma}$ designates
the Zariski closure of $R(B)$. We thus have a short exact sequence $1 \to \Gamma \to \tilde{\Gamma} \to \{ \pm 1 \} \to 1$,
which is split, by sending $-1$ to $(\rho^a(s),\rho^b(s))$ for some $s \in \mathcal{R}$, by the same
argument as in proposition \ref{propadhHW}. This extension is not trivial, and moreover we have $Z(\tilde{\Gamma}) = Z(\Gamma)$.
Indeed, by irreducibility of $R_1,R_2$ under the action of $P$, central elements in $\tilde{\Gamma} \subset \widetilde{\SO}_8 \times
\widetilde{\SO}_8(K)$ have
the form $(\la,\mu)$ for $\la,\mu \in K$, hence belong $Z(\widetilde{\SO}_8(K)^2) = Z(\Gamma)$, hence $Z(\tilde{\Gamma}) = Z(\Gamma)$.
We thus proved the following.

\begin{prop} Let $\Gamma, \tilde{\Gamma}$ denote the Zariski closures of $R(B),R(P)$, respectively.
Then $\tilde{\Gamma}$ is a split extension of $\{ \pm 1 \}$ by $\Gamma$ with
$Z(\tilde{\Gamma}) = Z(\Gamma)$.
\end{prop}

Notice that this extension, which is uniquely up to isomorphism by the action of $(s,s)$ on $\Gamma$,
for $s \in \mathcal{R}$, is not isomorphic (as a group) to the usual $\mathrm{Pin}_8(K)$ extension, because the centers
of $\mathrm{Pin}_8(K)$ and $\Spin_8(K)$ do not coincide.

\subsection{Proof of theorem 1 for $H_4$}

Theorem 1 of the introduction is known for $W \neq H_4$, as
it is a consequence of theorem \ref{maintheo} proved in \S 5.
We thus assume $W = H_4$.
Let $\approx'$ be defined as the introduction, meaning
$\rho^a \approx' \rho^b$ for the special representations
above and $\rho^1 \approx \rho^2 \Rightarrow \rho^1 \approx' \rho^2$. 
Recall the morphism $\Phi$ from section 2. By lemma \ref{lemlieH4} is factorizes through
$$
\Phi^+ : \mathcal{H}' \to \mathcal{I} \oplus \left(
\bigoplus_{\rho \in (\Irr(W)\setminus \Lambda\mathrm{Ref})/\approx} \mathcal{L}(\rho) \right)
$$
and we only need to prove that $\Phi^+$ is surjective. We know that
$\mathcal{H}(\rho) \simeq \mathcal{L}(\rho)$ for all $\rho \not\in \LRef$,
and one checks through the character table that, for $\rho \in \Irr(W)$,
$\mathcal{H}(\rho) \simeq \so_8(\k) \Leftrightarrow \rho \in \{ \rho^a, \rho^b \}$.
Then the arguments of proposition \ref{propidsimpl} and corollary \ref{corimpltheo} show that
$\Phi^+$ is surjective, which concludes the proof of theorem 1.

Finally, proposition \ref{propadhHW} implies the following for
$W \neq H_4$, the case $W = H_4$ being now straightforward.

\begin{theo} \label{theozarhecke}
Assume that $W$ is an irreducible complex reflection group which
is \emph{tackled}. The Zariski closures of the image of $P$ and $B^2$
in $H_W(q)^{\times}$ coincide, they are connected and reductive,
with Lie algebra $\mathcal{H} \otimes K$. They have index 2 in the Zariski closure of
the image of $B$, and the corresponding extension is split.
\end{theo}

It can be noticed that the Zariski closure of (the image of) $P$
is simply connected exactly for the groups $W$ which do not admit
$\rho \in \Irr(W) \setminus \LRef$ with $\eps \into S^2 \rho$.
This is the case for all irreducible groups of rank 2 by using :
lemma \ref{lemexcdim2} , the fact that $G(e,1,2)$ as well as its reflection
subgroups admit no irreducible representation of dimension greater than 2,
and inspection of the character table of the groups $G_{12}$, $G_{13}$, $G_{22}$ and
$\mathfrak{S}_3$.
Similarly, the groups $G(e,e,3)$ and $G(2e,e,3)$ have irreducible
representations of dimension $1,2,3$ or $6$. By lemma \ref{lemcasA4B4} below we cannot
have $\eps \into S^2 \rho$ with $\dim \rho = 6$ in rank 3, and the other
representations belong to $\LRef$ by lemma \ref{lemref23}.
Checking the exceptional groups and $\mathfrak{S}_4$ separately, it follows
that the corresponding algebraic groups are simply connected if $\rk W \leq 3$,
while many non-simply connected cases appear in rank 4 (e.g. $F_4, G_{29}, H_4,G_{31}$, etc.).

\subsection{Unitarisability in Coxeter types}


We assume that $W$ is a Coxeter group, and for
$q \in \C^{\times}$ we denote $\hat{H}(q)$ the
specialized (complex) Iwahori-Hecke algebra associated
to $q$. For $q$ close enough to 1, $\hat{H}(q)$
is isomorphic to the group algebra $\C W$, and
there is a natural 1-1 correspondance between its
irreducible representations and the ones of $W$
(we use \cite{GP} as our main reference on these topics).
We prove that, for $q$ close enough to 1 these representations
are unitarizable. This extends previous results
in type A (see \cite{WENZL}, and \cite{GT}).

The main idea in the proof below is a combination of deformation and descent
arguments along the following patterns.

$$
\xymatrix{ 
 &\C \ar[dl]_{z \mapsto \bar{z}} & &\R[q,q^{-1}] \ar[rr] \ar[ll] \ar[dl]_{q \mapsto q^{-1}} & & \R(q)  \ar[dl]_{q \mapsto q^{-1}}  \\
\C & & \R[q,q^{-1}] \ar[rr] \ar[ll]& & \R(q)  \\
& \R \ar[uu]|\hole \ar[ul]& &  \R[q+q^{-1}]\ar[rr] \ar[ll] \ar[uu]|\hole \ar[ul]& & \R(q+q^{-1})  \ar[uu]|\hole \ar[ul] \\
}
$$
{}
$$
\xymatrix{ 
 & \R(q) \ar[rr]^{q \mapsto e^{\mathrm{i} \pi h}} \ar[dl]_{q \mapsto q^{-1}} & & \C((h)) \ar[dl]_{h \mapsto -h} \\
 \R(q) \ar[rr]& & \C((h)) \\
 & \R(q+q^{-1}) \ar[rr] \ar[uu]|\hole \ar[ul]& & \C((h^2)) \ar[uu] \ar[ul] \\
}
$$

\begin{prop} \label{propunit} For $q \in \C$ with $|q| = 1$ and close enough to 1, every representation
of $\hat{H}(q)$ is unitarizable as a representation of $B$.
\end{prop}

\begin{proof} We can assume that the representation
considered here is irreducible. We let $Q_0 = \R(q)$
the field of rational fractions in one indeterminate
over $\R$. By \cite{GP} \S 9.3.8 this representation
is a specialization of $R_0 : B \to \GL_N(A_0)$
where $A_0 = \R[q,q^{-1}]$ is the ring of Laurent polynomials,
and the specialization $q=1$ provides the corresponding
representation of $W$. Since $W$ is a Coxeter group, up
to conjugation by an element of $\GL_N(\R)$ we can assume that this representation of $W$
is orthogonal w.r.t. the standard scalar product of $\R^N$.

We let $\eps_0 \in \Gal(Q_0|\R)$ defined by $\eps_0(q) = q^{-1}$
and $\eps \in \Aut(K_0)$ defined by $f(h) \mapsto f(-h)$.
These automorphisms stabilized $A_0$ and $A$, respectively. The
ring morphism $A_0 \into A$ defined by $q \mapsto \exp(\ii \pi h)$ induces
an embedding $Q_0 \into K_0$ equivariant under $\eps_0$ and $\eps$. The
ring morphism is compatible with the reductions at $q = 1$ and $h = 0$, meaning that
the following diagram is commutative.
$$
\xymatrix{
A_0 \ar[r] \ar[d]_{q=1} & A \ar[d]^{h=0} \\
\R \ar[r] & \C 
}
$$
By extension of scalars we deduce from $R_0$ a representation $K_0 \otimes R_0 : B \to
\GL_N(K_0)$. Since it is a representation of the generic Hecke algebra over
$K_0$, it is isomorphic to a monodromy representation $R : B \to \GL_N(A)$,
that can be chosen with $R(B) \subset U_N^{\eps}(A)$, where $U_N^{\eps}(A)$
denotes the formal unitary group associated to $\eps$, and $R_{h=0} = (R_0)_{q=1}$
(this follows from \cite{KRAMCRG}, prop. 2.4, since the reflections are both orthogonal and selfadjoint
w.r.t. the $W$-invariant (standard) scalar product chosen here). It follows that
there exists $P \in \GL_N(K_0)$ such that $R(g) P = P R_0(g)$ for all $g \in B$.
Up to multiplication by some power of $h$, we can assume that $P$ has coefficients in $A$
and that its reduction $\overline{P} \in Mat_N(\C)$ is nonzero. It follows that $\overline{P}$
is an intertwinner between the two absolutely irreducible representations $R_{h=0}$ and $R_{q=1}$.
Since these two representations are actually equal, $\overline{P}$ is some scalar,
that can be assumed to equal 1, up to multiplication by some real scalar, hence $P
\in \GL_N(A)$.

Let $b_1,\dots,b_r$ a set of generators for $B$. The mapping
$$
\Phi : J \mapsto \left( \eps_0( ^t R_0(b_i)) J - J R_0(b_i)^{-1}
\right)_{i=1..r}
$$
from $Mat_N(Q_0)$ to $Mat_N(Q_0)^r$ is $Q_0$-linear. The dimension of its kernel is thus
equal to the one of $K_0 \otimes \Phi : Mat_N(K_0) \to Mat_N(K_0)$. Considering
the corresponding linear system $\eps_0( ^t R_0(b_i)) J - J R_0(b_i)^{-1} = 0$
in $Mat_N(K_0)$ we see that it is in 1-1 correspondance with the space of $B$-module
morphisms between $(K_0 \otimes R_0)^{\eps}$ and 
$K_0 \otimes R_0^*$, where $R_0 : b \to \ ^tR(b^{-1})$ is the
dual representation of $R_0$ and $(K_0 \otimes R_0)^{\eps}$
maps $b \mapsto \eps(R_0(b))$.
These
two representations being absolutely irreducible, by Schur lemma this space has dimension at most
1, and a corresponding nonzero $J$ is necessarily invertible. Since $R$ is isomorphic to $R_0$
over $K_0$, such a $J$ exists. Indeed, from $R(b_i) = P R_0(b_i)P^{-1}$
and $ ^t \eps(R(b_i)) = R(b_i)^{-1}$,
we get $\eps( ^t P^{-1}) \eps ( ^t R_0(b_i) ) \eps( ^t P)
= P R_0(b_i)^{-1} P^{-1}$, meaning
$\eps ( ^t R_0(b_i) ) \eps( ^t P) P
= \eps( ^t P) P R_0(b_i)^{-1}$, and this is the desired
equation with
$J = J' = \eps( ^t P) P$. Since $\overline{P} = 1$ we have $J'_{h=0} = 1$.

It follows that there exists a nonzero $J_0 \in M_N(Q_0)$ solving
this linear system. Up to multiplication of $J_0$ by powers
of $q-1$, we can assume
$J_0 \in Mat_N(Q_0) \cap Mat_N(A)$, and moreover
$(J_0)_{h=0} \neq 0$. Since the reductions at $q=1$ of $R_0$
and $\eps_0 \, ^t R_0$ are (absolutely) irreducible and equal, we can actually
assume $(J_0)_{h=0} = 1$, up to rescaling $J_0$ by $(J_0)^{-1}_{h=0} \in \R$.

On the other hand, there exists $\mu \in K^{\times}$ with $J' = \mu J_0$.
Since $J'$ and $J_0$ belong to $\GL_N(A)$ 
with $(J')_{h=0} = (J_0)_{h=0} = 1$,
we get $\mu \in A$ and $\mu \equiv 1 \not\equiv 0$ modulo $h$.
As a consequence we get that $1+\eps(\mu)\mu^{-1} \in A$
is invertible. We let $c$ denote its inverse.
From $\eps(J') = ^t J'$ we get $\eps(\mu) \eps(J_0) = 
\mu ^t J_0$. 
Since $J_0 \neq 0$ we have $\eps(\mu) \mu^{-1} \in Q_0$
hence $c \in Q_0 \cap A$.

Let now $J = 2cJ_0 \in \GL_N(Q_0) \cap M_N(A)$. We check $^t J = \eps_0(J)$,
and $J_{h=0} = 2 \frac{1}{2} (J_0)_{h=0} = 1$.
It follows that $J$ can be specialized in a neighbourhood of $q = 1$ and,
when $|q| = 1$, it defines a unitary form which is positive definite
when $q$ is close enough to $1$, which concludes the proof.
\end{proof}

\subsection{Explicit unitary models : the example of $D_4$}

We make the remark that the proof given above is constructive,
in the sense that, starting from a given matrix model over $\R[q,q^{-1}]$,
or even $\R(q)$, one can algorithmically derive a matrix $J$ with coefficients in $\R(q)$
that affords a unitary structure for $|q| = 1$ and $q$ close to 1.
Indeed, the obtention of $J_0$ is merely a linear algebra matter
from matrix models for $R_0$ (plus chasing $(q-1)$ denominators) ;
moreover, since $\eps(\mu) \eps(J_0) = \mu J_0$ we
get the value of $\eps(\mu) \mu^{-1}$ from any nonzero coefficient
$m_{ij}$ of $J_0$ as $\eps(\mu)\mu^{-1} = \eps_0(m_{ji})m_{ij}^{-1}$
and get $c$ explicitely.

As an example, we carry out this procedure on the reflection representation
of the Coxeter group $D_4$. A convenient matrix model is given
by Hoefsmit matrices, with Artin generators $\sigma_1,\sigma_2,\sigma_3,\sigma_4$
following the convention that $\sigma_1,\sigma_2, \sigma_4$ commute one to the other.

$$
\sigma_1 \mapsto \begin{pmatrix} q & 0 & 0 & 0 \\ 0 & q & 0 & 0  \\
0 & 0 & \frac{q-q^{-1}}{2} & - \frac{q+q^{-1}}{2} \\
0 & 0 & - \frac{q+q^{-1}}{2} & \frac{q-q^{-1}}{2} 
\end{pmatrix}
\sigma_2 \mapsto \begin{pmatrix} q & 0 & 0 & 0 \\ 0 & q & 0 & 0  \\
0 & 0 & \frac{q-q^{-1}}{2} &  \frac{q+q^{-1}}{2} \\
0 & 0 &  \frac{q+q^{-1}}{2} & \frac{q-q^{-1}}{2} 
\end{pmatrix}
$$
{}
$$
\sigma_3 \mapsto \begin{pmatrix}
q & 0 & 0 & 0 \\ 0 & \frac{q^2-1}{q^3+q} & \frac{2q}{1+q^2} & 0 \\
0 & \frac{1+q^4}{q+q^3} & \frac{q^3-q}{1+q^2} & 0 \\
0 & 0 & 0 & q
\end{pmatrix}
\sigma_4 \mapsto \begin{pmatrix} \frac{q^2 - 1}{q+q^5} & \frac{q+q^3}{1+q^4} & 0 & 0 \\
\frac{q^6+1}{q^5+q} & \frac{q^5 - q^3}{q^4+1} & 0 & 0 \\
0 & 0 & q & 0 \\ 0 & 0 & 0 & q
\end{pmatrix}
$$
Using the above procedure, we find for $J$ the diagonal matrix
$1, q^2 + 1+q^{-2}, \frac{(q^2+q^{-2})(q^2 + q^{-2} - 1)}{2},\frac{(q^2+q^{-2})(q^2 + q^{-2} - 1)}{2}$,
which satisfies $\eps_0(J) = ^t J$. Specializing at a complex number with $|q| = 1$ written
$q = \exp(\ii x)$ with $|x| \leq \frac{\pi}{2}$, we get that the form has signature $(4,0)$
for $|x| < \frac{\pi}{6}$, $(1,3)$ for $\frac{\pi}{6} < |x| <  \frac{\pi}{4}$
and $(3,1)$ for $\frac{\pi}{4} < |x|$.

\subsection{Topological closure for Coxeter groups}

We still assume that $W$ is an irreductible Coxeter group. The
representations of the Hecke algebras are actually defined over
$\overline{\Q}_0[q,q^{-1}]$, where $\overline{\Q}_0 = \overline{\Q} \cap \R$. Like the definition
over $\R[q,q^{-1}]$, this is a consequence of the existence
of W-graphs for all Coxeter groups (see \cite{GP}). Obviously,
every $\rho \in \Irr(W)$ is defined over $\overline{\Q}_0$.

We now consider a transcendantal number $u \in \C$
with $|u| = 1$ and $u$ close enough to 1, so that $\hat{R}$
is unitarizable by proposition \ref{propunit},
and investigate the topological closure of $\hat{R}(P)$.
It is a compact (Lie) subgroup of $U_N$, hence lies inside $G(V_{\rho}) \cap U_N$,
where $G(V_{\rho})$ is the Zariski closure of $\hat{R}(P)$.

By irreducibility of the action of $P$, this group $G(V_{\rho}) \cap U_N$ is a maximal
compact subgroup of $G(V_{\rho})$ ; for, every intermediate
compact subgroup should be contained in some unitary group, but the image of $P$
acting irreductibly can be contained in only one such group. In particular,
this is a real form of $G(V_{\rho})$. The fact that it is the
topological closure of $\hat{R}(P)$ is an immediate consequence of the
following well-known lemma.
\begin{lemma}
If $H_1 \subset H_2 \subset U_N$ is an inclusion of (compact) closed
subgroups, then $H_1 = H_2$ if and only if they have the same Zariski
closure in $\GL_N(\C)$.
\end{lemma}
\begin{proof}
If $H_1 \neq H_2$, then there exists a continuous $f : U_N \to \C$
with $f(H_1) = \{ 0 \}$ and $f(g H_1) = \{ 1 \}$ for
some $g \in H_2 \setminus H_1$, that can be chosen $H_1$-invariant
by averaging w.r.t. some Haar measure on $H_1$. We denote $(z_{ij})$
the complex coordinates of $Mat_N(\C)$.
By the Stone-Weierstrass theorem we get a polynomial $P$ in
the $z_{ij}$ and $\bar{z}_{ij}$ with $|f - P| \leq 1/4$ and averaging
$P$ we get another polynomial $\tilde{P}$ which takes distinct constant values
over $H_1$ and $g H_1$, so we can assume $\tilde{P}(H_1) = 0$
and $\tilde{P}(gH_1) = 1$.

Now the identity $\bar{M} = \ ^t M^{-1}$ for $M \in U_N$ implies that
$\tilde{P}$ coincides over $U_N$ with some $Q/det^r$, where
$Q$ is a polynomial in the $(z_{ij})$ and $det$ is the determinant. In
particular, $Q$ vanishes on $H_1$ but not on $g H_1$, a contradiction
if $H_1$ and $H_2$ had the same Zariski closure.


\end{proof}

We choose an isomorphism $\psi : \C \stackrel{\sim}{\to} K$ preserving $\overline{\Q}$
that maps our transcendantal number $u$ to $e^{\ii \pi h}$,
by identifying both $\C$ and $K$ to the algebraic closure
of the field of rational fractions with coefficients in $\overline{\Q}$
over a continuous number of indeterminates, one of them being $u \in \C$ and
$e^{\ii \pi h} \in K$. Now $G(V_{\rho})$ is actually defined over
$\overline{\Q}_0(u)$ (by \cite{BOREL} prop 1.3 b) and
$R$ is isomorphic to the twisting of $\hat{R}$ by $\psi$.
It follows that the Zariski closure of $R(P)$ is deduced from
$G(V_{\rho})$ by $\psi$ followed by a conjugation in $\GL(V_{\rho}^h)$.
In particular $G(V_{\rho})$ is reductive and connected, and so
are its maximal compact subgroups. It follows that the topological
closure of $\hat{R}(P)$ is connected.


Let $\g$ denote the (complex) Lie algebra of the Lie or algebraic group
$G(V_{\rho})$. Through
$\psi$ we have $\g \otimes_{\C,\psi} K
\simeq \rho(\mathcal{H})\otimes_{\C} K$, where $\otimes_{\C,\psi}$
denotes the identification of $\C$ with $K$, so that $\C \otimes_{\C,\psi} K$
is naturally isomorphic both to $\C$ and to $K$. Indeed, the latter Lie algebra
$\rho(\mathcal{H}) \otimes K$
is the Lie algebra of the Zariski closure of $R(P)$ by the results
of section \ref{pfzarheck},
and $R$ is isomorphic to the twisting of $\hat{R}$ by $\psi$.
Let $\mathcal{H}^{\overline{\Q}_0} \subset \overline{\Q}_0 W$ denote the infinitesimal
Hecke algebra with coefficients in $\overline{\Q}_0$.
Then $\rho(\mathcal{H}) = \rho(\mathcal{H}^{\overline{\Q}_0})
\otimes_{\overline{\Q}_0} \C$, and, since $\psi$ is the identity
on $\overline{\Q}$,
$$
\g \otimes_{\C,\psi} K
\simeq (\rho(\mathcal{H}^{\overline{\Q}_0}) \otimes_{\overline{\Q}_0}\C)
\otimes_{\C,\psi} K.
$$
which implies $\g \simeq \rho(\mathcal{H}^{\overline{\Q}_0}) \otimes_{\overline{\Q}_0}\C
\simeq \rho(\mathcal{H})$, because $\psi$ is an isomorphism (or because
$\C$ is algebraically closed). Since the same arguments work
for $B^2$, we proved the following.

\begin{theo} \label{theotophecke} For $u \in \C^{\times}$ a transcendental number,
$|u| = 1$ and $u$ close to 1, if $\rho$ is an irreducible real
representation of (the irreducible finite Coxeter group) $W$,
then the closure of the image of $P$ (or $B^2$) in
the corresponding representation of $\widehat{H}(u)$ is connected and
has for Lie algebra a compact real form of $\rho(\mathcal{H})$.
\end{theo}

Finally, we recall from proposition \ref{compactform} that a compact real form of $\mathcal{H}$
is given by the real Lie algebra $\mathcal{H}^c$ generated inside
$\C W$ by the elements $\ii s$, $s \in \mathcal{R}$. The compact
real form of the theorem is thus isomorphic to $\rho(\mathcal{H}^c)$.

\section{Technical results on representations of $W$}
\subsection{Small-dimensional representations}

\begin{lemma} The irreducible representations $\bla$ of
$G(e,1,r)$ of dimension at most 8 satisfy $p(\bla) \leq 3$.
The possible $\la_{i_k} \neq \emptyset$ for $ 1 \leq k \leq p(\bla)$,
with the $i_k$ distinct indices are given by table \ref{tabsmallrep}.
\end{lemma}
\begin{proof}
If $p(\bla) = 1$, according to the branching rule the dimension of $\bla$
is the same as the dimension of the symmetric group $\mathfrak{S}_{r+1}$
associated to $\la$. It is well-known that, besides 1-dimensional
representations, $\mathfrak{S}_{r+1}$ admits no irreducible
representations of dimension less than $r$, except when $r = 3$.
The conclusion follows by case-by-case examination for $r \leq 8$.

If $p(\bla) = 2$, then $\dim(\bla) = \dim (\la,\mu)$,
where $(\la,\mu)$ is a representation of $G(2,1,r)$, that is the
Coxeter group of type $B_r$. Once again, it is sufficient to
examine the case $r \leq 9$, and the conclusion follows by case-by-case
examination.

Now assume $p(\bla) \geq 3$. We know $\dim(\bla) \geq p(\bla)!$,
hence $p(\bla) =3$. Moreover the binary representations
with $p(\bla)=3$ have dimension $3! = 6$. If $\bla$ is not
a binary representation, then $r \geq 4$ and its restriction to $G(e,1,4)$ contains
a representation of the same dimension as $([2],[1],[1],\emptyset,\dots)$,
whose dimension is 12. The conclusion follows.

\end{proof}
\begin{table}
$$
\begin{array}{|c|c|c|l|l|l|}
\hline
p(\bla) & r & \dim \bla & \la^{i_1}\mbox{ or }(\la^{i_1})' & \la^{i_2}\mbox{ or }(\la^{i_2})' & \la^{i_3}\mbox{ or }(\la^{i_3})' \\
\hline
1 & r \geq 1 & 1 & [r] & & \\
\cline{2-6}
 & 2 \leq r \leq 9 & r-1 & [r-1,1] & & \\
\cline{2-6}
 & 4 & 2 & [2,2] & & \\
\cline{2-6}
 & 5 & 5 & [3,2] & & \\
\cline{2-6}
 & 6 & 5 & [3,3] & & \\
\cline{2-6}
 & 5 & 6 & [3,1,1] & & \\
\hline
2 & 2 \leq r \leq 8 & r & [r-1] & [1] & \\
\cline{2-6}
 & 4 & 6 & [2] & [2] & \\
\cline{2-6}
 & 5 & 8 & [3,1] & [1] & \\
\hline
3 & 3 & 6 & [1] & [1] & [1] \\
\hline
\end{array} 
$$
\caption{Irreducible representations of $G(e,1,r)$ of dimension at most 8}
\label{tabsmallrep}
\end{table}

\begin{table}
$$
\begin{array}{|c|c|c|c|c|l|l|l|}
\hline
r & \dim \rho & e & A(\bla) & p(\bla) & \la^{i_1}\mbox{ or }(\la^{i_1})' & \la^{i_2}\mbox{ or }(\la^{i_2})' & \la^{i_3}\mbox{ or }(\la^{i_3})' \\
\hline
3 & 2 & * & 1 & 1 & [2,1] & & \\
\cline{3-8}
 & & 3|e  & 3& 3  & [1] & [1]& [1]\\
\cline{2-8}
 & 3 & * & 1 & 2 & [2] & [1] & \\
\hline
4 & 2 & * & 1 & 1 &[2,2] & & \\
\cline{2-8}
 & 3 & * & 1 & 1& [3,1] & & \\
\cline{3-8}
 & &  2 | e& 2 & 2 & [2] & [2] & \\
\hline
\end{array} 
$$
\caption{Irreducible representations of $G(e,e,r)$ of dimension 2 or 3 for $r \geq 3$.}
\label{tabsmallrep2}
\end{table}

We will need the following result to check the inductive
assumptions in case of multiplicity 2 components.

\begin{prop} \label{propsmallrep} If $r \geq 5$ then the dimension of all irreducible
representations of $G(e,e,r)$ differ from 2 and 3.
If $r = 2$ all irreducible representations of $G(e,e,r)$ have dimension 1 or 2.
If $r \in \{ 3 , 4\}$ and $\rho$ is an irreducible representation
of $G(e,e,r)$ of dimension 2 or 3 then $\rho$ is a component of some of the representations
$\bla$ of $G(e,1,r)$ listed in table \ref{tabsmallrep2}.
\end{prop}

\begin{proof}
If $r = 2$ the group $G(e,e,r)$ is a dihedral group and the result is
classical.
Assume that $r > 2$ and $\dim \rho \in \{2,3\}$.
If $A(\bla) = 1$ the result follows from table \ref{tabsmallrep},
hence we assume $A(\bla) \neq 1$.

From the inequality $\dim \rho \geq (p(\bla)-1)!$ of lemma \ref{lemineqsmall}
we get $p(\bla) \leq 3$. Moreover $1 \neq A(\bla) | p(\bla)$ by
lemma \ref{brisuregeer} hence $A(\bla) = p(\bla) \in \{2,3 \}$.
If $A(\bla) = p(\bla) < 3$ or $\dim \rho < 3$ then $\dim \bla = A(\bla) \dim \rho \leq 8$
and the result follows by inspection of
table \ref{tabsmallrep}. In case
$A(\bla) = p(\bla) = 3$, then there exists $i$ such that $\la^i\not\in \{ \emptyset,
[1] \}$, otherwise $\dim \rho \geq p(\bla)! = 6$ by lemma \ref{lemineqsmall}.
It follows that $\bla$ is binary with $r = p(\bla) = 3$, which completes
the proof.

\end{proof}

\subsection{Proof of lemma \ref{lemexcmrd}}

We note the following facts. If $W = G(e,e,r)$
with $r \geq 3$, since $W$ admits a single class of reflexions
we have $\XX(\rho) = \{ \un \}$ as soon as $\dim \rho > 1$.
The same holds for $W = G(2e,e,r)$ with $r \geq 3$
and $\dim \rho > 1$ as soon as the restriction of $\rho$ to
$W_0 = G(2e,2e,r)$ is \emph{not} irreducible : we cannot have
$\rho(s_i) = \pm 1$ since $\dim \rho > 1$, and $\rho(t^e) = \pm 1$
would contradict the irreducibility of $\rho$. Finally note that,
if this restriction $\rho_0$ \emph{is} irreducible (with $\dim \rho > 1$)
then $\rho^* \otimes \eps \simeq \rho \otimes \eta$ for some $\eta \in \XX(\rho)$ implies
that $\rho_0^* \otimes \eps \simeq \rho_0$.

If $r =2$, recall that all
$\rho \in \Irr(W)$ have dimension 2, hence we do not have
to consider this case in the sequel.

\begin{lemma} \label{lemcasF4} There exists $\rho \in \Irr(W)$ with $\dim \rho = 4$ and
and $\eta \in \XX(\rho)$ with $\eps \otimes \eta \into S^2 \rho^*$
if and only if $W$ is a Coxeter group of type $F_4$.
\end{lemma}
\begin{proof}
We first check using the character tables in CHEVIE that $F_4 = G_{28}$ is the
only exceptional type admitting such a representation. We then assume
by contradiction that $W$ has type $G(e,e,r)$, $r \geq 3$ and admits such a
$\rho$, that is $\eps \into S^2 \rho^*$.
Let $\bla$ be an irreducible representation of $G(e,1,r)$
such that $\rho$ embeds in the restriction of $\bla$. If $A(\bla) = 1$
then $\dim \bla = 4$ and we can use table \ref{tabsmallrep} to check
that there are no possibilities with $\bla^* \otimes \eps \simeq \bla$.
It follows that $A(\bla) \neq 1$. By lemma \ref{lemineqsmall} we have
$4 \geq (p(\bla)-1)!$ hence $p(\bla) \leq 3$. From
$1 \neq A(\bla) | p(\bla) \leq 3$ it follows that $A(\bla) = p(\bla) \in 
\{2,3\}$. If $\bla$ was binary we would have $p(\bla)! = \dim \bla = A(\bla)
\dim \rho = p(\bla) \dim \rho$ hence $4 = \dim \rho = (p(\bla)-1)!$,
a contradiction. Then $\bla$ is not binary, $4 \geq p(\bla)!$ by lemma
\ref{lemineqsmall} and $p(\bla) = A(\bla) = 2$. But then $\dim \bla = 8$
and $A(\bla) = 2$, a case excluded by table \ref{tabsmallrep}.

Now we assume that $W$ has type $G(2e,e,r)$ and that $\rho$ embeds in some
$\bla \in \Irr(G(2e,1,r))$. If the restriction of $\rho$
to $G(2e,2e,r)$ is irreducible then we are reduced to the previous
case. Otherwise it has two 2-dimensional components $\rho^+,\rho^-$.
It follows that $2 | A(\bla)$. Then $2=\dim \rho^+  \geq (p(\bla)-1)!$
implies $p(\bla) \leq 3$, and $2 | A(\bla)|p(\bla)$ implies
$A(\bla) = p(\bla) = 2$. Then $\dim \bla = 4$ contradicts $A(\bla) = 2$
by table \ref{tabsmallrep}, and this concludes the proof. 
\end{proof}

\begin{lemma} \label{lemcasA4B4} Assume there exists
$\rho \in \Irr(W)$ with $\dim \rho = 6$
and $\eta \in \XX(\rho)$ with $\eps \otimes \eta \into S^2 \rho^*$.
Then $\rho \in \LRef$, $W$ has rank 4, and $\rho$
factorizes through a representation of either
a Coxeter group of type $A_4,B_4,F_4,H_4$, or an exceptional
group of type $G_{29},G_{31}$. 
\end{lemma}
\begin{proof}
We first assume that $W = G(e,e,r)$ and take $\bla$
corresponding to $\rho$. We first exclude the case $A(\bla) \neq 1$,
by contradiction. From $6 \geq (p(\bla) -1)! $
we would get $p(\bla) \leq 4$. Then $p(\bla) = 4$
implies that $\bla$ is binary, $\dim \bla = 24$ and $A(\bla) = 4$.
Then $\rho$ factorizes through the representation $([1],[1],[1],[1])$
of $G(4,4,4)$, which we check to be of symplectic type. It
follows that $p(\bla) \leq 3$. Since $1 \neq A(\bla) | p(\bla)$
we have $p(\bla) = A(\bla) \in \{2,3 \}$. Then $p(\bla) = 3$
implies $\dim \bla = 18$; since $\bla$ is not binary
it follows that $\bla$ has dimension greater than $\dim ([2],[2],[2]) = 36$,
a contradiction. It follows that $p(\bla) =A(\bla) = 2$ and
$\dim \bla = 12$. Then $\bla$ has two non-empty parts of the same
shape $\la$, and we can assume $\la \supset [2]$. If $\la = [2]$
then $\dim \bla = 6 < 12$, and otherwise we can assume
either $\la \supset [2,1]$ or $\la \supset [3]$. But $\dim ([3],[3]) = 20>12$
and $\dim ([2,1],[2,1]) = 84>12$, a contradiction.

We thus have $A(\bla) = 1$ and $\dim \bla = 6$. From table \ref{tabsmallrep}
the additional condition $\bla^* \otimes \eps \simeq \bla$
implies that either $\bla = (\la,\emptyset,\dots)$
with $\la = [3,1,1]$, in which case $\rho$ factors through the
alternating square of the reflection representation of $\mathfrak{S}_5$,
or $p(\bla) = 2$, with three nonempty parts $\la_0=\la_g= \la_{e-g}=[1]$
with $e \neq 3g$. We consider the explicit models (see e.g.
\cite{MARINMICHEL}) for
this representation of $G(e,1,3)$, generated by $t,s_1,s_2$. For $\{i,j,k\} = \{ 1,2,3 \}$
we denote $e_{ijk}$ the tableau with $(i,j,k)$ placed in position
$(0,g,e-g)$. Then it is easily checked that
$x = e_{312}\wedge e_{321} + e_{123}\wedge e_{132} - e_{213} \wedge
e_{231}$ is fixed by $t$, and $s_i .x = -x$ for $i \in \{ 1 ,2 \}$.
Since $G(e,e,3)$ is generated by $s_1^t,s_1,s_2$ it follows
that $\eps \into \Lambda^2 \rho$, hence $\eps \not\into S^2 \rho^*$.

There remains to consider the case $W = G(2e,e,r)$. If the
restriction to $G(2e,2e,r)$ is irreducible we are done, otherwise
$\XX(\rho) = \{ \un \}$ hence $\rho^* \otimes \eps \simeq \rho$,
and moreover it is the direct sum of $\rho^+$ and $\rho^-$ of dimension 3.
By table \ref{tabsmallrep2} we have $p(\bla) = A(\bla) = 2$, and
the condition $\rho^*\otimes \eps \simeq \rho$ implies
that $\bla$ can be taken of the form $(\la_0,\emptyset,\dots,\la_{e-1},\emptyset,\dots)$
with $\{\la_0,\la_{e-1}\} = \{[2],[1,1]\}$. Then $\rho$
factorizes through either representations $([2],[1,1])$ or $([1,1],[2])$
of the Coxeter group of type $B_4$, which are the alternating square
of its reflection representation, and its tensor product by a
multiplicative character, respectively. It follows that all
these representations belong to $\LRef$. The assertion on
the exceptional groups is based on a case-by-case computer check.
\end{proof}

Finally we will prove that, for the groups $G(e,e,r)$ and $G(2e,e,r)$,
the case 
$\eps \otimes \eta \into S^2 \rho^*$.
for some $\eta \in \XX(\rho)$ with $\dim \rho = 8$ does not occur.
We first need a lemma.

\begin{lemma} \label{lemrestr8}The 8-dimensional irreducible representations of $G(e,e,r)$
are restrictions of a representation $\bla$ of $G(e,1,r)$ such that
$A(\bla) = 1$. The 8-dimensional irreducible representations of $G(2e,e,r)$
are restrictions of a representation $\bla$ of $G(2e,1,r)$ such that
$A(\bla) = 1$.
\end{lemma}
\begin{proof}
Let $\rho$ be a 8-dimensional irreducible representation of $G(e,e,r)$.
There exists a well-defined irreducible representation $\bla$ of $G(e,1,r)$
such that $\rho$ embeds in the restriction of $\bla$. Assume
by contradiction that $A(\bla) \neq 1$.

By lemma \ref{lemineqsmall} we know that $8 = \dim \rho \geq  (p(\bla)-1)!$,
hence $p(\bla) \leq 4$ and $p(\bla) \leq 3$ if $\dim \rho \neq 8$.
First assume that $\bla$
is binary. If $p(\bla) = 4$, then $\dim(\bla) = 24$
and $A(\bla) | 4$. Since $A(\bla) \neq 1$ then $2 | A(\bla)$
and $\dim( \rho) \in \{ 6,12 \}$, which contradicts $\dim(\rho)  =8$.

We thus can assume that $\bla$ is not binary. Then there exists
$\bmu \nearrow \bla$ with $p(\bmu) = p(\bla)$, hence
$\dim(\rho) \geq p(\bla)!$ and $p(\bla) \leq 3$. The case $p(\bla) = 1$
is ruled out by $A(\bla) \neq 1$. If $p(\bla) = 3$ that is $A(\bla) = 3$,
then $\dim(\rho) = 24$. Moreover, there exists $i,j,k$ distincts
and $\la \neq \emptyset$ such that $\la^i = \la^j = \la^k$. Since
$\bla$ is not binary, we can assume $\la \supset [2]$. But
$\dim([2],[2],[2]) = 36 \geq 24$, hence a contradiction. It
follows that $p(\bla) = 2$, $A(\bla) = 2$ and $\dim(\bla) = 16$.
Let $i,j$ distincts and $\la \neq \emptyset$ such that $\la = \la^i = \la^j$.
Since $([2],[2])$ has dimension 12 it follows that
$\la \supset [2,1]$ hence $\dim(\bla) \geq \dim([2,1],[2,1]) = 84 > 16$, a contradiction.

Now let $\rho$ be a irreducible representation of $G(2e,e,r)$,
and $\bla$ a irreducible representation of $G(2e,1,r)$ of which
it is a component. If its
restriction to $G(2e,2e,r)$ is irreducible we are reduced to the former
situation. Otherwise this restriction is the sum of two irreducible
components $\rho^+$ and $\rho^-$, and $2|A(\bla)$ hence $2 | p(\bla)$.
Since $4 = \dim \rho^+ \geq (p(\bla)-1)!$ we have $p(\bla) \leq 3$,
hence $p(\bla)=2$, $A(\bla) = 2$, $\dim \bla = 8 = \dim \rho$, which
implies $A(\bla) = 1$. We thus get a contradiction and the conclusion.

\end{proof}

\begin{prop} \label{lempasuncannyeer} Assume that $W$ has type
$G(e,e,r)$ or $G(2e,e,r)$, and let $\rho \in \Irr(W)$.
If $\dim \rho = 8$, then $S^2 \rho^*$ does not contain $\eps \otimes \eta$
for any $\eta \in \XX(\rho)$.
\end{prop}
\begin{proof}
In both kind of types for $W$,
such a representation $\rho$ has to be the restriction
of an 8-dimensional irreducible representation of some $G(d,1,r)$
by lemma \ref{lemrestr8}. 
In the case of $G(2e,e,r)$, the restriction of $\rho$ to
$G(2e,2e,r)$ has to be irreducible by the same lemma and
satisfy the same assumption, so we only
have to deal with the $G(e,e,r)$ case, and we can
assume by contradiction that $\eps \into
S^2 \rho^*$.

By the previous lemmas such a representation $\rho$ would be the restriction
of a 8-dimensional irreducible representation $\bla$ of $G(e,1,r)$ with $A(\bla) = 1$.
In order that $\rho^* \otimes \eps \simeq \rho$ the set $\{ \la^i \ | \ 
i \in [0,e-1] \}$ has to be stable under $\la \mapsto \la'$. The only
possibility in table \ref{tabsmallrep} is that $p(\bla) = 2$ and that there exists $i,j$ such that
$\la^i = [2,1]$ and $\la^j = [1]$. In particular $r = 4$.

We thus can assume $\la^0 = [2,1]$, and $\la^i = [1]$ with $i \neq 0$.
Then $\rho^* \otimes \eps \simeq \rho$ implies $\la^{[-i]} = (\la^i)'$
hence $e$ is even and $i = e/2$. From the matrix formulas (see e.g.
\cite{MARINMICHEL} formula 3.3) such a representation factorizes
through $G(e,e,4) \onto G(2,2,4)$ and more precisely is induced by the
restriction to the Coxeter group $D_4=G(2,2,4)$ of
$([2,1],[1])$. It is
easily checked that $\eps$ embeds in the alternating square of this
representation, not in the symmetric one. This concludes
the proof.
\end{proof}

In order to conclude the proof of lemma \ref{lemexcmrd},
we then examine (using CHEVIE) the irreducible 8-dimensional representations
of the exceptional groups, getting that only $H_4$ admits
such a representation (actually two of them) with
$\eps\otimes \eta \into  S^2\rho^*$ for some $\eta \in \XX(\rho)$.

\newpage
\tableofcontents

\newpage

\vfill
\eject


\begin{thebibliography}{Ma06b}

\bibitem[Bd]{benard}  M. Benard, {\it Schur indices and splitting fields of the
Unitary reflection groups} J. Algebra {\bf 38} (1976), 318--342.
\bibitem[Bes]{bessis} D. Bessis,
{\it Sur le corps de d{\'e}finition d'un groupe de r{\'e}flexions
complexe}, Comm. Algebra {\bf 25} (1997), 2703--2716.
\bibitem[Bor]{BOREL} A. Borel, {\it Linear algebraic groups}, 2nd edition, GTM 126, Springer 1991.
\bibitem[Ar]{ARIKI} S. Ariki, {\it Representation theory of a Hecke algebra for $G(r,p,n)$}, {\it J. Algebra \bf 177} (1995), 164--185.
\bibitem[ArKo]{ARIKIKOIKE} S. Ariki, K. Koike, {\it A Hecke algebra of $\Z/r\Z \wr  \mathfrak{S}_n$ and construction of its irreducible representations},
{\it Adv. Math. \bf 106} (1994), 216--243.
\bibitem[Bo456]{BOURB456} N. Bourbaki, {\it Groupes et alg\`ebres de Lie, chapitres 4,5,6}, Hermann, 1968.
\bibitem[Br]{BRIESKBOURB} E. Brieskorn, {\it Sur les groupes de tresses [d'apr\`es
V.I. Arnold]}, S\'eminaire Bourbaki 1971/1972, expos\'e 401, LNM 317, Springer-Verlag, 1973.
\bibitem[BMR]{BMR} M. Brou\'e, G. Malle, R. Rouquier, {\it Complex reflection groups, braid groups, Hecke algebras}, J. Reine Angew. Math. {\bf 500} (1998) 127-190.
\bibitem[FLW]{FLW} M. Freedman, M. Larsen, Z. Wang, {\it The two-eigenvalue problem and density of Jones representation
of braid groups}, Comm. Math. Phys. {\bf 228} (2002), 177--199.
\bibitem[FH]{FH}  W. Fulton, J. Harris,
{\it Representation theory. A first course}, GTM 129, Springer, Berlin-Heidelberg-New York, 1991.
\bibitem[GP]{GP} M. Geck, G. Pfeiffer, {\it Characters of Finite Coxeter Groups and Iwahori-Hecke algebras}, London Math. Soc. Monographs new series 21, Clarendon Press, 2000.
\bibitem[Ko]{KOHNOKOSZUL} T. Kohno, {\it S\'erie
de Poincar\'e-Koszul associ\'ee aux groupes de tresses pures}, Invent. Math {\bf 82} (1985), 57--76. 
\bibitem[Ma01]{THESE} I. Marin, {\it Re\-pr\'e\-sen\-ta\-tions lin\'eaires des tresses infinit\'esimales}, Th\`ese de l'universit\'e Paris XI-Orsay, 2001.
\bibitem[Ma03a]{QUOT} I. Marin, {\it Quotients infinit\'esimaux du groupe de tresses}, Ann. Inst. Fourier (Grenoble) {\bf 53} (2003), 1323--1364.
\bibitem[Ma03b]{HECKINF} I. Marin, {\it Infinitesimal Hecke Algebras}, Comptes Rendus Math\'ematiques {\bf 337} S\'erie I, 297--302 (2003).
\bibitem[Ma04]{IRRED} I. Marin, {\it Irr\'eductibilit\'e g\'en\'erique des 
produits tensoriels de monodromies}, Bull. Soc. Math. Fr. {\bf 132}, 201--232 (2004).  
\bibitem[Ma05]{ASSOC} I. Marin, {\it On the representation theory of braid groups}, preprint arxiv:RT/0502118 (v3), 2005.
\bibitem[Ma06]{GT} I. Marin, {\it Caract\`eres de rigidit\'e du groupe
de Grothendieck-Teichm\"uller}, Compos. Math. {\bf 142} (2006), 657--678. 
\bibitem[Ma07a]{LIETRANSP} I. Marin, {\it L'alg\`ebre de Lie des transpositions}, J. Algebra {\bf 310} (2007), 742--774.
\bibitem[Ma07b]{KRAMMINF} I. Marin, {\it Sur les repr\'esentations de Krammer g\'en\'eriques}, Ann. Inst. Fourier (Grenoble) {\bf 57}, 1883--1925 (2007).
\bibitem[Ma07c]{BRANCH} I. Marin, {\it Branching properties for the groups $G(de,e,r)$}, preprint arxiv:0711.1845 (v2).
\bibitem[Ma07d]{KRAMCRG} I. Marin, {\it Krammer representations for complex braid groups}, preprint arxiv:math/0711.3096.
\bibitem[Ma08]{galglie} I. Marin, {\it Group algebras of finite groups as Lie algebras}, preprint arxiv:math/0809.0074, to appear in
Comm. Algebra.
\bibitem[Ma09]{arrrefl} I. Marin, {\it Reflection groups acting on their
hyperplanes}, J. Algebra {\bf 322} (2009), 2848--2860.
\bibitem[MM]{MARINMICHEL} I. Marin, J. Michel, {\it Automorphisms of
complex reflection groups}, preprint arxiv:math/0701266.
\bibitem[Sta]{STANLEY} R. Stanley, {\it Relative invariants of finite groups generated by pseudoreflections},
J. Algebra {\bf 49} (1977), 134--148. 
\bibitem[Ste]{STEINBERG} R. Steinberg, {\it Differential equations
invariant under finite reflection groups}, Trans. A.M.S. {\bf 112} (1964) 392--400.
\bibitem[Sp]{SPRINGER} T. Springer, {\it Linear algebraic groups, 2nd edition}, Birkh\"auser, 1998.
\bibitem[Wz]{WENZL} H. Wenzl, {\it Hecke algebras of type $A_n$ and subfactors}, Invent. Math. {\bf 92} (1988), 349--383.
\end{thebibliography}
\end{document}